%% file: sbim-cm.tex
\pdfoutput=1
\RequirePackage{silence}
\documentclass[a4paper,svgnames]{amsart}
\usepackage[hmarginratio={1:1},vmarginratio={1:1},lmargin=60.0pt,tmargin=60.0pt]{geometry}

\usepackage{booktabs}
\allowdisplaybreaks

\usepackage[T1]{fontenc}
\usepackage{latexsym,exscale,amsfonts,amssymb,mathtools}
\usepackage{amsmath,amsthm,amsfonts,amssymb,amscd,textcomp,bbm}
\usepackage{mathrsfs,stackrel}
\usepackage{etoolbox}
\usepackage{tikz,tikz-cd}
\usetikzlibrary{matrix,arrows.meta,decorations.markings,shapes.multipart,decorations.pathreplacing,decorations.pathmorphing,shapes.geometric}
\usepackage{aliascnt}
\usepackage{ytableau}
\usepackage{xparse}
\usepackage{cite}

\usepackage{dynkin-diagrams}
\usepackage{nicematrix}
\usepackage{tikz-3dplot}
\usepackage{mathrsfs}
\DeclareMathAlphabet{\mathscrbf}{OMS}{mdugm}{b}{n}

\usepackage{array}
\newcolumntype{C}{>{$}c<{$}}

\usepackage{xcolor,colortbl}

\definecolor{mygray}{gray}{0.6}
\definecolor{mygraydark}{gray}{0.4}
\definecolor{mygraylight}{gray}{0.85}
\definecolor{spinach}{RGB}{46,139,87}
\definecolor{tomato}{RGB}{255,99,71}
\definecolor{orchid}{RGB}{143,40,194}
\definecolor{neon}{RGB}{77,77,255}
\definecolor{lightneon}{RGB}{110,110,255}
\definecolor{pumpkin}{RGB}{224,180,80}
\definecolor{citron}{RGB}{190,180,90}

\definecolor{lava}{RGB}{207,16,32}
\definecolor{cream}{RGB}{255,253,208}
\definecolor{verdigris}{RGB}{67,179,174}
\definecolor{Black}{RGB}{0,0,0}
\definecolor{mydarkblue}{RGB}{10,10,170}
\definecolor{darkspinach}{RGB}{20,70,20}
\definecolor{darktomato}{RGB}{155,40,30}
\definecolor{darkorchid}{RGB}{50,10,100}
\definecolor{darklava}{RGB}{150,8,16}

\definecolor{zero}{RGB}{0,0,0}
\definecolor{one}{RGB}{255,0,0}
\definecolor{two}{RGB}{0,255,0}
\definecolor{three}{RGB}{0,0,255}

\usepackage{todonotes}

\usepackage{enumitem}
\setlist[enumerate]{itemsep=0.15cm,label=\emph{\upshape(\alph*)}}
\setlist[enumerate,2]{itemsep=0.15cm,label=\emph{\upshape(\roman*)}}
\setlist[enumerate,3]{itemsep=0.15cm,label=\emph{\upshape(\Alph*)}}

\let\emph\relax
\DeclareTextFontCommand{\emph}{\bfseries\em}

\renewcommand{\dots}{\text{...}}
\renewcommand{\vdots}{\rotatebox{90}{\text{...}}}
\renewcommand{\ddots}{\raisebox{0.175cm}{\rotatebox{-45}{\text{...}}}}

\newcommand{\placeholder}{{}_{-}}
\newcommand{\mystrut}{\rule[-0.2\baselineskip]{0pt}{0.9\baselineskip}}

\newcommand{\cf}{\text{cf.}}
\newcommand{\etc}{\text{etc.}}

\newcommand{\muta}{\text{mutatis mutandis}}

\newcommand{\C}{\mathbb{C}}
\newcommand{\R}{\mathbb{R}}

\newcommand{\N}{\mathbb{Z}_{\geq 0}}

\newcommand{\Z}{\mathbb{Z}}
\newcommand{\Zv}{\Z^{[\vpar]}}
\newcommand{\Nv}{\N^{[\vpar]}}
\newcommand{\K}{\mathbb{K}}

\newcommand{\algstuff}[1]{\mathrm{#1}}
\newcommand{\vanideal}[1]{\algstuff{J}_{#1}}
\newcommand{\vanset}[1]{\algstuff{V}_{#1}}
\newcommand{\vanpar}[1]{\algstuff{V}'_{#1}}

\newcommand{\munit}{\mathbbm{1}}

\newcommand{\setstuff}[1]{\mathrm{#1}}
\newcommand{\catstuff}[1]{\mathbf{#1}}
\newcommand{\varstuff}[1]{\mathtt{#1}}

\newcommand{\obstuff}[1]{\mathtt{#1}}
\newcommand{\morstuff}[1]{\mathrm{#1}}

\newcommand{\idmor}{\morstuff{id}}

\newcommand{\End}{\setstuff{End}}
\newcommand{\Hom}{\setstuff{Hom}}

\newcommand{\rank}{N}
\newcommand{\levelthing}{M}
\newcommand{\level}{e}
\newcommand{\sym}{\mathfrak{S}}

\newcommand{\sln}[1][\rank]{\mathfrak{sl}_{#1}}

\newcommand{\crg}[2]{G(#1,#1,#2)}
\newcommand{\crgg}[2]{G(#1,1,#2)}

\newcommand{\rooty}{\eta}
\newcommand{\rankroot}{\zeta}
\newcommand{\verl}[2]{\catstuff{V}(#2)^{#1}}

\newcommand{\verg}[2]{\catstuff{a}_{#1,#1,#2}}
\newcommand{\verG}[2]{\catstuff{A}_{#1,#1,#2}}

\newcommand{\dcenter}[1]{\catstuff{Z}(#1)}
\newcommand{\ggroup}[1]{[#1]_{\oplus}}
\newcommand{\ggroupc}[1]{[#1]_{\oplus}^{\C}}
\newcommand{\simples}[1]{\setstuff{Si}(#1)}

\newcommand{\vpar}{\varstuff{v}}
\newcommand{\qpar}{\varstuff{q}}
\newcommand{\Cq}{\C^{\qpar}}
\newcommand{\Cv}{\C^{\vpar}}
\newcommand{\vnum}[1]{[#1]_{\vpar}}

\newcommand{\SUN}{\mathrm{SU}_{\rank}}
\newcommand{\Usln}{U_{\qpar}(\mathfrak{sl}_{\rank})}

\newcommand{\Ursln}{U_{\rooty}(\mathfrak{sl}_{\rank})}
\newcommand{\slncat}{\catstuff{Rep}_{\qpar}(\sln)}
\newcommand{\slncatr}{\catstuff{Rep}_{\rooty}(\sln)}
\newcommand{\Ll}{\obstuff{L}}

\newcommand{\ba}{\mathbf{a}}

\newcommand{\bk}{\mathbf{k}}
\newcommand{\bm}{\mathbf{m}}
\newcommand{\bg}{\boldsymbol{\gamma}}
\newcommand{\bl}{\boldsymbol{\lambda}}
\newcommand{\bsum}{\Sigma}
\newcommand{\bs}{\boldsymbol{\sigma}}
\newcommand{\transpose}[1]{\overline{#1}}

\newcommand{\iunit}{\mathtt{i}}
\newcommand{\junit}{\Upsilon}

\newcommand{\fu}{\varstuff{X}}
\newcommand{\pxy}[1]{\setstuff{U}_{#1}}
\newcommand{\Zfun}{Z}
\newcommand{\Efun}{E^{-}}

\newcommand{\tnumbers}[2]{p_{#1}^{#2}}

\newcommand{\br}{\beta}
\newcommand{\rev}{\mathrm{rev}}
\newcommand{\Tr}{\mathrm{Tr}}

\newcommand{\family}{\mathcal{F}}

\newcommand{\frob}{\mathrm{Frob}}
\newcommand{\prefourier}{\tilde{\mathbb{S}}}
\newcommand{\cmcell}{\Gamma}

\newcommand{\turnsymbol}{\to}
\newcommand{\stabilizer}[1]{\mathrm{stab}_{#1}}
\newcommand{\Stabilizer}[1]{\mathrm{Stab}_{#1}}

\newcommand{\regular}{R}
\newcommand{\idemp}{e}
\newcommand{\Idemp}{E}

\newcommand{\nodes}{I}
\newcommand{\weyl}{W}
\newcommand{\hecke}{H(\tilde{A}_{\rank{-}1})}
\newcommand{\nhecke}[1][\infty]{T_{#1}}

\newcommand{\kl}{\theta}
\newcommand{\kltwo}{\Theta}
\newcommand{\rkl}[2]{h^{#1}_{#2}}
\newcommand{\lkl}[2]{{}^{#1}_{#2}h}
\newcommand{\RKL}[2]{c^{#1}_{#2}}
\newcommand{\LKL}[2]{{}^{#1}_{#2}c}
\newcommand{\scalar}{\varkappa}

\newcommand{\basisH}[1][\infty]{\algstuff{H}\hspace*{-.015cm}{}^{#1}}
\newcommand{\Hbasis}[1][\infty]{{}^{#1}\algstuff{H}}
\newcommand{\basisC}[1][\infty]{\algstuff{C}\hspace*{-.015cm}{}^{#1}}
\newcommand{\Cbasis}[1][\infty]{{}^{#1}\algstuff{C}}

\newcommand{\lcell}{\mathsf{L}}
\newcommand{\rcell}{\mathsf{R}}
\newcommand{\tcell}{\mathsf{J}}

\newcommand{\killideal}[1]{\algstuff{I}_{#1}}

\newcommand{\M}{\algstuff{M}}
\newcommand{\matrep}[2]{\M_{#1}(#2)}

\newcommand{\vertices}{\mathtt{V}}
\newcommand{\edges}{\mathtt{E}}

\newcommand{\point}{\catstuff{P}}
\newcommand{\alg}{\algstuff{A}^{D}}

\usepackage{xltabular}
\usepackage{tabulary}
\usepackage{booktabs}

\usepackage[all]{xy}
\usepackage{tikz}
\usetikzlibrary{cd}
\usetikzlibrary{decorations}
\usetikzlibrary{decorations.markings}
\usetikzlibrary{decorations.pathreplacing}
\usetikzlibrary{decorations.pathmorphing}
\usetikzlibrary{arrows.meta,shapes,positioning,matrix,calc}
\usetikzlibrary{shapes.callouts}
\tikzset{anchorbase/.style={baseline={([yshift=-0.5ex]current bounding box.center)}},
tinynodes/.style={font=\tiny,text height=0.25ex,text depth=0.05ex},
smallnodes/.style={font=\scriptsize,text height=0.75ex,text depth=0.15ex},
}

\tikzset{->-/.style={decoration={
markings,
mark=at position .5 with {\arrow{>}}},
postaction={decorate}}
}

\tikzset{-<-/.style={decoration={
markings,
mark=at position .5 with {\arrow{<}}},
postaction={decorate}}
}

\tikzset{-><-/.style={decoration={
markings,
mark=at position .4 with {\arrow{<}},
mark=at position .6 with {\arrow{>}}},
postaction={decorate}}
}

\tikzset{square/.style={regular polygon,regular polygon sides=4}}
\tikzset{triangle/.style={regular polygon,regular polygon sides=3}}
\tikzset{pentagon/.style={regular polygon,regular polygon sides=5}}

\def\NewTheorem#1{%
\newaliascnt{#1}{equation}%
\newtheorem{#1}[#1]{#1}%
\aliascntresetthe{#1}%
\expandafter\def\csname #1autorefname\endcsname{#1}%
}
\def\equationautorefname~#1\null{(#1)\null}

\numberwithin{equation}{subsection}

\NewTheorem{Proposition}
\NewTheorem{Theorem}
\NewTheorem{Corollary}
\AtEndEnvironment{Corollary}{\null\hfill$\square$}%
\NewTheorem{Lemma}
\theoremstyle{definition}
\NewTheorem{Definition}
\AtEndEnvironment{Definition}{\null\hfill$\Diamond$}%
\NewTheorem{Classification Problem}
\AtEndEnvironment{Classification Problem}{\null\hfill$\Diamond$}%
\NewTheorem{Notation}
\AtEndEnvironment{Notation}{\null\hfill$\Diamond$}%
\NewTheorem{Example}
\AtEndEnvironment{Example}{\null\hfill$\Diamond$}%
\NewTheorem{Examples}
\AtEndEnvironment{Examples}{\vskip-10mm\null\hfill$\Diamond$}%
\NewTheorem{Conjecture}

\theoremstyle{remark}
\NewTheorem{Remark}
\AtEndEnvironment{Remark}{\null\hfill$\Diamond$}%
\NewTheorem{Question}
\AtEndEnvironment{Question}{\null\hfill$\Diamond$}%
\NewTheorem{Assumption}

\usepackage[hypertexnames=false]{hyperref}
\usepackage{bookmark}
\hypersetup{
pdftoolbar=true,
pdfmenubar=true,
pdffitwindow=false,
pdfstartview={FitH},
pdftitle={On Hecke and asymptotic categories for complex reflection groups},
pdfauthor={Abel Lacabanne, Daniel Tubbenhauer and Pedro Vaz},
pdfsubject={},
pdfcreator={Abel Lacabanne, Daniel Tubbenhauer and Pedro Vaz},
pdfproducer={Abel Lacabanne, Daniel Tubbenhauer and Pedro Vaz},
pdfkeywords={},
pdfnewwindow=true,
colorlinks=true,
linkcolor=mydarkblue,
citecolor=teal,
filecolor=magenta,
urlcolor=orchid,
linkbordercolor=lava,
citebordercolor=teal,
urlbordercolor=orchid,
linktocpage=true
}

\def\makeautorefname#1#2{\csdef{#1autorefname}{#2}}
\makeautorefname{section}{Section}%
\makeautorefname{subsection}{Section}%
\makeautorefname{subsubsection}{Section}%

\def\changed#1{#1}
\def\ochanged#1{#1}
\def\changedd#1{#1}
\def\ochangedd#1{#1}
\def\changeddd#1{{#1}}
\def\ochangeddd#1{{#1}}

\definecolor{darkpurple}{HTML}{2E1A47}
\definecolor{lavender}{HTML}{9B80B5}
\definecolor{palelavender}{HTML}{D9B7E2}

\begin{document}
\title[On Hecke and asymptotic categories for \changedd{a family of} complex reflection groups]{On Hecke and asymptotic categories for \changedd{a family of} complex reflection groups}
\author[A. Lacabanne, D. Tubbenhauer and P. Vaz]{Abel Lacabanne, Daniel Tubbenhauer and Pedro Vaz}

\address{A.L.: Laboratoire de Math{\'e}matiques Blaise Pascal (UMR 6620), Universit{\'e} Clermont Auvergne, Campus Universitaire des C{\'e}zeaux, 3 place Vasarely, 63178 Aubi{\`e}re Cedex, France,\newline \href{http://www.normalesup.org/~lacabanne}{www.normalesup.org/$\sim$lacabanne},
\href{https://orcid.org/0000-0001-8691-3270}{ORCID 0000-0001-8691-3270}}
\email{abel.lacabanne@uca.fr}

\address{D.T.: The University of Sydney, School of Mathematics and Statistics F07, Office Carslaw 827, NSW 2006, Australia, \href{http://www.dtubbenhauer.com}{www.dtubbenhauer.com}, \href{https://orcid.org/0000-0001-7265-5047}{ORCID 0000-0001-7265-5047}}
\email{daniel.tubbenhauer@sydney.edu.au}

\address{P.V.: Institut de Recherche en Math{\'e}matique et Physique, 
Universit{\'e} catholique de Louvain, Chemin du Cyclotron 2,  
1348 Louvain-la-Neuve, Belgium, \href{https://perso.uclouvain.be/pedro.vaz}{https://perso.uclouvain.be/pedro.vaz}, \href{https://orcid.org/0000-0001-9422-4707}{ORCID 0000-0001-9422-4707}}
\email{pedro.vaz@uclouvain.be}

\begin{abstract}
Generalizing the dihedral picture for G(M,M,2), we construct Hecke algebras (and \changed{present a strategy for constructing Hecke} categories) and asymptotic counterparts. We think of these as associated with the complex reflection group G(M,M,N).
\end{abstract}


\addtocontents{toc}{\protect\setcounter{tocdepth}{1}}

\maketitle

\tableofcontents

\section{Introduction}

Weyl groups are among the most important objects in algebra, as they govern the representation theory of their associated reductive group. Weyl groups are real reflection groups and special cases of \emph{complex reflection groups}, and it is an interesting question what kind of `reductive group' is associated with a complex reflection group. These `reductive groups' are \changed{(probably) not reductive groups themselves, but are} believed to exist in a certain sense \changed{and share many combinatorial and representation theoretical features of reductive groups. T}hey were famously named by Brou{\'e}--Malle--Michel \cite{BrMaMi-spets} after the Greek island Spetses: such a `reductive group' is called a \emph{spets}.

We cannot provide a definitive answer to what spetses are, but recent developments indicate that Soergel bimodules, also known as \emph{Hecke categories}, can often serve as a replacement whenever an associated geometric or Lie theoretic picture is missing.

In this paper, we focus on the complex reflection groups of type $\crg{\levelthing}{\rank}$ for $\levelthing \geq \rank$ and suggest that they have an associated Hecke algebra and category. These arise from Chebyshev polynomials associated with root systems, have Kazhdan--Lusztig (KL) type combinatorics, include asymptotic categories, are related to Calogero--Moser (CM) families, and encode Fourier matrices for $\crg{\levelthing}{\rank}$.

\subsection{From dihedral group to the main results}\label{SS:IntroA}

The dihedral group $\crg{\levelthing}{2}$ of order $2\levelthing$ is 
one of the simplest examples of a complex reflection group that is not a Weyl group (unless $\levelthing$ is small), 
yet it still exhibits behavior typically associated with objects from Lie theory. 

Let us list a few of these, all of which are for the middle cell:
\begin{enumerate}[label=(\Alph*)]

\item The KL basis is determined by the coefficients of the Chebyshev polynomials \cite{duCl-positivity-finite-hecke,El-two-color-soergel,Tu-sandwich-cellular}, 
and the Chebyshev polynomials determine the characters of the simple representations of $\mathrm{SL}_{2}$.

\item The $\N$-representations of the dihedral group, or the $2$-representations of the dihedral Hecke category, are indexed by ADE Dynkin diagrams \cite{KiMaMaZi-2-reps-soergel,MaTu-soergel}. This 
is closely related to (but does not quite match) with three different objects: irreducible conformal field theories (CFT for short) for $\mathrm{SU}_{2}$ \cite{Pa-cft}, subgroups of quantum $\mathrm{SU}_{2}$ 
\cite{KiOs-ade-verlinde,Os-module-categories} and module categories of the $\mathrm{SL}_{2}$ Verlinde category \cite{Os-module-categories}.

\item The Drinfeld centers of the asymptotic categories associated to the dihedral Hecke algebras are modular categories whose $S$-matrices coincide  
with the so-called Fourier matrices \cite{RoTh-center}, the base change between `unipotent characters' and `unipotent character sheaves', constructed in the 90s \cite{Lu-exotic-fourier} in an ad hoc fashion. Hence, 
one can think of these centers as `unipotent character sheaves', matching 
\cite{Lu-truncated} which proposed `unipotent character sheaves' associated to Coxeter groups that encode the Fourier matrices.

\item The simple representations of the (complexified) asymptotic Hecke algebra are given by the KL family associated to the cell \cite{Lu-leading-coeff-hecke}.

\end{enumerate}
In this paper, we generalize all the above to the case $\rank\geq 2$ (for $\rank=1$ the group $\crg{\levelthing}{\rank}$ is trivial and we from now on assume that $\rank>1$).
The outline is as follows.
\begin{enumerate}[label=(\alph*)]

\item \textit{Generalizing (A).} We define a certain subalgebra $\nhecke$ of the Hecke algebra 
of affine type $A_{\rank-1}$, and then a finite dimensional quotient $\nhecke[\level]$, the \emph{Nhedral Hecke algebra} of level $\level$, obtained from $\nhecke$ by annihilating certain KL basis elements. 
\textit{This algebra is immediately truncated to the analog of the aforementioned middle cell.} As we will show, the algebra 
$\nhecke$ has its own KL theory and representation theory akin to the dihedral Hecke algebra. 
Its KL basis is determined by the coefficients of the Chebyshev polynomials associated with simple representations of $\mathrm{SL}_{\rank}$.

\item \textit{Generalizing (\ochangedd{B}).} Returning to $\nhecke[\level]$, we \changed{observe (some aspects proven and others conjectural)} that most of Zuber's generalized ADE Dynkin diagrams \cite{Zu-gen-dynkin-diagrams} give rise 
to $\N$-representations of $\nhecke[\level]$. This 
is closely related (but does not quite match) with three different objects: irreducible CFT for $\mathrm{SU}_{\rank}$ as e.g. in \cite{PeZu-from-CFT-to-graphs}, subgroups of quantum $\mathrm{SU}_{\rank}$ 
as e.g. in \cite{Oc-classification-sun} and module categories of the $\mathrm{SL}_{\rank}$ Verlinde category as, for example, in \cite{Ga-SU3}.

\item \textit{Generalizing (\ochangedd{C}).} \changed{(With a (*), see below.)} There \changed{should be} a categorification of $\nhecke[\level]$, called \emph{Nhedral Hecke category} (or Nhedral Soergel bimodules) of level $\level$. 
This category \changed{should be} positively graded. \changed{We define what we believe is its degree zero part and we call it the \emph{asymptotic category} $\verg{\levelthing}{\rank}$}. We show that 
the Drinfeld center of $\verg{\levelthing}{\rank}$ is a modular category and compute its $S$ and $T$ matrices. We show that the $S$ matrix coincides 
with Malle's Fourier matrix for $\crg{\levelthing}{\rank}$ \cite{Ma-unipotente-grade}.

\item \textit{Generalizing (\ochangedd{D}).} We then define a matrix category $\verG{\levelthing}{\rank}$ over $\verg{\levelthing}{\rank}$, the \emph{big asymptotic category}, which, by construction, 
is Morita equivalent to $\verg{\levelthing}{\rank}$. We explain how this category is related to a CM family for $\crg{\levelthing}{\rank}$\changed{, which plays the role of the middle cell for the dihedral group}.

\end{enumerate}
Let us summarize the main points of the paper with two overview diagrams. The first diagram illustrates how the main algebras and categories are related:
\begin{gather*}
\begin{tikzpicture}
\node (A) at (-2.6, 2.5) [draw, rectangle, minimum width=10em, align=center, neon,fill=neon!25] {\color{black}\shortstack{Nhedral Hecke\\category (*)}};
\node (B) at (4, 2.5) [draw, rectangle, minimum width=10em,align=center, tomato,fill=tomato!25] {\color{black}\shortstack{Asymptotic\\category $\verg{\levelthing}{\rank}$}};
\node (C) at (10.6, 2.5) [draw, rectangle, minimum width=10em,align=center,orchid,fill=orchid!25] {\color{black}\shortstack{Big asymptotic\\category $\verG{\levelthing}{\rank}$}};
\node (D) at (-2.6, 0) [draw, rectangle, minimum width=10em,align=center, neon,fill=neon!25] {\color{black}\shortstack{Nhedral Hecke\\algebra}};
\node (E) at (4, 0) [draw, rectangle, minimum width=10em,align=center, tomato,fill=tomato!25] {\color{black}\shortstack{Decategorification\\ of $\verg{\levelthing}{\rank}$}};
\node (F) at (10.6, 0) [draw, rectangle, minimum width=10em,align=center,orchid,fill=orchid!25] {\color{black}\shortstack{Asymptotic Hecke al-\\gebra for a CM cell}};
\draw[<->] (A) -- node[above] {degree zero} node[below] {equivalence \autoref{QSH}} (B);
\draw[<->] (B) -- node[above] {Morita} node[below] {equivalence} (C);
\draw[<->] (D) -- node[above] {degree zero} node[below] {isomorphism \autoref{QSH}} (E);
\draw[<->] (E) -- node[above] {Morita} node[below] {equivalence} (F);
\draw[dashed, ->, neon] (A) -- node[right] {decat.} (D);
\draw[dashed, ->, tomato] (B) -- node[right] {decat.} (E);
\draw[dashed, ->, orchid] (C) -- node[right] {decat.} (F);
\end{tikzpicture}
\end{gather*}
The ``\autoref{QSH}'' is \changed{explained} in \autoref{SS:Remarks} below, while (*) means that we believe this is definable and interesting to define, though we do not do this here since the relevant technology is missing while writing this paper.
\begin{gather*}
\scalebox{1.25}{
\begin{tikzpicture}
\shade[radius=0.85cm, top color=darkorchid, bottom color=lightneon, middle color=tomato, opacity=0.5] (0,0) circle (0.85cm);
\node (0) at (0, 0) [draw, circle, align=center, radius=0.9cm] {\color{black}\shortstack{Nhedral\\picture}};
\node (1) at ({60*0}:4.3cm) [draw, rectangle, minimum width=9em, minimum height=3em, align=center,fill=mygray!25] {\color{black}\shortstack{A. Chebyshev\\polynomials}};
\node (2) at ({45*1}:3.5cm) [draw, rectangle, minimum width=9em, minimum height=3em, align=center, orchid,fill=orchid!25] {\color{black}\shortstack{\ochangedd{F. CM family}}};
\node (3) at ({45*3}:3.5cm) [draw, rectangle, minimum width=9em, minimum height=3em, align=center, orchid,fill=orchid!25] {\color{black}\shortstack{\ochangedd{E. Unipotent chars.}\\\ochangedd{for \crg{\levelthing}{\rank}}}};
\node (4) at ({60*3}:4.3cm) [draw, rectangle, minimum width=9em, minimum height=3em, align=center, neon,fill=neon!25] {\color{black}\shortstack{D. Irreducible\\CFT}};
\node (5) at ({45*5}:3.5cm) [draw, rectangle, minimum width=9em, minimum height=3em, align=center, neon,fill=neon!25] {\color{black}\shortstack{\ochangedd{C}. Subgroups of\\quantum $\mathrm{SU}(\rank)$}};
\node (6) at ({45*7}:3.5cm) [draw, rectangle, minimum width=9em, minimum height=3em, align=center, neon,fill=neon!25] {\color{black}\shortstack{\ochangedd{B}. Module cats for\\Verlinde cats}};
\draw[->] (0) -- node[above] {\autoref{SS:CPoly}} (1);
\draw[->, orchid] (0) -- node[right] {\autoref{SS:BigCat}} (2);
\draw[->, orchid] (0) -- node[left] {\autoref{SS:FamilyUnipotent}\,} (3);
\draw[->, neon] (0) -- node[above] {\autoref{SS:NhedralIntegral}} (4);
\draw[->, neon] (0) -- node[left] {\autoref{SS:NhedralIntegral}} (5);
\draw[->, neon] (0) -- node[right] {\autoref{SS:NhedralIntegral}} (6);
\end{tikzpicture}
}
\end{gather*}
In this picture we indicate \changedd{the broader framework of Nhedral combinatorics. Though e.g. irreducible CFTs are beyond the scope of this paper and we do not discuss them any further, \autoref{S:Nhedral} contains the combinatorial material which relates to this topic. This picture is not intended as a rigorous statement but rather as an illustration of the numerous (potential) connections.}

\subsection{Related works and speculations}

\changed{The only infinite family of
complex reflection groups is $G(M,D,N)$ with $D\mid M$.
Among these, there are two infinite families of so-called spetsial complex reflection groups, which represent the two different extremes $D=M$ and $D=1$, yet both might permit some Hecke category combinatorics.} Additionally, there is a twisted variant that appears to fit into the Hecke category framework, and we comment on all three here.

\smallskip
\noindent{\color{mygray}\rule[1.0ex]{\linewidth}{0.5pt}}

The first infinite family is $\crg{\levelthing}{\rank}$, as discussed in this paper, \changed{ and we interpret the Nhedral picture (Chebyshev polynomials {\etc}) as the combinatorics of the associated Hecke category.}
Furthermore, for $\rank=3$, \ochangedd{(A) and (B)} above were originally generalized in \cite{MaMaMiTu-trihedral}. The paper \cite{MaMaMiTu-trihedral} also provides a classification of $\N$-representations of $\nhecke[\level]$ for small $\level$, a task we expect to be achievable for at least small $\rank$ as well.
Moreover, exciting so-called exotic nilCoxeter and deformed affine NilHecke algebras are obtained in \cite{ElJuYo-nilcoxeter,ElJuYo-nilcoxeter2} (with some results also for $\rank>3$), but we do not know how to relate this work to ours.

\smallskip
\noindent{\color{mygray}\rule[1.0ex]{\linewidth}{0.5pt}}

The other infinite family of spetsial complex reflection groups consists of the groups $\crgg{\levelthing}{\rank}$, where the corresponding Hecke algebra is the Ariki--Koike algebra \cite{ArKo-hecke-algebra,BrMa-hecke}. Fourier matrices and families of unipotent characters for these groups have been studied quite extensively; see, for example, \cite{Ma-unipotente-grade,Cu-fusion-algebras,BoRo-cyclic,La-drinfeld-double,La-fourier-quantum-gln}.
However, we are not aware of any general Hecke category combinatorics for $\crgg{\levelthing}{\rank}$. The paper \cite{LaTuVa-annular-webs-levi} identifies the Ariki--Koike algebra as a subalgebra of webs on an annulus, suggesting the possibility of using annular foams, as in \cite{RoWa-sym-homology}, to describe a Hecke category. For $\rank=1$, a Hecke category has been proposed in \cite{GoTh-soergel-cyclic}, but it seems quite different from the annular web picture.

\smallskip
\noindent{\color{mygray}\rule[1.0ex]{\linewidth}{0.5pt}}

There is also the story of Fourier matrices associated with groups having an automorphism, such as the Ree groups \changed{or $\crgg{\levelthing}{\rank}$}; see, for example, \cite{Ma-unipotente-grade,GeMa-fourier-frobenius,La-crossed-S-matrix}.
The corresponding Fourier matrices are not symmetric and often not integral, so they may only arise via equivariant module categories and not via modular categories.
The automorphism `twists' the setting in a certain sense and one might hope that 
the `twisted' 2-representations of the Hecke category as e.g. \cite{MaMaMiTuZh-soergel-2reps} 
can be used to define degree zero module categories \changed{where} crossed $S$-matrices in the sense of \cite{De-modular-crossed} might play a role.

\subsection{Additional remarks}\label{SS:Remarks}

For some proofs \changed{and some proposed constructions} in this paper we need the following.
\begin{center}
\fcolorbox{tomato!100}{cream!75}{%
\begin{minipage}{0.8\linewidth}
\centering
The \emph{quantum Satake hypothesis (QSH)} is:
\begin{gather}\label{QSH}\tag{QSH}
\text{\cite[Theorem 5.35]{El-q-satake} is true.}
\end{gather}
Whenever we assume \autoref{QSH} holds, we state this explicitly.
\end{minipage}%
}
\end{center}

\begin{Remark}
The paper \cite{El-q-satake} proves the quantum Satake hypothesis for $\rank=2,3$\ochanged{. So we have no further assumptions in these two cases.}
\end{Remark}

\begin{Remark}\label{R:GeoSatake}
\changedd{We expect that several statements to follow (such as \autoref{L:Embedding}, \autoref{R:h-element}, \autoref{R:KL-basis-Chebyshev}, \autoref{L:KLaction}, \autoref{P:KLbasisInfinity}, or \autoref{L:integrality}) will not depend on \autoref{QSH} itself, except perhaps \autoref{L:integrality} in the finite case; see \autoref{R:GeoSatake2} for some more comments. On the other hand, we anticipate that \autoref{QSH} will play a role in the categorified context.}
\end{Remark}

\begin{Remark}
We postponed several proofs to 
\autoref{app}. If the reader is missing a proof, then they should be able to find it there. Moreover, in \autoref{S:slnthings} and \autoref{S:Nhedral}, we will generalize some of the main results of \cite{MaMaMiTu-trihedral}, following their exposition. Some results will have proofs that work {\muta} and we will be brief with these, and we will point out when the arguments are sufficiently different.
\end{Remark}

\begin{Remark}
The paper is readable in black-and-white but we recommend reading it in color.
\end{Remark}

\begin{Remark}
Code for some 
of the calculations in this paper is available on GitHub\ochangedd{, see} \cite{LaTuVa-Nhedral-code}.
\end{Remark}

\subsection{Table of notation and general conventions}
\label{SS:notations}

The following is the list of the most important concepts.

\renewcommand\arraystretch{1.15}
\begin{xltabular}{\textwidth}{p{1.95cm}>{\ttfamily}p{5.5cm}p{8.2cm}}\toprule
Symbol & Name & Description
\\
\midrule\endhead
\bottomrule\endfoot
$\placeholder$
&
Placeholder
&
Used as a placeholder symbol
\\
\hline\mystrut
$\ggroup{\placeholder}$
&
Grothendieck group
&
The additive Grothendieck group of $\placeholder$ (for semisimple categories this equals the abelian one)
\\
\hline\mystrut
$\placeholder{}^{\turnsymbol}$
&
A cyclic action
&
See \autoref{SS:QuantumGroup}.
\\
\hline\mystrut
$\transpose{\placeholder}$
&
The dual
&
If $\bm=(m_{1},\dots,m_{\rank{-}1})$, then
$\transpose{\bm}=(m_{\rank{-}1},\dots,m_{1})$
\\
\hline\mystrut
$(\placeholder)_{i}$
&
Degree $i$ part
&
If the category $\placeholder$ is graded, then this denotes the degree $i$ part
\\
\hline\mystrut
$A(\placeholder)$
&
Adjacency matrix
&
The adjacency matrix of a graph $\placeholder$
\\
\hline\mystrut
$\verG{\levelthing}{\rank}$
&
The big asymptotic category
&
The big asymptotic category for $\crg{\levelthing}{\rank}$ defined in \autoref{D:BigCat}
\\
\hline\mystrut
$\verg{\levelthing}{\rank}$
&
The asymptotic category
&
The asymptotic category for $\crg{\levelthing}{\rank}$ defined in \autoref{D:AsymptoticCat}
\\
\hline\mystrut
$\RKL{x}{y},\LKL{x}{y}$
&
Nhedral KL elements
&
See \autoref{SS:NHeckeAlgebras}
\\
\hline\mystrut
$d^{\bk}_{\bm}$
&
Change-of-basis coefficients
&
See \autoref{Eq:dNumbers}
\\
\hline\mystrut
$\dim_{c}\placeholder$
&
Categorical dimension
&
The categorical dimension of $\placeholder$
\\
\hline\mystrut
$\dim_{\rooty}\placeholder$
&
Quantum dimension
&
The quantum dimension of $\placeholder$
\\
\hline\mystrut
$\dim_{\K}\placeholder$
&
Usual dimension
&
The dimension of $\placeholder$ over $\K$
\\
\hline\mystrut
$\level$
&
Level
&
The level, a number $\level\in\N$ that we fix
\\
\hline\mystrut
$\rooty$
&
Root of unity
&
The primitive $2\levelthing$th root of unity $\exp(\iunit\pi/\levelthing)$
\\
\hline\mystrut
$\rooty^{1/\rank}$
&
Root of unity
&
The $\rank$th root of $\rooty$ given by $\exp(\iunit\pi/\rank\levelthing)$
\\
\hline\mystrut
$\family$
&
Family of unipotent characters
&
A family of unipotent characters for $\crg{\levelthing}{\rank}$, see \autoref{SS:FamilyUnipotent}
\\
\hline\mystrut
$\family_{0}$
&
Principal series
&
The principal series of $\family$
\\
\hline\mystrut
$\crg{\levelthing}{\rank}$
&
Complex reflection group
&
The imprimitive complex reflection group of order $\levelthing^{\rank{-}1}\rank!$
\\
\hline\mystrut
$\cmcell$
&
CM cell
&
The two-sided CM cell associated with $\family$, see 
\autoref{SS:BigCat}
\\
\hline\mystrut
$\rkl{x}{y},\lkl{x}{y}$
&
Nhedral Bott--Samelson elements
&
See \autoref{SS:NHeckeAlgebras}
\\
\hline\mystrut
$\iunit$
&
Imaginary unit
&
The usual square root of $-1$
\\
\hline\mystrut
$\nodes$
&
Vertices
&
$\nodes=\{0,\dots,\rank-1\}$, the vertices of the affine 
type $A_{\rank-1}$ Dynkin diagram
\\
\hline\mystrut
$\vanideal{\level}$
&
Vanishing ideal
&
Defined in \autoref{D:KoornwinderV}
\\
\hline\mystrut
$\levelthing$
&
Level (plus Coxeter number)
&
This is $\level+\rank$, the `level' of $\crg{\levelthing}{\rank}$ such that
$\crg{\levelthing}{2}$ is the dihedral group of order $2\levelthing$
\\
\hline\mystrut
$\M$
&
Nhedral representation
&
Representation of the Nhedral Hecke algebra; potentially decorated with symbols
\\
\hline\mystrut
$\rank$
&
Rank (potentially plus one)
&
The rank, a number $\rank\in\Z_{\geq 2}$ that we fix
\\
\hline\mystrut
$\scalar$
&
Quantum $\rank-1$ factorial
&
$\vnum{\rank-1}!$
\\
\hline\mystrut
$\tnumbers{\rank}{\level}$
&
Simplicial polytopic numbers
&
$\tnumbers{\rank}{\level}=\binom{\level+\rank{-}1}{\rank{-}1}=\binom{\levelthing-1}{\rank{-}1}$
\\
\hline\mystrut
$\qpar$
&
Quantum parameter
&
The quantum generic parameter
\\
\hline\mystrut
$\vnum{k}$
&
Quantum numbers
&
The $k$th quantum number, $\vnum{k}=\frac{\qpar^{k}-\qpar^{-1}}{\qpar-\qpar^{-1}}$
\\
\hline\mystrut
$\vnum{k}!$
&
Quantum factorials
&
The $k$th quantum factorial, $\vnum{k}!=\vnum{1}\dots\vnum{k}$
\\
\hline\mystrut
$\mathrm{rk}\,\placeholder$
&
Rank
&
The rank of the category $\placeholder$ (the number of indecomposable objects)
\\
\hline\mystrut
$\slncatr$
&
Representation category
&
The fusion category of $\Ursln$-representations
\\
\hline\mystrut
$\slncat$
&
Representation category
&
The category of (type $1$) $\Usln$-representations
\\
\hline\mystrut
$S$
&
An $S$-matrix
&
The matrix $S$ of a modular category involved in the action of the modular group
\\
\hline\mystrut
$\overline{S}$
&
Complex conjugate
&
Entry-wise complex conjugate of $S$
\\
\hline\mystrut
$\simples{\placeholder}$
&
Simple objects
&
The set of isomorphism classes of simple objects of $\placeholder$
\\
\hline\mystrut
$\bsum\bm$
&
Sum of the entries
&
For $\bm=(m_{1},\dots,m_{r})$ we let
$\bsum\bm=m_{1}+\dots+m_{r}$
\\
\hline\mystrut
$\stabilizer{\placeholder}$
&
Stabilizer
&
Cardinal of a stabilizer
\\
\hline\mystrut
$T$
&
A $T$-matrix
&
The matrix $T$ of a modular category involved in the action of the modular group
\\
\hline\mystrut
$\nhecke[\placeholder]$
&
Nhedral Hecke algebra
&
The Nhedral Hecke algebra of level $\placeholder$, see 
e.g. \autoref{D:Nhedral}
\\
\hline\mystrut
$\theta$
&
Ribbon structure
&
The ribbon structure on $\slncatr$ defined via $\rooty^{1/\rank}$
\\
\hline\mystrut
$\kl_{i}$
&
Nhedral generators
&
The generators of $\nhecke[\placeholder]$
\\
\hline\mystrut
$\pxy{\bm}$
&
Chebyshev polynomials for $\sln[\rank]$
&
The higher versions of Chebyshev polynomials, see \autoref{D:CKpoly}, due to Koornwinder and Eier--Lidl
\\
\hline\mystrut
$\Usln$
&
Quantum $\sln$
&
The generic quantum group for quantum parameter $\qpar$
\\
\hline\mystrut
$\Ursln$
&
Another quantum $\sln$
&
The quantum group specialized quantum parameter $\rooty$
\\
\hline\mystrut
$\vanset{\level}$
&
Koornwinder variety
&
Defined in \autoref{D:KoornwinderV}
\\
\hline\mystrut
$\vanpar{\level}$
&
Reparametrization of $\vanset{\level}$
&
Defined in \autoref{D:KoornwinderPrime}
\\
\hline\mystrut
$\verl{\level}{\rank}$
&
Verlinde category
&
The Verlinde category for $\sln$ defined using $\rooty$, with ribbon structure $\theta$
\\
\hline\mystrut
$\weyl$
&
Affine Weyl group
&
The Weyl group for affine type $A_{\rank-1}$
\\
\hline\mystrut
$\weyl_{\changed{\mathtt{i}}}$
&
Parabolic subgroup
&
The parabolic subgroup for the vertices \changed{$\nodes\setminus\{\mathtt{i}\}$}
\\
\hline\mystrut
$\chi_{c}$
&
Central character
&
The ``colors'' associated with simple $\Usln$-representations
\\
\hline\mystrut
\changedd{$X^{+}$}
&
\changedd{Dominant weights}
&
\changedd{The set of dominant weights of type $A_{n-1}$}
\\
\hline\mystrut
\changedd{$X^{+}(\level)$}
&
\changedd{Level $e$ cut-off}
&
\changedd{The level $e$ cut-off of dominant weights of type $A_{n-1}$}
\\
\hline\mystrut
$\fu_{i}$
&
Fundamental variables
&
The variables used for the fundamental $\Usln$-representations
\\
\hline\mystrut
$\dcenter{\placeholder}$
&
Drinfeld center
&
The Drinfeld center of a category $\placeholder$
\\
\hline\mystrut
$\Zfun_{i}$
&
Koornwinder's Z-functions
&
Defined in \autoref{D:Zfunctions}
\\
\hline\mystrut
\changed{$\rankroot$}
&
\changed{Root of unity}
&
\changed{The primitive $\rank$th root of unity $\exp(2\iunit\pi/\rank)$}
\end{xltabular}
\label{table:notation}

\changeddd{Unless otherwise specified, the following general conventions apply throughout:}
\begin{enumerate}

\item Our conventions for rings and fields, \changedd{where we use the generic symbol $\K$}, are: \changeddd{superscripts} mean we adjoin a certain element to the ring/field, e.g. $\Cq=\C(\qpar)$. We write square brackets if we adjoin these elements as Laurent polynomials, e.g. $\Zv=\Z[\vpar,\vpar^{-1}]$.

\item All $\K$-vector spaces are finite dimensional in this paper, 
or free of finite rank if $\K$ is a ring \changed{and all categories are $\K$-linear}.

\item \changed{All modules considered are left modules, and representations of quantum enveloping algebras are always of type 1 (as defined in, e.g., \cite[Section 1.4]{AnPoWe-representation-qalgebras}).}
    
\item \changed{We often say ``is a XYZ'' instead of ``can be equipped with the structure of an XYZ'' to avoid overloading statements. For instance, we say a category is modular.}
    
\item \changed{Similarly, we say, for example, ``there are only finitely many simple objects'' instead of ``there are only finitely many isomorphism classes of simple objects.''}

\end{enumerate}

\noindent\textbf{Acknowledgments.}
We like to thank Elijah Bodish, C{\'e}dric Bonnaf{\'e}, Ben Elias \ochanged{and Simon Riche} for their invaluable insights and thoughtful discussions, and ChatGPT for help during proofreading. \ochangedd{We are indebted to the referee for their meticulous review, which resulted in a number of corrections and meaningful enhancements to the manuscript. All these contributions} significantly enriched our work. We also thank a flower of hexagons for giving us the wrong impression which killed a different project and led to this paper.

AL is grateful for the support and hospitality of the Sydney Mathematical Research Institute (SMRI). DT was sponsored by the ARC Future Fellowship FT230100489, and PV was supported by the Fonds de la Recherche Scientifique-FNRS under Grant No. CDR-J.0189.23.

\section{Some \texorpdfstring{$\sln$}{sln} combinatorics}\label{S:slnthings}

We start by fixing some notation regarding $\sln$. Most of the material is known, but our exposition for some parts is new.

\subsection{Root combinatorics}

Denote by $(\varepsilon_{1},\dots,\varepsilon_{\rank})$ the standard basis of $\R^{\rank}$, which we equip with the standard symmetric bilinear form $(\varepsilon_{i},\varepsilon_{j})=\delta_{i,j}$. Denote by $\sym_{\rank}$ the symmetric group on $\rank$ letters, which acts naturally on $\R^{\rank}$ by permutation \changed{of coordinates}. We let $E=\{(x_{1},\ldots,x_{\rank})\in \R^{\rank}\mid x_{1}+\cdots+x_{\rank} = 0\}$ and let $\alpha_{i} = \varepsilon_{i+1}-\varepsilon_{i}$ for $1\leq i < n$. The vectors $\alpha_{1},\ldots,\alpha_{\rank{-}1}$ are the simple roots and we fix the coroots $\alpha_{i}^{\vee}\in E^{\ast}$ such that $\langle \alpha_{i},\alpha_{j}^{\vee} \rangle = a_{ij}$, where $(a_{ij})$ is the usual Cartan matrix of type $A_{\rank-1}$ and $\langle\placeholder,\placeholder\rangle$ is the duality pairing. The weight lattice is $X=\{\lambda\in E \mid \langle \lambda,\alpha_{i}^{\vee}\rangle \in \Z \text{ for all } 1\leq i < \rank\}$ and the dominant weights are $X^{+}=\{\lambda\in E \mid \langle \lambda,\alpha_{i}^{\vee}\rangle \in \N \text{ for all } 1\leq i < \rank\}$. We also denote by $\omega_{1},\ldots,\omega_{N-1}\in E$ the fundamental weights which are defined through the equalities $\langle \omega_{i},\alpha_j^{\vee} \rangle = \delta_{i,j}$, and by $\rho=\omega_1+\cdots+\omega_{\rank{-}1}$ the sum of the fundamental weights, or equivalently the half-sum of positive roots. Using the basis of fundamental weights, we identify $X$ with $\Z^{\rank{-}1}$ and $X^{+}$ with $\N^{\rank{-}1}$.

\begin{Notation}
\changedd{We will repeatedly sum over entries of tuples, and we use the following shorthand notation: $\bsum\bm=m_{1}+\dots+m_{k}$ for $\bm=(m_{1},\dots,m_{k})$. 
A weight $\lambda\in X^{+}$ can be identified with a tuple $\bl=(\lambda_{1},\dots,\lambda_{\rank{-}1})$, and we make this identification tacitly. We will also write $\bsum\lambda$ instead of $\bsum\bl$ for $\lambda_1+\dots+\lambda_{\rank{-}1}$.}
\end{Notation}

We will also work with cut-offs $X^{+}(e)$ of the weight lattice, which will depend on the level $\level\in \N$, which are defined by $X^{+}(\level) = \{\lambda \in X^{+}\mid \langle\lambda,\alpha_{1}^{\vee}+\dots+\alpha_{\rank-1}^{\vee} \rangle \leq \level\}$. Therefore, a weight $\lambda=\lambda_{1}\omega_{1}+\cdots+\lambda_{\rank{-}1}\omega_{\rank{-}1}$ is in $X^{+}(\level)$ if and only if each $\lambda_{i}$ is nonnegative and \changed{$\bsum\lambda\leq\level$.}

The following picture is stolen from \cite{MaMaMiTu-trihedral}. We stole it because it summarizes our conventions for $\rank=3$:
\begin{gather*}\label{eq:weight-picture}
\scalebox{1.2}{
\begin{tikzpicture}[anchorbase, xscale=0.35, yscale=0.5]
\draw[thin, white, fill=darkpurple, opacity=0.6] (0,0) to (3,3) to (-3,3) to (0,0);
\draw[thin, white, fill=lavender, opacity=0.5] (0,0) to (4,4) to (-4,4) to (0,0);
\draw[thin, white, fill=palelavender, opacity=0.4] (0,0) to (5,5) to (-5,5) to (0,0);
\draw[thin, densely dashed, mygray] (0,0) to (7,7);
\draw[thin, densely dashed, mygray] (0,0) to (-7,7);
\node at (0,6.75) {\text{{\tiny $X^+=\N^{2}$}}};
\draw[thick, densely dotted, black] (4,3) node [right] {\text{{\tiny $e=3$}}} to (-4,3);
\draw[thick, densely dotted, black] (-5,4) node [left] {\text{{\tiny $e=4$}}} to (5,4);
\draw[thick, densely dotted, black] (6,5) node [right] {\text{{\tiny $e=5$}}} to (-6,5);
\draw[thick, black, ->] (0,0) to (3,1.5) node [right] {\text{{\tiny $\alpha_1$}}};
\draw[thick, black, ->] (0,0) to (-3,1.5) node [left] {\text{{\tiny $\alpha_2$}}};
\node at (0,0) {$\bullet$};
\node at (0,2) {$\bullet$};
\node at (3,3) {$\bullet$};
\node at (-3,3) {$\bullet$};
\node at (0,4) {$\bullet$};
\node at (1,1) {$\bullet$};
\node at (-2,2) {$\bullet$};
\node at (1,3) {$\bullet$};
\node at (-2,4) {$\bullet$};
\node at (4,4) {$\bullet$};
\node at (2,2) {$\bullet$};
\node at (-1,1) {$\bullet$};
\node at (-1,3) {$\bullet$};
\node at (2,4) {$\bullet$};
\node at (-4,4) {$\bullet$};
\node at (-5,5) {$\bullet$};
\node at (-3,5) {$\bullet$};
\node at (-1,5) {$\bullet$};
\node at (1,5) {$\bullet$};
\node at (3,5) {$\bullet$};
\node at (5,5) {$\bullet$};
\node at (-6,6) {$\bullet$};
\node at (-4,6) {$\bullet$};
\node at (-2,6) {$\bullet$};
\node at (0,6) {$\bullet$};
\node at (2,6) {$\bullet$};
\node at (4,6) {$\bullet$};
\node at (6,6) {$\bullet$};
\draw[thin, ->] (-5.5,4.75) node [left] {\text{{\tiny $(2,2)$}}} to [out=0, in=165] (-.2,4.2);
\draw[thin, ->] (5.5,3.75) node [right] {\text{{\tiny $(1,2)$}}} to [out=180, in=15] (-.8,3.2);
\draw[thin, ->] (-5.5,2.75) node [left] {\text{{\tiny $(0,2)$}}} to [out=0, in=165] (-2.2,2.2);
\end{tikzpicture}
}
.
\end{gather*}
\changed{The shaded regions are, in order from dark to light, $X^{+}(3)$, $X^{+}(4)\setminus X^{+}(3)$ and $X^{+}(5)\setminus X^{+}(4)$.}

\subsection{Quantum group generically and semisimplified}\label{SS:QuantumGroup}

Our conventions follow \cite[Chapters 4-7]{Ja-lectures-qgroups}.

For $\qpar$ a formal variable,
let $\Usln$ denote the \emph{quantum enveloping ($\Cq$-)algebra 
associated to $\sln$}, with the Hopf algebra structure as chosen in 
\cite[Section 4.8]{Ja-lectures-qgroups}.
Let us denote its category of finite dimensional 
representations by $\slncat$. This category is semisimple
and has the same combinatorics as the corresponding category for $\sln$ itself, so all the below follows from classical theory. In particular, the simple $\Usln$-representations 
(we also write $\sln$-representations for short)
are parameterized by the integral positive Weyl chamber:
\begin{gather*}
\left\{
\Ll_{\bm}|\bm=(m_{1},\dots,m_{\rank{-}1})\in X^{+}
\right\}.
\end{gather*}
\changed{Let $\transpose{\bm}=(m_{\rank{-}1},\dots,m_{1})$.
As a matter of fact, we have $(\Ll_{\bm})^{*}\cong\Ll_{\transpose{\bm}}$ for all $\bm\in X^{+}$.}

The highest weight here is $m_{1}\omega_{1}+\dots+m_{\rank{-}1}\omega_{\rank{-}1}$.
Let $\ggroup{\placeholder}$ denote the (additive)
Grothendieck group, which is a $\Z$-algebra. We have a $\Z$-basis given by (the set of the elements)
\begin{gather*}
\changed{[\Ll_{\bm}]
\in\ggroup{\slncat}}.
\end{gather*}
Scalar extension gives a $\C$-algebra:
\begin{gather*}
\ggroupc{\slncat}=\ggroup{\slncat}\otimes_{\Z}\C.
\end{gather*}
\changed{Let $\Z[\fu_{i}]=\Z[\fu_{i}\mid i\in\{1,\dots,\rank{-}1\}]$ where we use the 
$\fu_{i}$ as variables. Recall that the fundamental $\sln$-representations $\otimes$-generated $\slncat$, meaning that every simple $\sln$-representation appears as a direct summand of some $\otimes$-tensor product of the fundamental $\sln$-representations. They also commute, $\Ll_{\omega_{i}}\otimes\Ll_{\omega_{j}}\cong\Ll_{\omega_{j}}\otimes\Ll_{\omega_{i}}$ and do not satisfy any additional relation. Thus, we can see them as variables in the polynomial ring $\Z[\fu_{i}]$. That is, there is an isomorphism of rings
\begin{gather*}
\ggroup{\slncat}\xrightarrow{\cong}\Z[\fu_{i}],
[\Ll_{\omega_{i}}]\mapsto\fu_{i},
\end{gather*}
which we will use to identify $\Ll_{\omega_{i}}$ and $\fu_{i}$. (That is, abusing notation, we see them simultaneously as variables and \changedd{isomorphism classes of} $\sln$-representations.)}

Write
$\fu^{k}_{i}=\fu^{\otimes k}_{i}$, for short, and also use $\fu_{i}\fu_{j}=\fu_{j}\fu_{i}$. For $\bk\in\N^{\rank{-}1}$, write $\fu^{\bk}=\fu_{1}^{k_{1}}\cdots\fu_{\rank{-}1}^{k_{\rank{-}1}}$. \changed{In this notation, we can state the following useful fact.}

\begin{Lemma}\label{L:GrothendieckRing}
We have two bases of $\ggroup{\slncat}$ given by
$\{\ochanged{[\Ll_{\bm}]}|\bm\in X^{+}\}$
and $\{[\fu^{\bk}]|\bk\in\N^{\rank{-}1}\}$. Moreover, 
as $\Z$-algebras $\ggroup{\slncat}\cong\Z[\fu_{i}]$.
\end{Lemma}

\begin{proof}
By classical theory.
\end{proof}

\autoref{L:GrothendieckRing} motivates the definition of the \emph{change-of-basis coefficients}:
\begin{gather}\label{Eq:dNumbers}
[\Ll_{\bm}]
=
{\textstyle\sum_{\bk}}\,
d^{\bk}_{\bm}\cdot
[\fu^{\bk}],
\quad d^{\bk}_{\bm}\in\Z.
\end{gather}
Note that this sum is finite since $d^{\bk}_{\bm}=0$ unless $\bsum\bk\leq\bsum\bm$.
The numbers $d^{\bk}_{\bm}$ can be computed inductively, as explained in \autoref{SS:CPoly} below, and we have $d^{\transpose{\bk}}_{\transpose{\bm}}=d^{\bk}_{\bm}$ 
and $d^{\bm}_{\bm}=1$.

\changed{The center of $\SUN$ is $\Z/\rank\Z$, which agrees with the weight lattice modulo the root lattice. We refer to the image of a weight
in this quotient as its ``color'', following e.g. \cite{MaMaMiTu-trihedral}. One can check that all weights within $\Ll_{\bm}$ have the same color, which motivates:}

\begin{Definition}\label{definition:color-code}
Let us define colors associated 
to the simple $\Usln$-representations:
\begin{gather}\label{eq:color-code}
\chi_{c}(\Ll_{\bm})=m_{1}+2m_{2}+\dots+(\rank{-}1)m_{\rank{-}1}\in\Z/\rank\Z.
\end{gather}
We call $\chi_{c}(\Ll_{\bm})$ 
the \emph{central character} of $\Ll_{\bm}$.
\end{Definition}

\begin{Remark}
\changedd{The term central character comes from the fact that the center of the group $SL_{\rank}(\C)$ is isomorphic to $\Z/\rank\Z$.}
\end{Remark}

\begin{Example}\label{E:tetrahedron}
\changed{For $\rank=2$, this is the parity coloring of $\N$\changedd{. F}or $\rank=3$ see \cite[Example 3.16]{MaMaMiTu-trihedral}\changedd{. F}or $\rank=4$, we have the following pictures \changed{of the $\sln[4]$-weights in $X^{+}(e)$, \changedd{which} we identify with the respective $\Ll_{\bm}$. (The pictures display 3D graphs, designed to create a 3D effect with higher nodes extending into the page. \changedd{The bottom left picture also gives the coordinates.})} 
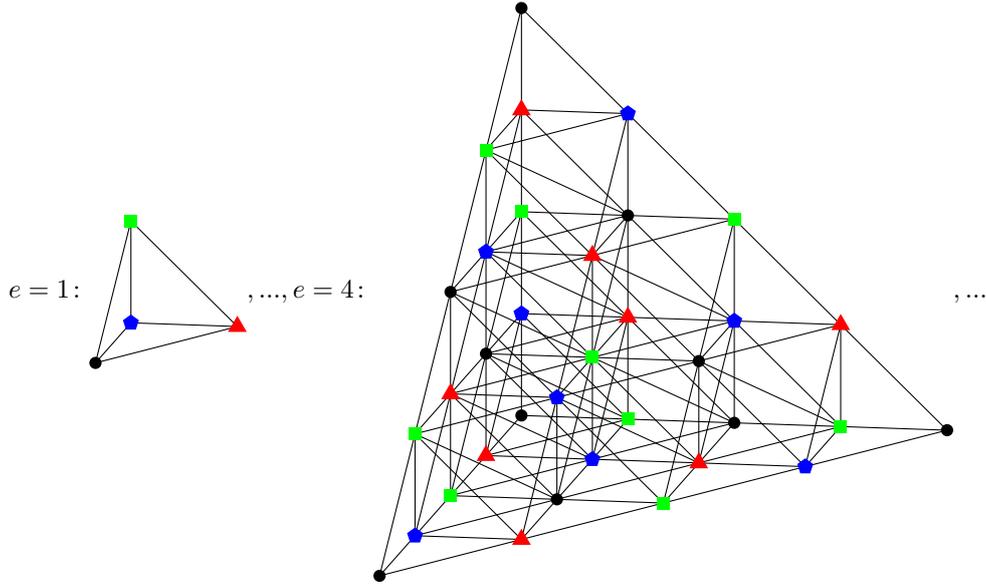
\begin{gather*}
\level=1\colon
\begin{gathered}
\input{figures/tetrahedron-1},
\\
\tdplotsetmaincoords{2}{0}
\begin{tikzpicture}[tdplot_main_coords,scale=1.4,anchorbase]
  \coordinate (o1) at (4/3,1/3,1/3);
  \coordinate (o2) at (1/3,4/3,1/3);
  \coordinate (o3) at (1/3,1/3,4/3);
  \draw (0,0,0) -- (o1);
  \draw (0,0,0) -- (o3);
  \draw (o1) -- (o2);
  \draw (o1) -- (o3);
  \draw (o2) -- (o3);
  \node[] at (o3){\fcolorbox{black}{white}{\scalebox{0.5}{(0,0,1)}}};
  \draw (0,0,0) -- (o2);
  \node[] at (0,0,0){\fcolorbox{black}{white}{\scalebox{0.5}{(0,0,0)}}};
  \node[] at (o1){\fcolorbox{black}{white}{\scalebox{0.5}{(1,0,0)}}};
  \node[] at (o2){\fcolorbox{black}{white}{\scalebox{0.5}{(0,1,0)}}};
\end{tikzpicture},
\end{gathered}
\qquad\qquad\qquad
\level=4\colon
\input{figures/tetrahedron-4}
.
\end{gather*}
The color $0$ corresponds to black circles, $1$ to red triangles, $2$ to green squares and $3$ to blue pentagons. For example, for $\level=1$ the black circle represents the trivial $\sln[4]$-representation $\Ll_{0,0,0}$, the red triangle the defining $\sln[4]$-representation $\Ll_{1,0,0}$, the blue pentagon its dual $\Ll_{0,0,1}$ and the green square the 6 dimensional simple $\sln[4]$-representation $\Ll_{0,1,0}$.}

\changed{In these pictures the edges correspond to the action (by tensoring) of the direct sum of the fundamental $\sln[4]$-representations $T=\Ll_{1,0,0}\oplus\Ll_{0,1,0}\oplus\Ll_{0,0,1}$, seen in $\slncatr$ (defined a few line below). For example,$T\otimes\Ll_{0,0,0}\cong\Ll_{1,0,0}\oplus\Ll_{0,1,0}\oplus\Ll_{0,0,1}$, which in the graph is represented by the three edges adjacent to the black circle.}
\end{Example}

These colors define a grading by $\Z/\rank\Z$ on $\slncat$\changed{, in the sense of \cite[Definition 5.9]{EtNiOs-tensor-categories},} \changedd{and we denote by $\slncat_{i}$ the full subcategory of objects of color $i$.}

\begin{Lemma}
All simple summands of $\fu^{\bk}$ have 
central character $\chi_{c}(\Ll_{\bk})$.
\end{Lemma}

\begin{proof}
All summands of $\fu^{\bk}$ have the same color\changed{,} and $\Ll_{\bk}$ is a summand of $\fu^{\bk}$.
\end{proof}

Specializing $\qpar$ to \ochangedd{the} primitive $2\levelthing$th root of unity $\rooty\ochanged{=\exp(\iunit\pi/\levelthing)}$ (using the integral form), where $\levelthing=\level+\rank$, we obtain 
$\Ursln$, as defined and used in e.g. \cite{Lu-hopf-universal-enveloping}, \cite{AnPoWe-representation-qalgebras}.
There is an associated category of representations, which we can semisimplify as usual, see e.g. \cite{AnPa-fusion-lie-algebras} for details. This category is denoted by $\slncatr$.
The Grothendieck group of this category is a $\Z$-algebra by the fusion product\changed{.}

\begin{Lemma}\label{L:Verlinde}
We have two bases of $\ggroup{\slncatr}$ given by
$\left\{
[\Ll_{\bm}]\mid \bm\in X^{+}(\level)\right\}$ and 
$\left\{
[\fu^{\bk}]\mid \bk\in X^{+}(\level)
\right\}$.
\end{Lemma}

\begin{proof}
Well-known by \cite{AnPa-fusion-lie-algebras}, see for example \cite[Example 8.18.5]{EtGeNiOs-tensor-categories}.
\end{proof}

The following are known as \emph{simplicial polytopic numbers}, since they count points in simplices.

\begin{Definition}
We define $\tnumbers{\rank}{\level}=\binom{\level+\rank{-}1}{\rank{-}1}=\binom{\levelthing-1}{\rank{-}1}$.
\end{Definition}

We have shifted them when compared to the usual definition in the sense that our $\tnumbers{\rank}{\level}$ is what is often denoted $P_{\rank{-}1}(\level-1)$.

\begin{Example}
For $\rank=2$ we have the linear numbers $\tnumbers{2}{\level}=\level+1$, for $\rank=3$ we have the triangular numbers $\tnumbers{3}{\level}=\tfrac{(\level+1)(\level+2)}{2}$, see e.g. \cite[A000217]{Oeis}.
\end{Example}

\begin{Lemma}\label{L:VerlindeAgain}
We have $\ggroup{\slncatr}\cong\Z^{\oplus\tnumbers{\rank}{\level}}$ as free $\Z$-modules.
\end{Lemma}

\begin{proof}
By \autoref{L:Verlinde}, this is just a count.
\end{proof}

The object $\Ll_{e\omega_{1}}$ is invertible in $\slncatr$ and\changed{,} for any $\bm\in X^{+}(\level)$\changed{,} we have $\Ll_{e\omega_{1}}\otimes \Ll_{\bm} \simeq \Ll_{\bm^{\turnsymbol}}$\changed{,} where $\bm^{\turnsymbol} = (\level-\bsum\bm,m_{1},\ldots,m_{\rank-2})$. \changed{The map $\bm\mapsto\bm^{\turnsymbol}$} defines an action of $\Z/\rank\Z$ on $X^{+}(e)$ and we will denote by $\stabilizer{\bm}$ the size of the stabilizer of $\bm$.

\changedd{Given $\bm \in X^{+}(\level)$, the action $\placeholder{}^{\turnsymbol}$ is more clearly understood on the extended tuple $(m_1,\dots,m_{\rank-1},e-\bsum\bm)$: it is the cyclic shift of the $\rank$ coordinates $(m_1,\dots,m_{\rank-1},e-\bsum\bm)\mapsto (e-\bsum\bm,m_1,\dots,m_{\rank-1})$.}

\begin{Example}
\changedd{For $\rank=3$ and $\level=3$, the action $\placeholder{}^{\turnsymbol}$ on $X^{+}(\level)$ is a rotation:}
\begin{gather*}
\begin{tikzpicture}[anchorbase, xscale=.7, yscale=1]
	\draw [thick] (0,0) to (1,1) to (0,2) to (1,3) to (3,3);
	\draw [thick] (0,2) to (-2,2) to (-3,3);
	\draw [thick] (0,0) to (-1,1) to (0,2) to (-1,3) to (-3,3);
	\draw [thick] (0,2) to (2,2) to (3,3);
	\draw [thick] (1,1) to (-1,1) to (-2,2) to (-1,3) to (1,3) to (2,2) to (1,1);
	\node at (0,0) {$\bullet$};
    \node at (0,-0.3) {$(0,0)$};
	\node at (0,2) {$\bullet$};
	\node at (3,3) {$\bullet$};
    \node at (3.7,3) {$(3,0)$};
	\node at (-3,3) {$\bullet$};
    \node at (-3.7,3) {$(0,3)$};
	\node at (1,1) {$\bullet$};
	\node at (-2,2) {$\bullet$};
	\node at (1,3) {$\bullet$};
	\node at (2,2) {$\bullet$};
	\node at (-1,1) {$\bullet$};
	\node at (-1,3) {$\bullet$};
	\draw [->, purple] (3,3.2) to [out=170, in=10] (-3,3.2);
	\draw [->, purple] (-3,2.8) to [out=280, in=170] (-.3,0);
	\draw [->, purple] (.3,0) to [out=10, in=260] (3,2.8);
	\draw [<-, blue] (1,2.8) to (1,1.2);
	\draw [<-, blue] (.8,1.1) to (-1.7,1.9);
	\draw [<-, blue] (-1.7,2.1) to (.8,2.9);
	\draw [->, spinach] (-1,2.8) to (-1,1.2);
	\draw [->, spinach] (-.8,1.1) to (1.7,1.9);
	\draw [->, spinach] (1.7,2.1) to (-.8,2.9);
\end{tikzpicture}
.
\end{gather*}
\changedd{For $\rank=4$, the action $\placeholder{}^{\turnsymbol}$ on the tetrahedron $X^{+}(\level)$ is a rotation composed with a reflection, corresponding to a $4$-cycle on the extremal vertices of the tetrahedron. Explicitly with $\level=22$, one has $(2,3,11)^{\turnsymbol} = (6,2,3)$ and $(1,10,1)^{\turnsymbol} = (10,1,10)$. Note that the extended tuple $(1,10,1,10)$ is $2$-periodic and therefore $(1,10,1)$ has a stabilizer of order $\rank/2$.}
\end{Example}

As for $\slncat$, we have a grading by $\Z/\rank\Z$ on $\slncatr$ by colors \ochangedd{and we denote by $\slncatr_{i}$ the full subcategory of objects of color $i$}. For later use we recall:

\begin{Lemma}\label{L:Grading}
The grading on $\slncatr$ by colors induces a 
grading on its Drinfeld center, also by colors.
\end{Lemma}

\begin{proof}
Easy and omitted.
\end{proof}

We denote the full subcategory of $\slncatr$ of objects of color $i$ by $\slncatr_{i}$, and similarly for the Drinfeld center.

\subsection{Braiding and ribbon structure}

The fusion category $\slncatr$ is equipped with a braided structure $\br_{\placeholder,\placeholder}$. This structure is obtained through the use of the universal $R$-matrix of $\Ursln$, as one can adapt from e.g. \cite[Chapter 7]{Ja-lectures-qgroups}, and a choice of an $\rank$th root of $\rooty$, which we denote by $\rooty^{1/\rank}$.
\changedd{Denote by $\rankroot$ the $\rank$th root of unity $\rooty^{2\levelthing/\rank}=\exp(2\iunit\pi/\rank)$.} \changed{We also endow $\slncatr$ with a spherical structure, the corresponding quantum trace $\Tr$ and the corresponding ribbon structure $\theta$.
(For the definition of a ribbon structure, we refer to \cite[Definition 8.10.1]{EtGeNiOs-tensor-categories}.)} We will denote by $S_{\bm,\bm^{\prime}}$ the values of the $S$-matrix as defined in \cite[Definition 8.13.2]{EtGeNiOs-tensor-categories}.

\begin{Lemma}\label{L:Smat-sln}
For any $\bm,\bm^{\prime}\in X^{+}(\level)$, we have the following formula for the $S$-matrix of $\slncatr$:
\begin{gather*}
S^{\sln[]}_{\bm,\bm^{\prime}}=\frac{\sum_{w\in \sym_{\rank}}(-1)^{l(w)}\rooty^{2(\bm+\rho,w(\bm^{\prime}+\rho))}}{\sum_{w\in \sym_{\rank}}(-1)^{l(w)}\rooty^{2(\rho,w(\rho))}}.
\end{gather*}
For any $\bm\in X^{+}(\level)$, the ribbon $\theta_{\bm}$ on $\Ll_{\bm}$ is given by multiplication by the scalar $\rooty^{(\bm,\bm+2\rho)}$.
\end{Lemma}

\begin{proof}
A proof of the lemma can be found in e.g. \cite[Theorem 3.3.20]{BaKi-lecture-tensor}.
\end{proof}

Using Brugui{\`e}res's criterion \cite[Section 5]{Br-cat-modulaires}, the fusion category $\slncatr$ is modular if and only if $\rooty^{1/\rank}$ is a primitive $2\levelthing\rank$th root of unity.

\begin{Lemma}\label{L:GradingBraiding}
Let $\bm\in X^{+}(\level)$. Then $\br_{\Ll_{\level\omega_{1}},\Ll_{\bm}}\circ\br_{\Ll_{\bm},\Ll_{\level\omega_{1}}} = \rankroot^{-\chi_{c}(\Ll_{\bm})}\cdot\idmor$.
\end{Lemma}

\begin{proof}
Since $\Ll_{\bm^{\turnsymbol}}=\Ll_{\bm}\otimes\Ll_{\level\omega_{1}}$ and since $\theta$ is a ribbon element, we have $\theta_{\bm^{\turnsymbol}} = \changed{(\theta_{\bm}\otimes\theta_{\level\omega_{1}})}\circ\br_{\Ll_{\level\omega_{1}},\Ll_{\bm}}\circ\br_{\Ll_{\bm},\Ll_{\level\omega_{1}}}$. Using \autoref{L:Smat-sln}, we therefore obtain $\br_{\Ll_{\level\omega_{1}},\Ll_{\bm}}\circ\br_{\Ll_{\bm},\Ll_{\level\omega_{1}}}=\rooty^{(\bm^{\turnsymbol},\bm^{\turnsymbol}+2\rho)-(\bm,\bm+2\rho)-(\level\omega_{1},\level\omega_{1}+2\rho)}\cdot\idmor$. One may finally show that $(\bm^{\turnsymbol},\bm^{\turnsymbol}+2\rho)-(\bm,\bm+2\rho)-(\level\omega_{1},\level\omega_{1}+2\rho) = -2\level\chi_{c}(\Ll_{\bm})/\rank$.
\end{proof}

\changedd{Since the exponent of $\rankroot$ in \autoref{L:GradingBraiding} is the opposite of the color,} the grading by color is then retrieved using the \changedd{double} braiding \ochangedd{similar to} \cite[Lemma 8.22.1]{EtGeNiOs-tensor-categories}.

\subsection{Chebyshev-like polynomials}\label{SS:CPoly}

The following polynomials are due to Koornwinder \cite{Ko}, and 
Eier--Lidl \cite{EiLi-orthogonal-polynomials}.
We mildly change the conventions for convenience.

\begin{Definition}\label{D:CKpoly}
For each $\bm\changed{\in X^{+}}$ we define \emph{Chebyshev polynomials (of the second kind) for $\sln$}, denoted by
$\pxy{\bm}\in\Z[\fu_{i}]$\changed{,} as
\begin{gather*}
\pxy{\bm}=
{\textstyle\sum_{\bk}}\,
d_{\bm}^{\bk}
\cdot
\changed{[\fu^{\bk}]},
\end{gather*}
with $d_{\bm}^{\bk}\in\Z$ as in \autoref{Eq:dNumbers}. 

By convention, $\pxy{\bm}$ and $\Ll_{\bm}$ with negative subscripts $m_{i}$ are zero.
\end{Definition}

\begin{Lemma}\label{L:ChebyshevRecursion}
Let $w_{j}^{i}$ denote the weights of $\Ll_{\omega_{i}}$.
\changed{For $\bm\in X^+$, we} have the following Chebyshev-like recursion relations 
\begin{gather*}
\pxy{\bm}(\fu_{1},\dots,\fu_{\rank-1})=\pxy{\transpose{\bm}}(\fu_{\rank-1},\dots,\fu_{1})
,\quad
\fu_{i}\pxy{\bm}
=
\sum_{\ochanged{j}}\pxy{\bm+w_{j}^{i}}.
\end{gather*}
Together with the starting conditions for $\changed{\bsum\bm}=0,1$ (as in \autoref{example:sl3-polys} for $\rank=4$), these 
recursion relations determine the polynomials 
$\pxy{\bm}$ for all $\bm$.
\end{Lemma}

\begin{proof}
By construction and classical theory \changed{(such as the dual Pieri rule)}.
\end{proof}

\begin{Example}\label{example:sl3-polys}
The examples for $\rank=2$ are the Chebyshev polynomials of the second kind normalized using the variable $x/2$ instead of $x$, while \cite[Example 2.6]{MaMaMiTu-trihedral} lists examples for $\rank=3$. 

The next case is $\rank=4$, 
we have three fundamental variables, $\fu_{1}$, $\fu_{2}$ and $\fu_{3}$, associated to the simple $4$ (the vector representation), $6$ and $4$ dimensional $\sln[4]$-representations.
We have the following recursion. If $m_{i}<0$, then $\pxy{m_{1},m_{2},m_{3}}=0$, and $\pxy{0,0,0}=1$ and:
\begin{gather*}
\pxy{m_{1},m_{2},m_{3}}=
\fu_{1}\cdot\pxy{m_{1}-1,m_{2},m_{3}}
-\pxy{m_{1}-2,m_{2}+1,m_{3}}-\pxy{m_{1}-1,m_{2}-1,m_{3}+1}-\pxy{m_{1}-1,m_{2},m_{3}-1},
\\
\pxy{m_{1},m_{2},m_{3}}=
\left\{\begin{aligned}
&\fu_{2}\cdot\pxy{m_{1},m_{2}-1,m_{3}} - \pxy{m_{1}+1,m_{2}-2,m_{3}+1} - \pxy{m_{1}-1,m_{2}-1,m_{3}+1} \\&- \pxy{m_{1}+1,m_{2}-1,m_{3}-1} - \pxy{m_{1}-1,m_{2},m_{3}-1} - \pxy{m_{1},m_{2}-2,m_{3}},
\end{aligned}\right.
\\
\pxy{m_{1},m_{2},m_{3}} = \fu_{3}\cdot\pxy{m_{1},m_{2},m_{3}-1} - \pxy{m_{1},m_{2}+1,m_{3}-2} - \pxy{m_{1}+1,m_{2}-1,m_{3}-1} - \pxy{m_{1}-1,m_{2},m_{3}-1}.
\end{gather*}
Now $\pxy{0,0,0}=1$ and
\begin{gather*}
\begin{tikzpicture}[baseline=(current bounding box.center)]
\matrix (m) [matrix of math nodes, row sep=0.2cm, column
sep=0.0cm, text height=1.5ex, text depth=0.25ex, ampersand replacement=\&, font=\scriptsize,
nodes={anchor=south},
row 3/.style={nodes={minimum height=1cm}},
row 4/.style={nodes={minimum height=1.5cm}},
row 5/.style={nodes={minimum height=1.5cm}},
] {
\changed{\bsum\bm = 1} \& 
\pxy{1,0,0}=\fu_{1},
\;
\pxy{0,1,0}=\fu_{2},
\;
\pxy{0,0,1}=\fu_{3},
\\
\changed{\bsum\bm = 2} \&
\begin{gathered}
\pxy{2,0,0}=\fu_{1}^{2}-\fu_{2},
\;
\pxy{1,1,0}=\fu_{1}\fu_{2}-\fu_{3},
\;
\pxy{1,0,1}=\fu_{1}\fu_{3}-1,
\\
\pxy{0,2,0}=\fu_{2}^{2}-\fu_{1}\fu_{3},
\;
\pxy{0,1,1}=\fu_{2}\fu_{3}-\fu_{1},
\;
\pxy{0,0,2}=\fu_{3}^{2}-\fu_{2},
\end{gathered}
\\
\changed{\bsum\bm = 3} \&
\begin{gathered}
\pxy{3,0,0}=\fu_{1}^{3}-2\fu_{1}\fu_{2}+\fu_{3},
\;
\pxy{2,1,0}=\fu_{1}^{2}\fu_{2}-\fu_{1}\fu_{3}-\fu_{2}^{2}+1,
\;
\pxy{2,0,1}=\fu_{1}^{2}\fu_{3}-\fu_{2}\fu_{3}-\fu_{1},
\;
\pxy{1,2,0}=\fu_{1}\fu_{2}^{2}-\fu_{1}^{2}\fu_{3}-\fu_{2}\fu_{3}+\fu_{1},
\\
\pxy{1,1,1}=\fu_{1}\fu_{2}\fu_{3}-\fu_{1}^{2}-\fu_{3}^{2},
\;
\pxy{1,0,2}=\fu_{1}\fu_{3}^{2}-\fu_{1}\fu_{2}-\fu_{3},
\;
\pxy{0,3,0}=\fu_{2}^{3}-2\fu_{1}\fu_{2}\fu_{3}+\fu_{1}^{2}+\fu_{3}^{2}-\fu_{2},
\\
\pxy{0,2,1}=\fu_{2}^{2}\fu_{3}-\fu_{1}\fu_{3}^{2}-\fu_{1}\fu_{2}+\fu_{3},
\;
\pxy{0,1,2}=\fu_{2}\fu_{3}^{2}-\fu_{1}\fu_{3}-\fu_{2}^{2}+1,
\;
\pxy{0,0,3}=\fu_{3}^{3}-2\fu_{2}\fu_{3}+\fu_{1}.
\;
\end{gathered}
\\
};
\draw[thin, densely dotted] ($(m-1-1)+(-.5,-.25)$) to ($(m-1-1)+(15.1,-.25)$);
\draw[thin, densely dotted] ($(m-1-1)+(-.5,-1.1)$) to ($(m-1-1)+(15.1,-1.1)$);
%
\draw[thin, densely dotted] ($(m-1-1)+(.5,.25)$) to ($(m-1-1)+(.5,-2.3)$);
\end{tikzpicture}
\end{gather*}
For $\changed{\bsum\bm = 4}$ there are already $15$ polynomials and for $\changed{\bsum\bm = 5}$ there are $21$, so we omit to put them here. Instead let us list the polynomials for the symmetric powers $U^{(k)}=\pxy{(k,0,0)}$ varying $k$ \changed{(the shading indicates when the constant coefficient is non zero)}:
\begin{center}
\begin{tabular}{C|C|C|C|C|C|C}
k & \cellcolor{gray!50}0 & 1 & 2 & 3 & \cellcolor{gray!50}4 & 5 \\
\hline
U^{(k)} & \cellcolor{gray!50}1 & \fu_{1} & \fu_{1}^{2}{-}\fu_{2} & \fu_{1}^{3}{-}2\fu_{1}\fu_{2}{+}\fu_{3} & \cellcolor{gray!50}\fu_{1}^{4}{-}3\fu_{1}^{2}\fu_{2}{+}2\fu_{1}\fu_{3}{+}\fu_{2}^{2}{-}1 & \fu_{1}^{5}{-}4\fu_{1}^{3}\fu_{2}{+}3\fu_{1}^{2}\fu_{3}{+}3\fu_{1}\fu_{2}^{2}{-}2\fu_{2}\fu_{3}{-}2\fu_{1}
\\
\end{tabular}
\end{center}
\begin{center}
\begin{tabular}{C|C|C|C|C|C|C|C}
k & 6 & 7 & \cellcolor{gray!50}8 & 9 & 10 & 11 & \cellcolor{gray!50}12 \\
\hline
U^{(k)} & \fu_{1}^{6}{\pm}\dots{+}2\fu_{2} & \fu_{1}^{7}{\pm}\dots{-}2\fu_{3} & \cellcolor{gray!50}\fu_{1}^{8}{\pm}\dots{+}1 & \fu_{1}^{9}{\pm}\dots{+}3\fu_{1} & \fu_{1}^{10}{\pm}\dots{-}3\fu_{2} & \fu_{1}^{11}{\pm}\dots{+}3\fu_{3} & \cellcolor{gray!50}\fu_{1}^{12}{\pm}\dots{-}1 \\
\end{tabular}
\end{center}
We will see in the proof of \autoref{L:ChebConstant} below how one can compute these fairly efficiently. 
\end{Example}

\begin{Lemma}\label{L:ChebConstant}
The polynomial $\pxy{\bm}$ has a nonzero 
constant term only if $\chi_{c}(\Ll_{\bm})=0$. 
Moreover, \changed{for a given $e\in\N$, all Chebyshev polynomial\changedd{s} $\pxy{\bm}$ with $\bsum\bm=\level+1$ have a zero constant term if and only if $e\equiv 0 \bmod \rank$}.
\end{Lemma}

\begin{proof}
This proof is much more involved than in \cite[Lemma 2.8]{MaMaMiTu-trihedral}.

The trivial representation is of color $0$, which implies that if $\chi_{c}(\Ll_{\bm})\neq 0$ then $U_{\bm}$ has a zero constant coefficient.

The vanishing of the constant term is understood using the recursion \cite[(5.18)]{Be-Chebyshev-Laplace-Beltrami}, which is nicely expressed using partitions instead of highest weights. This recursion implies that the Chebyshev polynomial $\pxy{\bm}$ has a non zero constant coefficient if and only if the residues modulo $\rank$ of $(m_{1}+\dots+m_{\rank-1}+\rank-1,m_{2}+\dots+m_{\rank-1}+\rank-2,\dots,m_{\rank-1}+1,0)$ are all different. Therefore, if $\bsum\bm \equiv 1 \mod \rank$, the residue $0$ appears at least twice and $\pxy{\bm}$ has zero constant term. Conversely, given $\level\equiv k \bmod\rank$ with $0<k<\rank$, one may check that for $\bm=\omega_{k}+e\omega_{\rank-1}$, the Chebyshev polynomial $\pxy{\bm}$ has a nonzero constant term.
\end{proof}

\begin{Remark}
\changed{To be self-contained, the alternative recursion that we mention above comes from the following observation: It is 
remarkably easy to find a recursion for the $k$th symmetric power of the vector $\sln$-representation, which correspond\changedd{s} to $\pxy{(k,0,\dots,0)}$ in our notation. For $\rank=2$ this recursion is the standard recursion since all simple $\sln[2]$-representations are symmetric powers. For $\rank>2$ let $U^{(k)}(\fu_{1},\dots,\fu_{\rank-1})=\pxy{(k,0,\dots,0)}(\fu_{1},\dots,\fu_{\rank-1})$
denote these polynomials. The recursion then takes the form
\begin{gather*}
U^{(k+\rank)}-\fu_{1}U^{(k+\rank-1)}+\dots(-1)^{\rank-1}\fu_{\rank-1}U^{(k+1)}+(-1)^{\rank}U^{(k)}=0,
\end{gather*}
with some additional starting conditions. \changedd{The main observation made in \cite[(5.18)]{Be-Chebyshev-Laplace-Beltrami} is that this recursion also takes a nice form in the partition notation for the highest weights.} In this recursion setting all variable\changedd{s} to zero (the constant term) gives
\begin{gather*}
U^{(k+\rank)}+(-1)^{\rank}U^{(k)}=0,
\end{gather*}
and the claim then follows easily.
}
\end{Remark}

A classical result of Kostant \cite{Ko-Macdonald-eta-function} 
gives a formula for certain powers of the Dedekind $\eta$-function
by summing over simple $\sln$-representations $\Ll_{\bm}$, and \changed{the coefficients $\epsilon_{\bm}$ which appear} are in $\{0,1,-1\}$.
The following is a fun side observation:

\begin{Proposition}
The constant coefficient of the Chebyshev polynomial $\pxy{\bm}$ is equal to the trace of a(ny) Coxeter element acting on the zero weight space of $\Ll_{\bm}$ which in turn is equal to $\epsilon_{\bm}$.
\end{Proposition}

\begin{proof}
This follows by comparing \cite[(5.18)]{Be-Chebyshev-Laplace-Beltrami} and \cite[Theorem 1.2]{AdFr-rim-hook}.
\end{proof}

\subsection{Koornwinder variety}\label{SS:Koornwinder}

All proofs in this section are much more intricate than their rank $2$ or $3$ counterparts in \cite{MaMaMiTu-trihedral}. We define:

\begin{Definition}\label{D:KoornwinderV}
Let $\vanideal{\level}$ be the ideal generated by
\begin{gather*}
\left\{
\pxy{\bm}\mid\bsum\bm=e+1
\right\}
\subset\Z[\fu_{i}].
\end{gather*}
We call $\vanideal{\level}$ the \emph{vanishing ideal} of level $\level$. Associated 
to it is the \emph{Koornwinder variety} of level $\level$
\begin{gather*}
\vanset{\level}=
\left\{
\bg=(\gamma_{1},\dots,\gamma_{\rank{-}1})\in\C^{\rank{-}1}\mid 
p(\bg)=0\;\text{for all}\; p\in\vanideal{e}
\right\}\subset\C^{\rank{-}1}
\end{gather*}
which we consider as a complex variety.
\end{Definition}

\begin{Remark}
Following history, one could also call $\vanset{\level}$ the
\emph{Chebyshev--Eier--Koornwinder--Lidl variety}, cf. \cite{Ko}, which 
discusses the case $\rank=3$ and \cite{EiLi-orthogonal-polynomials}, which 
discusses the general case, but that is a mouthful.
\end{Remark}

\begin{Example}
All roots of Chebyshev polynomials $\pxy{\bm}$ have their first coordinate in the interior of the \emph{$\rank$-cusped hypocycloid} \changed{defined by the parametric equation}
\begin{gather*}
x(\theta)=(\rank-1)\cos(\theta)+\cos\big((\rank-1)\theta\big)\text{ and }y(\theta)=(\rank-1)\sin(\theta)-\sin\big((\rank-1)\theta\big).
\end{gather*}
This is a folk result, and can, for example, be explicitly found in
\cite[Section 3]{Ka-traces}. Recall that these are plane curves generated by following a point on a circle of radius $1$ that rolls within a circle of radius $\rank$. Here are the pictures \changed{(thoughout, we identify the plane with complex numbers)}:
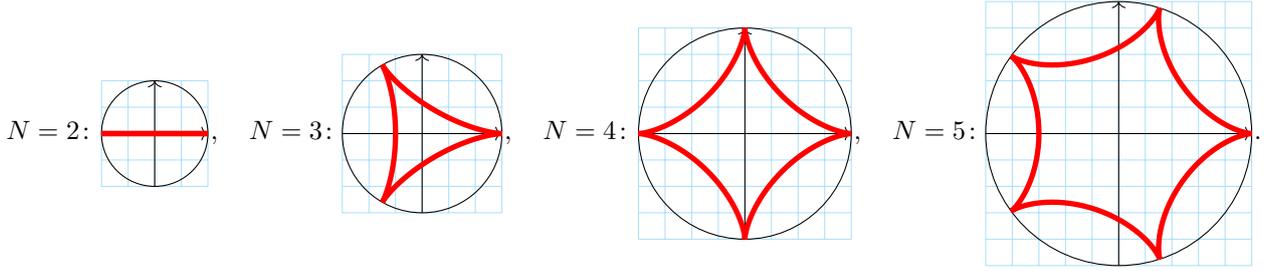
\begin{gather*}
\rank=2\colon
\begin{tikzpicture}[anchorbase,scale=0.35]
\def\a{1} \def\b{2}
\draw[cyan!30,very thin] (-\b,-\b) grid (\b,\b);
\draw[->] (-\b,0) -- (\b,0);
\draw[->] (0,-\b) -- (0,\b);
\draw (0,0) circle (\b);
\draw[line width=2pt,red] plot[samples=100,domain=0:\a*360,smooth,variable=\t] ({(\b-\a)*cos(\t)+\a*cos((\b-\a)*\t/\a},{(\b-\a)*sin(\t)-\a*sin((\b-\a)*\t/\a}); <-- edited (to add \a*360)
\end{tikzpicture}
,\quad
\rank=3\colon
\begin{tikzpicture}[anchorbase,scale=0.35]
\def\a{1} \def\b{3}
\draw[cyan!30,very thin] (-\b,-\b) grid (\b,\b);
\draw[->] (-\b,0) -- (\b,0);
\draw[->] (0,-\b) -- (0,\b);
\draw (0,0) circle (\b);
\draw[line width=2pt,red] plot[samples=100,domain=0:\a*360,smooth,variable=\t] ({(\b-\a)*cos(\t)+\a*cos((\b-\a)*\t/\a},{(\b-\a)*sin(\t)-\a*sin((\b-\a)*\t/\a}); <-- edited (to add \a*360)
\end{tikzpicture}
,\quad
\rank=4\colon
\begin{tikzpicture}[anchorbase,scale=0.35]
\def\a{1} \def\b{4}
\draw[cyan!30,very thin] (-\b,-\b) grid (\b,\b);
\draw[->] (-\b,0) -- (\b,0);
\draw[->] (0,-\b) -- (0,\b);
\draw (0,0) circle (\b);
\draw[line width=2pt,red] plot[samples=100,domain=0:\a*360,smooth,variable=\t] ({(\b-\a)*cos(\t)+\a*cos((\b-\a)*\t/\a},{(\b-\a)*sin(\t)-\a*sin((\b-\a)*\t/\a}); <-- edited (to add \a*360)
\end{tikzpicture}
,\quad
\rank=5\colon
\begin{tikzpicture}[anchorbase,scale=0.35]
\def\a{1} \def\b{5}
\draw[cyan!30,very thin] (-\b,-\b) grid (\b,\b);
\draw[->] (-\b,0) -- (\b,0);
\draw[->] (0,-\b) -- (0,\b);
\draw (0,0) circle (\b);
\draw[line width=2pt,red] plot[samples=100,domain=0:\a*360,smooth,variable=\t] ({(\b-\a)*cos(\t)+\a*cos((\b-\a)*\t/\a},{(\b-\a)*sin(\t)-\a*sin((\b-\a)*\t/\a}); <-- edited (to add \a*360)
\end{tikzpicture}
.
\end{gather*}
To be completely explicit, let $\rank=4$ and $\level=2$. The \changed{points in the corresponding Koornwinder variety are}:
\begin{gather*}
(0,-1,0),(0,1,0),(\sqrt{3}\iunit,-2,-\sqrt{3}\iunit),(-\sqrt{3}\iunit,-2,\sqrt{3}\iunit),(\sqrt{3},2,\sqrt{3}),(-\sqrt{3},2,-\sqrt{3}),\\(e^{2\iunit\pi/8},0,e^{14\iunit\pi/8}),(e^{6\iunit\pi/8},0,e^{10\iunit\pi/8}),(e^{10\iunit\pi/8},0,e^{6\iunit\pi/8}),(e^{14\iunit\pi/8},0,e^{2\iunit\pi/8})\changed{\in\C^{3}}.
\end{gather*}
We have the following plot of the first coordinates:
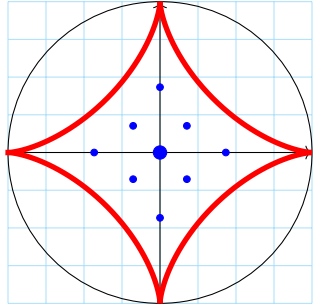
\begin{gather*}
\begin{tikzpicture}[anchorbase,scale=0.5]
\def\a{1} \def\b{4}
\draw[cyan!30,very thin] (-\b,-\b) grid (\b,\b);
\draw[->] (-\b,0) -- (\b,0);
\draw[->] (0,-\b) -- (0,\b);
\draw (0,0) circle (\b);
\draw[line width=2pt,red] plot[samples=100,domain=0:\a*360,smooth,variable=\t] ({(\b-\a)*cos(\t)+\a*cos((\b-\a)*\t/\a},{(\b-\a)*sin(\t)-\a*sin((\b-\a)*\t/\a}); <-- edited (to add \a*360)
\draw[blue,fill=blue] (0,0) circle (5.0pt);
\draw[blue,fill=blue] ({sqrt(3)},0) circle (2.5pt);
\draw[blue,fill=blue] (0,{sqrt(3)}) circle (2.5pt);
\draw[blue,fill=blue] ({-sqrt(3)},0) circle (2.5pt);
\draw[blue,fill=blue] (0,{-sqrt(3)}) circle (2.5pt);
\draw[blue,fill=blue] (45:1) circle (2.5pt);
\draw[blue,fill=blue] (135:1) circle (2.5pt);
\draw[blue,fill=blue] (225:1) circle (2.5pt);
\draw[blue,fill=blue] (315:1) circle (2.5pt);
\end{tikzpicture}
\end{gather*}
The point at zero is illustrated thick since two \changed{points in the Koornwinder variety share zero as the first coordinate ($(0,-1,0)$ and $(0,1,0)$)}. Similar conventions are used throughout, i.e. the thickness of points indicates their multiplicity.
\end{Example}

\begin{Example}
For $\rank=6$ and $\level=6$, the Koornwinder variety has $\tnumbers{6}{6}=462$ points\changed{, as one can check with MAGMA for example}. We have the following \changed{plot of the first coordinate, resp. of the second coordinate}:
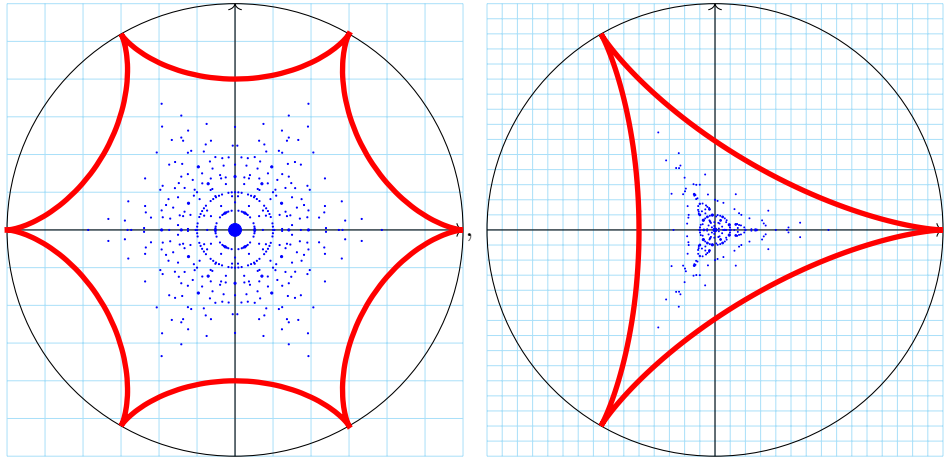
\begin{gather*}
\input{figures/large-pictures}    
\end{gather*}
We have not included the third coordinate since the points are on the real line. The fourth and fifth coordinates are the complex conjugates of the second and the first.
Note that the radius of the circle of the $i$th coordinate equals the dimension of the $\sln$-representation $\Ll_{\omega_{i}}$.
\end{Example}

\begin{Definition}\label{D:Zfunctions}
For $1\leq i\leq\rank-1$, we introduce \emph{Koornwinder's Z-functions}
$\Zfun_{i}\colon E\to\C$ defined by
\begin{gather*}
\Zfun_{i}(\bs) = \sum_{j}\exp\left(\iunit(\bs,w^{i}_{j})\right),
\end{gather*}
where we recall that $(w_{j}^{i})_{j}$ are the weights of the fundamental representation $\Ll_{\omega_{i}}$.
\end{Definition}

\begin{Remark}
The map $\Zfun_{1}$ is $2\pi Y$ periodic, where $Y$ is the root lattice, and is invariant under the action of $\sym_{\rank}$. The fundamental domain of $E/2\pi Y$ for this action is
\begin{gather*}
D=\left\{\sum_{i=1}^{\rank-1}\lambda_{i}\alpha_{i} \Big| 2\lambda_{i}\geq \lambda_{i-1}+\lambda_{i+1} \text{ for all } 1\leq i <N, \lambda_{1}+\lambda_{\rank-1}\leq 2\pi\right\}.
\end{gather*}
Note that, in contrast to \cite{Ko} (this is $\rank=3$), the map $\Zfun_{1}$ is not injective on the fundamental domain $D$ for $\rank > 3$. For example, for $\rank=4$, $D$ is a tetrahedron and the image of the six edges under $\Zfun_{1}$ are the following: 
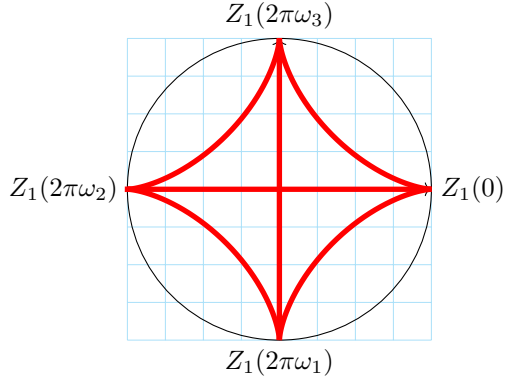
\begin{gather*}
\begin{tikzpicture}[anchorbase,scale=0.5]
\def\a{1} \def\b{4}
\draw[cyan!30,very thin] (-\b,-\b) grid (\b,\b);
\draw[->] (-\b,0) -- (\b,0);
\draw[->] (0,-\b) -- (0,\b);
\draw (0,0) circle (\b);
\draw[line width=2pt,red] plot[samples=100,domain=0:\a*360,smooth,variable=\t] ({(\b-\a)*cos(\t)+\a*cos((\b-\a)*\t/\a},{(\b-\a)*sin(\t)-\a*sin((\b-\a)*\t/\a}); <-- edited (to add \a*360)
\draw[line width=2pt,red] (-\b,0) -- (\b,0);
\draw[line width=2pt,red] (0,-\b) -- (0,\b);
\node[right] at (\b,0){$\Zfun_{1}(0)$};
\node[left] at (-\b,0){$\Zfun_{1}(2\pi\omega_{2})$};
\node[below] at (0,-\b){$\Zfun_{1}(2\pi\omega_{1})$};
\node[above] at (0,\b){$\Zfun_{1}(2\pi\omega_{3})$};
\end{tikzpicture}
\end{gather*}
The edge joining the vertices $0$ and $2\pi\omega_{2}$ and the one joining the vertices $2\pi\omega_{1}$ and $2\pi\omega_{3}$ are mapped to the parts of the real and the imaginary axis inside the hypocycloid.
\end{Remark}

To describe $\vanset{\level}$ using the functions $\Zfun_{i}$, we introduce a different parametrization.

\begin{Definition}\label{D:KoornwinderPrime}
Define $\vanpar{\level}$ as the set 
$\vanpar{\level}=\left\{\tfrac{2\pi}{\level+\rank}(\bk+\rho)\mid\bk\in X^{+}(\level)\right\}$. For $\bs=\frac{2\pi}{\level+\rank}(\bk+\rho)\in\vanpar{\level}$, we denote by $\bs^{\turnsymbol}$ the element $\frac{2\pi}{\level+\rank}(\bk^{\turnsymbol}+\rho)\in\vanpar{\level}$.
\end{Definition}

\begin{Lemma}\label{L:KoornwinderVarietyZ}
Given $\bs\in\vanpar{\level}$, we have $(\Zfun_{1}(\bs),\dots,\Zfun_{\rank-1}(\bs))\in \vanset{\level}$.
\end{Lemma}

\begin{proof}
See \autoref{Pf:KoornwinderVarietyZ}.
\end{proof}

\begin{Example}
For $\rank=3$, the coordinates of the vectors of $\vanpar{\level}$ in the basis $(\alpha_{1},\alpha_{2})$ coincide with the expression as in \cite[(2-11)]{MaMaMiTu-trihedral}.
For $\rank=4$, we get
\begin{gather*}
\vanpar{\level}=
\left\{
\tfrac{2\pi}{4(\level+4)}
(3k_{1}+2k_{2}+1k_{3}+6,
2k_{1}+4k_{2}+2k_{3}+8,
1k_{1}+2k_{2}+3k_{3}+6)
\mid 0\leq k_{1}+k_{2}+k_{3}\leq\level
\right\},
\end{gather*}
from \autoref{D:KoornwinderPrime}, where the coordinates are given in the basis $(\alpha_{1},\alpha_{2},\alpha_{3})$.
\end{Example}

In the next theorem we use the M{\"o}bius function $\mu$.

\begin{Theorem}\label{T:KoornwinderVariety}
We have the following:
\begin{enumerate}
\item As $\Z$-algebras we have $\ggroup{\slncatr}\cong\Z[\fu_{i}]/\vanideal{\level}$.
\item $\#\vanset{\level}=\tnumbers{\rank}{\level}$.
\item The number of points in $\vanset{\level}$ with stabilizer of size $m$ under multiplication by $\rankroot$ is
\begin{gather*}
\frac{\rank}{\levelthing}\sum_{k\mid g}\mu(k)\binom{\levelthing/mk}{\rank/mk},
\end{gather*}
where $g=\gcd(\rank/m,\levelthing/m)$.
\end{enumerate}
\end{Theorem}

\begin{proof}
\textit{Part (a).} \changed{As explained in \cite[Section 5.a]{Be-Chebyshev-Laplace-Beltrami}, imported to our setting using \cite[Section 5.d]{Be-Chebyshev-Laplace-Beltrami}, the Chebyshev polynomial $\pxy{\bm}$ is the character of the classical analog of $\Ll_{\bm}$. Therefore $[\Ll_{\bm}] \mapsto \pxy{\bm}$ is a well-defined ring isomorphism $\ggroup{\slncat}\cong\Z[\fu_{i}]$. The claimed isomorphism then follows from \autoref{L:Verlinde} and \cite[Proposition 1.4]{AnSt}, which shows that the negligible ideal is tensor generated by $\{\Ll_{\bm} \mid \bsum\bm=\level+1\}$ (this could also be derived from quantum Racah formula as e.g. in \cite[Corollary 8]{Sa-quantum-roots-of-unity}).}

\changedd{\textit{Part (b).} This is proven in \autoref{Pf:KoornwinderVariety-b}.}

\changedd{\textit{Part (c).} This is proven in \autoref{Pf:KoornwinderVariety-c}.}
\end{proof}

\section{Nhedral Hecke algebras and categories}\label{S:Nhedral}

While reading this section we recommend having 
\autoref{L:GrothendieckRing} and \autoref{T:KoornwinderVariety}.(a) \changeddd{in mind.}
that we are going to ``color''.

\subsection{Nhedral Hecke algebras}\label{SS:NHeckeAlgebras}

Let $\nodes=\{0,\dots,\rank{-}1\}$, identified with the set of vertices of the affine type $A_{\rank{-}1}$ Dynkin diagram. The associated Weyl group is $\weyl=\langle
s_{i}|i\in\nodes\rangle/(\text{relations in \autoref{Eq:Dynkin}})$.

\begin{Notation}
To \changed{ease} notation, we write $i$ instead of $s_{i}$ for the standard generators of $\weyl$.
We will also often write e.g. $\mathtt{145}$ for \changed{the set} $\{1,4,5\}$.
\end{Notation}

Summarized in one picture:
\begin{gather}\label{Eq:Dynkin}
\tilde{A}_{\rank{-}1}=A^{(1)}_{\rank{-}1}\colon
\begin{tikzpicture}[scale=1.0,anchorbase]
\foreach \r [remember=\r as \rr] in {0,...,4} {
\node[circle,inner sep=1.8pt,fill=purple] (\r) at (360/7*\r:1){};
\node at (360/7*\r:1.3){$\r$};
\ifnum\r>0\draw(\rr)--(\r);\fi
}
\node[circle,inner sep=1.8pt,fill=purple] (6) at (360/7*6:1){};
\node at (360/7*6:1.3){$\rank{-}1$};
\draw[dashed](4) arc [start angle=205, end angle=290, radius=1] -- (6);
\draw(6)--(0);
\end{tikzpicture}
,
\quad
\text{relations}\colon
ii=1\text{ and }
\begin{gathered}
010=101
\\
121=212
\\
\dots
\\
(\rank{-}1)0(\rank{-}1)=0(\rank{-}1)0
\end{gathered}
.
\end{gather}
For $k\in\{0,\dots,\rank{-}1\}$\changed{, a \emph{$k$-color} is 
a subset $C\subset\nodes$ of size $k$}. For such a $C$ its \emph{ingredient colors} are the $D\subset\nodes$ of size $k-1$ such that $\#C\Delta D=1$, i.e. they differ by precisely one element. The \changed{unique $0$-color} is $\emptyset$, also called \emph{white}. 

\begin{Remark}
The colors correspond to \emph{finite parabolic subgroups.} 
\changed{Note that $\nodes$ is not a color, which could be called \emph{black}, as the corresponding subgroup would be $\weyl$ itself and therefore infinite (so black does not exist)}.
\end{Remark}

The \changed{$(\rank{-}1)$-colors} are the subsets 
$C\subset\nodes$ of size $\rank{-}1$, and since we use them often we call these \emph{top colors}. These correspond to maximal finite parabolic subgroups of $\weyl$.

\begin{Example}
For $\rank=3$, the color analogy can be used at its fullest, either in the RYB model (as in \cite{MaMaMiTu-trihedral}, the convention we follow) or the RGB model:
Here $\mathtt{0}$ is red, $\mathtt{1}$ is either yellow or green, and $\mathtt{2}$ is blue.
Moreover, $\mathtt{01}$ is orange or yellow, $\mathtt{02}$ is purple or magenta, and $\mathtt{12}$ is green or cyan.
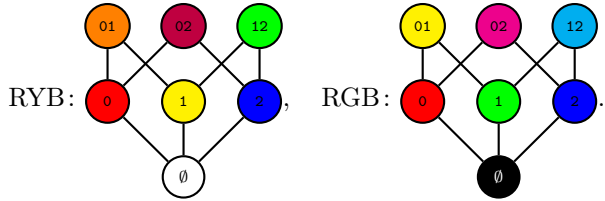
\begin{gather*}
\text{RYB}\colon
\begin{tikzpicture}[anchorbase,-,auto,node distance=1cm,
thick,main node/.style={circle,draw,font=\tiny}]
\node[main node,fill=white] (1) {$\emptyset$};
\node[main node,fill=yellow] (2) [above of=1] {\,$\mathtt{1}$\,};
\node[main node,fill=red] (3) [left of=2] {\,$\mathtt{0}$\,};
\node[main node,fill=blue] (4) [right of=2] {\,$\mathtt{2}$\,};
\node[main node,fill=purple] (5) [above of=2] {$\mathtt{02}$};
\node[main node,fill=orange] (6) [left of=5] {$\mathtt{01}$};
\node[main node,fill=green] (7) [right of=5] {$\mathtt{12}$};
\path[every node/.style={font=\sffamily\small}]
(1) edge (2)
edge (3)
edge (4)
(2) edge (6)
edge (7)
(3) edge (5)
edge (6)
(4) edge (5)
edge (7);
\end{tikzpicture}
,\quad
\text{RGB}\colon
\begin{tikzpicture}[anchorbase,-,auto,node distance=1cm,
thick,main node/.style={circle,draw,font=\tiny}]
\node[main node,fill=black] (1) {\color{white}$\emptyset$};
\node[main node,fill=green] (2) [above of=1] {\,$\mathtt{1}$\,};
\node[main node,fill=red] (3) [left of=2] {\,$\mathtt{0}$\,};
\node[main node,fill=blue] (4) [right of=2] {\,$\mathtt{2}$\,};
\node[main node,fill=magenta] (5) [above of=2] {$\mathtt{02}$};
\node[main node,fill=yellow] (6) [left of=5] {$\mathtt{01}$};
\node[main node,fill=cyan] (7) [right of=5] {$\mathtt{12}$};
\path[every node/.style={font=\sffamily\small}]
(1) edge (2)
edge (3)
edge (4)
(2) edge (6)
edge (7)
(3) edge (5)
edge (6)
(4) edge (5)
edge (7);
\end{tikzpicture}
.
\end{gather*}
(Note that white is black in RGB.)
Beyond $\rank=3$ the color analogy gets a bit shaky, but is still useful to keep in mind.
\end{Example} 

Label the vertices of an $\rank$gon from $0$ to $\rank{-}1$ counterclockwise.
Let us put the top colors on the vertices of such a regular $\rank$gon, so that $C$ corresponds to the vertex $i$ with $C\cup i=\nodes$. Let $\Z/\rank\Z\cong\langle\rho_{0},\dots,\rho_{\rank{-}1}|\rho_{i}\rho_{j}=\rho_{i+j(\bmod\rank)}\rangle$ act on this configuration by rotation, i.e. $\rho=\rho_{1}\colon i\to i+1(\bmod\rank)$. For example, for $\rank=4$ and $\rank=5$:
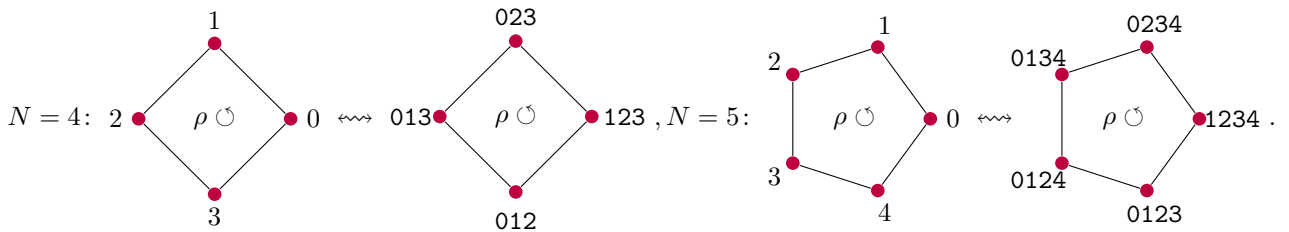
\begin{gather*}
\rank=4\colon
\begin{tikzpicture}[scale=1.0,anchorbase]
\node at (360/4*0:1.3){$0$};
\node at (360/4*1:1.3){$1$};
\node at (360/4*2:1.3){$2$};
\node at (360/4*3:1.3){$3$};
\node[circle,inner sep=1.8pt,fill=purple] (0) at (360/4*0:1){};
\node[circle,inner sep=1.8pt,fill=purple] (1) at (360/4*1:1){};
\node[circle,inner sep=1.8pt,fill=purple] (2) at (360/4*2:1){};
\node[circle,inner sep=1.8pt,fill=purple] (3) at (360/4*3:1){};
\draw(0)--(1);
\draw(1)--(2);
\draw(2)--(3);
\draw(3)--(0);
\node at (0,0) {$\rho\circlearrowleft$};
\end{tikzpicture}
\leftrightsquigarrow
\begin{tikzpicture}[scale=1.0,anchorbase]
\node at (360/4*0:1.38){$\;\mathtt{123}$};
\node at (360/4*1:1.32){$\mathtt{023}$};
\node at (360/4*2:1.38){$\mathtt{013}$};
\node at (360/4*3:1.38){$\mathtt{012}$};
\node[circle,inner sep=1.8pt,fill=purple] (0) at (360/4*0:1){};
\node[circle,inner sep=1.8pt,fill=purple] (1) at (360/4*1:1){};
\node[circle,inner sep=1.8pt,fill=purple] (2) at (360/4*2:1){};
\node[circle,inner sep=1.8pt,fill=purple] (3) at (360/4*3:1){};
\draw(0)--(1);
\draw(1)--(2);
\draw(2)--(3);
\draw(3)--(0);
\node at (0,0) {$\rho\circlearrowleft$};
\end{tikzpicture}
,
\rank=5\colon
\begin{tikzpicture}[scale=1.0,anchorbase]
\node at (360/5*0:1.3){$0$};
\node at (360/5*1:1.3){$1$};
\node at (360/5*2:1.3){$2$};
\node at (360/5*3:1.3){$3$};
\node at (360/5*4:1.3){$4$};
\node[circle,inner sep=1.8pt,fill=purple] (0) at (360/5*0:1){};
\node[circle,inner sep=1.8pt,fill=purple] (1) at (360/5*1:1){};
\node[circle,inner sep=1.8pt,fill=purple] (2) at (360/5*2:1){};
\node[circle,inner sep=1.8pt,fill=purple] (3) at (360/5*3:1){};
\node[circle,inner sep=1.8pt,fill=purple] (4) at (360/5*4:1){};
\draw(0)--(1);
\draw(1)--(2);
\draw(2)--(3);
\draw(3)--(4);
\draw(4)--(0);
\node at (0,0) {$\rho\circlearrowleft$};
\end{tikzpicture}
\leftrightsquigarrow\hspace*{-0.25cm}
\begin{tikzpicture}[scale=1.0,anchorbase]
\node at (360/5*0:1.38){$\;\mathtt{1234}$};
\node at (360/5*1:1.32){$\mathtt{0234}$};
\node at (360/5*2:1.38){$\mathtt{0134}$};
\node at (360/5*3:1.38){$\mathtt{0124}$};
\node at (360/5*4:1.32){$\mathtt{0123}$};
\node[circle,inner sep=1.8pt,fill=purple] (0) at (360/5*0:1){};
\node[circle,inner sep=1.8pt,fill=purple] (1) at (360/5*1:1){};
\node[circle,inner sep=1.8pt,fill=purple] (2) at (360/5*2:1){};
\node[circle,inner sep=1.8pt,fill=purple] (3) at (360/5*3:1){};
\node[circle,inner sep=1.8pt,fill=purple] (4) at (360/5*4:1){};
\draw(0)--(1);
\draw(1)--(2);
\draw(2)--(3);
\draw(3)--(4);
\draw(4)--(0);
\node at (0,0) {$\rho\circlearrowleft$};
\end{tikzpicture}
.
\end{gather*}

\begin{Notation}
We will always have a set of size $\rank-1$ inside a set of size $\rank$, and we will index the subsets by the missing value. Similarly, the parabolic subgroups of $\weyl$ corresponding to these subsets will be indexed by the missing value. For example, if $\rank=4$, then $\mathtt{023}$ will be denoted $\mathtt{1}$ and the corresponding parabolic subgroup by $\weyl_{\mathtt{1}}$. 
\end{Notation}

Recall that $\vpar$ is a generic parameter. The quantum numbers in this variable will be denoted by $\vnum{a}$ etc.

\begin{Definition}\label{D:Nhedral}
The \emph{Nhedral Hecke algebra $\nhecke$ of level $\infty$} is 
the associative unital ($\Cv$-)algebra generated by 
$\rank$ elements $\{\kl_{i}|i\in\nodes\}$
subject to the following relations.
\begin{gather}\label{Eq:fundamental-relations}
\kl_{i}^{2}=\vnum{\rank}!\cdot\kl_{i}, \quad
\kl_{k+i+j}\kl_{k+i}\kl_{k}
=
\kl_{k+i+j}\kl_{k+j}\kl_{k}
\text{ for all }
i,j,k\in\nodes,
\end{gather}
where indices are taken modulo $\rank$. The second type of relation is called \emph{fundamental commutativity}.
\end{Definition}

\begin{Remark}\label{R:Analogy}
\changedd{We have the following rough analogy, motivating \autoref{D:Nhedral}.
This algebra sits inside the affine Hecke algebra (see \autoref{L:Embedding}), capturing the part associated with the bottom \ochangeddd{Kazhdan--Lusztig} cell. Through decategorified geometric Satake, tensoring with a \changeddd{simple} representation of color $i$\changeddd{, call it $\Ll_{\bm}$}, acts on representations of color $j$, producing representations of color $i+j$. This operation corresponds to \changeddd{multiplying with} a \ochangeddd{Kazhdan--Lusztig} basis element $\Theta_{w}$ for \changeddd{a particular maximal double coset in $W_{i+j}\backslash W/W_{j}$ depending on $\bm$ and $j$}. As a result, a single representation gives rise to different elements in the Hecke algebra depending on the starting color.}

\changedd{This analogy is only approximate: multiplication in the Hecke algebra corresponds to composition of functors \changeddd{of tensoring with the corresponding representation} only up to scalar\changeddd{. (For the reader familiar with \cite{MaMaMiTu-trihedral}, for $\rank=2$ and $\rank=3$ there is a diagrammatic incarnation of this functor that changes the color of faces.)} For example, the trivial representation maps to an idempotent $\theta_{i}$, but $\theta_{i}^{2}=\vnum{\rank}!\cdot\theta_{i}$. This mismatch arises because geometric Satake more precisely lands in the spherical affine Hecke algebra.}
\end{Remark}

\changedd{Recall that $\fu_{i}$ denotes the isomorphism class of the representation $\Ll_{\omega_{i}}$ of $\slncat$.} The fundamental commutativity is then mimicking 
$\fu_{j}\fu_{i}=\fu_{i}\fu_{j}$ \changedd{when compared with \autoref{R:Analogy}}.
Moreover, given an expression $\kl_{j}\dots\kl_{i}$, we call $i$ the \emph{starting} and $j$ the \emph{ending color}.

\begin{Example}
Let $\rank=4$. Then, up to changing the starting color, we 
have
\begin{gather*}
\kl_{3}\kl_{1}\kl_{0}
=
\kl_{3}\kl_{2}\kl_{0}
\leftrightsquigarrow
\fu_{2}\fu_{1}=\fu_{1}\fu_{2},
\\
\kl_{0}\kl_{1}\kl_{0}
=
\kl_{0}\kl_{3}\kl_{0}
\leftrightsquigarrow
\fu_{3}\fu_{1}=\fu_{1}\fu_{3},
\\
\kl_{1}\kl_{2}\kl_{0}
=
\kl_{1}\kl_{3}\kl_{0}
\leftrightsquigarrow
\fu_{3}\fu_{2}=\fu_{2}\fu_{3},
\end{gather*}
as the fundamental commutativity relations.
\end{Example}

\begin{Remark}\label{R:shift}
\changed{Recall from \autoref{SS:IntroA} that we think of Nhedral Hecke algebras as having a degree zero part corresponding to (the Grothendieck ring of) $\verg{\levelthing}{\rank}$. Essentially, this means that there is a scaling of the Nhedral Hecke algebras such that: a) specializing $\vpar$ to zero results in a coefficient-free structure; b) the relations are simplified; and c) the combinatorics of words in these algebras mirror the behavior of the fundamental variables $\fu_{i}$. We will explore all of this below.}

\ochangedd{Note that $\mathbf{a}=\frac{1}{2}\rank(\rank-1)$ is half of the span of $\vnum{\rank}!$. So after \changed{shifting by $\vpar^{\mathbf{a}}$}, at $\vpar=0$ the relation $\kl_{i}^{2}=\vnum{\rank}!\cdot\kl_{i}$ becomes $\kl_{i}^{2}=\kl_{i}$. We will see a natural explanation for the value of $\mathbf{a}=\frac{1}{2}\rank(\rank-1)$ in \autoref{R:AValue}.}
\end{Remark}

Every $i$ defines a word of length $\frac{1}{2}\rank(\rank{-}1)$ as follows.
Say $i=\rank{-}1$, the general case being similar, then
\begin{gather*}
\changed{w_{i}=\big(012\dots(\rank-2)\big)\dots(012)(01)(0)\in\weyl_{i}.}
\end{gather*}
Note that the $w_{i}$ are \changed{reduced expressions for} the longest elements 
in their associated parabolic subgroups $\weyl_{i}$.

Let $\hecke$ denote the Hecke algebra of affine type $A_{\rank{-}1}$. Following the conventions of
\cite{So-tilting-a}, this algebra is generated by 
$\{\kltwo_{i}|i\in\nodes\}$ subject to (reading indices modulo $\rank$):
\begin{gather*}
\kltwo_{i}^{2}=\vnum{2}\kltwo_{i},
\quad
\kltwo_{i}\kltwo_{j}\kltwo_{i}-\kltwo_{i}=\kltwo_{j}\kltwo_{i}\kltwo_{j}-\kltwo_{j}
\text{ for }|i-j|=1,
\quad
\kltwo_{i}\kltwo_{j}=\kltwo_{j}\kltwo_{i}
\text{ for }|i-j|>1.
\end{gather*}
Let $\kltwo_{w}$ denote the KL basis element for $w\in\weyl$, see e.g. \cite{So-tilting-a}.
We have:

\begin{Lemma}\label{L:Embedding}
(This assumes \autoref{QSH}.) The algebra homomorphism given by 
\begin{gather*}
\kl_{i}\mapsto\kltwo_{w_{i}},
\end{gather*}
defines an embedding $\nhecke\hookrightarrow\hecke$ of algebras.
\end{Lemma}

The proof of \autoref{L:Embedding} is different from the one given in \cite[Lemma 3.2]{MaMaMiTu-trihedral}, \changedd{and is postponed to \autoref{Pf:Embedding}.}

\begin{Remark}\label{R:GeoSatake2}
\changedd{Returning to \autoref{R:GeoSatake}, these statements are primarily uncategorified and should follow from combinatorial equivalences akin to the geometric Satake equivalence. However, we remain uncertain about how to formalize this, and after consulting several experts, none of us could interpret, for example, the fundamental commutativity relation geometrically. We identify three immediate issues. First, the combinatorial statement that most closely resembles what we require, namely, that certain Kazhdan--Lusztig basis elements multiply like representations, applies to the spherical Hecke algebra, whereas our work involves the affine Hecke algebra (though this is a relatively minor issue). Second, the fundamental commutativity relation in our context involves three elements, whereas the usual Satake equivalence concerns the commutativity of the tensor product, which involves only two factors. As such, the fundamental commutativity is more appropriately viewed as a 2-categorical analog, a perspective that is directly incorporated into \autoref{QSH}. Finally, \autoref{QSH} is specifically concerned with type A phenomena, while the geometric Satake correspondence is not type-specific, leading to subtle but crucial differences. Ultimately, we remain unsure how to eliminate the assumption \autoref{QSH} from these results.}
\end{Remark}

The Nhedral KL combinatorics works as follows. For $\bk=(k_{1},\dots,k_{\rank-1})\in X^{+}$ and starting color $i$, let 
\begin{gather*}
\rkl{\bk}{i}=\kl_{i_{\bsum\bk}}\dots\kl_{i_{1}}\kl_{i_{0}}
\end{gather*}
with $i_{0}=i$ and $i_{r+1}=\rho^{j}(i_{r})=i_{r}+j$ \changed{such that $k_{j}=\#\{0 \leq r < \bsum\bk\mid i_{r+1}=\rho^{j}(i_{r})\}$ for all $j$. (This looks cumbersome, but is not difficult to explain with an example, see \autoref{E:hkl} below.)} 
We will call $\bsum\bk$ the ending color.
Reversing the order, let
\begin{gather*}
\lkl{\bk}{i}=\kl_{i_{0}}\kl_{i_{1}}\dots\kl_{i_{\bsum\bk}}
\end{gather*}
be defined similarly as $\rkl{\bk}{i}$ but using 
$\rho^{-j}$ instead of $\rho^{j}$.

\begin{Lemma}\label{L:well-defined-paths}
For \ochanged{any color} $i\in \nodes$ and any $\bk\in X^{+}$, the \ochanged{elements} $\rkl{\bk}{i}$ and $\lkl{\bk}{i}$ only \changeddd{depend} on $\bk$ and not on the chosen sequence $i=i_{0},i_{1},\dots,i_{\bsum\bk}$.
\end{Lemma}

\begin{proof}
\changed{Using the fundamental commutativity relation, we can swap two successive differences of colors. Therefore, any word representing $\rkl{\bk}{i}$ is equivalent to the word with increasing successive differences of colors.}
\end{proof}

\begin{Lemma}\label{L:samebasis_h}
\changedd{For any $i\in \nodes$ and $\bk\in X^{+}$, let $j$ be the ending color of $\rkl{\bk}{i}$. We have $\rkl{\bk}{i}=\lkl{\bk}{j}$.}
\end{Lemma}

\begin{proof}
\changedd{This follows immediately from \autoref{L:well-defined-paths}.}
\end{proof}

\begin{Remark}
\label{R:h-element}
Via \autoref{QSH}, the element $\rkl{\bk}{i}$ is associated to $\fu_{1}^{k_{1}}\dots\fu_{\rank-1}^{k_{\rank-1}}$ because its definition involves $k_{j}$ times the application of $\rho^{j}$. 
\end{Remark}

\begin{Example}\label{E:hkl}
Let us choose $0$ as the starting color. For $\rank=4$ and $\bsum\bk=2$, we have
\begin{gather*}
\rkl{2,0,0}{0} = \kl_{2}\kl_{1}\kl_{0} \leftrightsquigarrow \fu_{1}^{2},\quad
\rkl{0,2,0}{0} = \kl_{0}\kl_{2}\kl_{0} \leftrightsquigarrow \fu_{2}^{2},\quad
\rkl{0,0,2}{0} = \kl_{2}\kl_{3}\kl_{0} \leftrightsquigarrow \fu_{3}^{2},\\
\rkl{1,1,0}{0} = \kl_{3}\kl_{1}\kl_{0} \leftrightsquigarrow \fu_{2}\fu_{1},\quad
\rkl{1,0,1}{0} = \kl_{0}\kl_{1}\kl_{0} \leftrightsquigarrow \fu_{3}\fu_{1},\quad
\rkl{0,1,1}{0} = \kl_{1}\kl_{2}\kl_{0} \leftrightsquigarrow \fu_{3}\fu_{2}.
\end{gather*}
Note that, for example, $\rkl{1,1,0}{0} = \kl_{3}\kl_{2}\kl_{0}$ by \autoref{Eq:fundamental-relations}, and on the representation side by $\fu_{2}\fu_{1}=\fu_{1}\fu_{2}$. \changed{With $1$ as the starting color, we have $\rkl{2,0,0}{1} = \kl_{3}\kl_{2}\kl_{1}$.}
\end{Example}

\changed{Let} $\scalar=\vnum{\rank-1}!$ and \changed{recall the numbers} $d^{\bk}_{\bm}$ from
\autoref{Eq:dNumbers}. \changed{Motivated by \autoref{Eq:dNumbers}, f}or each $\bm\in X^{+}$, we define 
\emph{(right) Nhedral KL basis elements}:
\begin{gather*}
\RKL{\bm}{i}=\sum_{\bk}\scalar^{-\bsum\bk}d^{\bk}_{\bm}\cdot\rkl{\bk}{i}.
\end{gather*}
Note that the three sums are finite, 
because $d^{\bk}_{\bm}=0$ unless $\bsum\bk\leq\bsum\bm$. Similarly, using $\lkl{\bk}{i}$ instead of $\rkl{\bk}{i}$, we define the (left) Nhedral KL basis elements $\LKL{\bm}{i}$. Note that the ending color of every term $\rkl{\bk}{i}$ appearing in the sum have the same ending color.

\begin{Remark}
\label{R:KL-basis-Chebyshev}
Via \autoref{QSH}, the element $\RKL{\bm}{i}$ is associated to the Chebyshev polynomial $\pxy{\bm}(\fu_{1},\dots,\fu_{\rank-1})$.
\end{Remark}

\begin{Example}
For $\rank=4$ and $\bsum\bk=2$, we have
\begin{gather*}
\RKL{2,0,0}{0} = \scalar^{-2}\kl_{2}\kl_{1}\kl_{0}-\scalar^{-1}\kl_{2}\kl_{0} \leftrightsquigarrow \pxy{2,0,0}=\fu_{1}^{2}-\fu_{2},\ 
\RKL{0,2,0}{0} = \scalar^{-2}\kl_{0}\kl_{2}\kl_{0}-\scalar^{-2}\kl_{0}\kl_{1}\kl_{0} \leftrightsquigarrow \pxy{0,2,0}=\fu_{2}^{2}-\fu_{3}\fu_{1},\\
\RKL{0,0,2}{0} = \scalar^{-2}\kl_{2}\kl_{3}\kl_{0}-\scalar^{-1}\kl_{2}\kl_{0} \leftrightsquigarrow \pxy{0,0,2}=\fu_{3}^{2}-\fu_{2},\ 
\RKL{1,1,0}{0} = \scalar^{-2}\kl_{3}\kl_{1}\kl_{0}-\scalar^{-1}\kl_{3}\kl_{0} \leftrightsquigarrow \pxy{1,1,0}=\fu_{2}\fu_{1}-\fu_{3},\\
\RKL{1,0,1}{0} = \scalar^{-2}\kl_{0}\kl_{1}\kl_{0}-\kl_{0} \leftrightsquigarrow \pxy{1,0,1}=\fu_{3}\fu_{1}-1,\ 
\RKL{0,1,1}{0} = \scalar^{-2}\kl_{1}\kl_{2}\kl_{0}-\scalar^{-1}\kl_{1}\kl_{0} \leftrightsquigarrow \pxy{0,1,1}=\fu_{3}\fu_{2}-\fu_{1},
\end{gather*}
with $0$ as the starting color.
\end{Example}

\begin{Lemma}\label{L:samebasis_c}
\changedd{For any $i\in \nodes$ and $\bk\in X^{+}$, let $j$ be the ending color of $\RKL{\bk}{i}$. We have $\RKL{\bk}{i}=\LKL{\bk}{j}$.}
\end{Lemma}

\begin{proof}
\changedd{Same as \autoref{L:samebasis_h}.}
\end{proof}

\begin{Lemma}\label{L:KLaction}(This assumes \autoref{QSH})
For all $i\in \nodes$ and $\bm\in X^{+}$, let $j$ be the ending color of $\RKL{\bm}{i}$. Let $0\leq k < \rank$. With the notation in \autoref{L:ChebyshevRecursion}, we have
\begin{gather*}
\kl_{j+k}\RKL{\bm}{i} = 
\begin{cases}
\vnum{\rank}!\RKL{\bm}{i} & k=0,\\
\scalar\sum_{l}\RKL{\bm+w_{l}^{k}}{i} & \text{otherwise},
\end{cases}
\end{gather*}
where terms with negative entries are zero. Similarly for the left KL elements.
\end{Lemma}

\begin{Remark}
Coming back to \autoref{R:shift}, \changed{after shifting,} at $\vpar=0$, \autoref{L:KLaction} becomes
\begin{gather*}
\kl_{j+k}\RKL{\bm}{i} = 
\begin{cases}
\RKL{\bm}{i} & k=0,\\
0 & \text{otherwise},
\end{cases}
\end{gather*}
since $\scalar$ has a smaller span than $\rank(\rank-1)$. This justifies the ``degree $0$ isomorphism'' from $\nhecke[\level]$ in \autoref{D:NHedral-level} to the Grothendieck ring of the asymptotic category $\verg{\levelthing}{\rank}$ in \autoref{D:AsymptoticCat}.
\end{Remark}

\begin{Proposition}\label{P:KLbasisInfinity}(This assumes \autoref{QSH})
Each of the four sets 
\begin{gather*}
\basisH=
\{1\}
\cup
\{\rkl{\bk}{i}|\bk\in X^+,\, i\in \nodes\},
\quad
\Hbasis=
\{1\}\cup
\{\lkl{\bk}{i}|\bk\in X^+,\, i\in \nodes\},
\\
\basisC=
\{1\}
\cup\{\RKL{\bm}{i}|\bm\in X^+,\, i\in \nodes\}
\quad
\Cbasis=
\{1\}\cup
\{\LKL{\bm}{i}|\bm\in X^+,\, i\in \nodes\}
\end{gather*} 
is a basis of $\nhecke$. The first two are called \emph{Nhedral Bott--Samelson bases}, the final two \emph{Nhedral KL bases}.
\end{Proposition}

Recall the definition of left, right and two-sided cells for $\nhecke$ using the Nhedral KL basis, see for example \cite[Definition 3.10]{MaMaMiTu-trihedral}. The unit forms its own cell, that we call the trivial cell.

\begin{Proposition}\label{P:Cells}
\changed{There is one nontrivial two-sided cell for the algebra $\nhecke$, namely
\begin{gather*}
\tcell=
\left\{\RKL{\bm}{i}|\bm\in X^+,i\in \nodes\right\}
=
\left\{\LKL{\bm}{i}|\bm\in X^+,i\in \nodes\right\}.
\end{gather*}
The left and right cells contained in $J$ are
\begin{gather*}
\lcell_{i}=\left\{\RKL{\bm}{i}|\bm\in X^+\right\},
\quad\quad
{}_{i}\rcell=\left\{\LKL{\bm}{i}|\bm\in X^+\right\},
\quad\quad \hbox{ for $i\in \nodes$.}
\end{gather*}
Therefore, there are $N$ nontrivial left and right cells.
}
\end{Proposition}

\begin{proof}
As in \cite[Proof of Proposition 3.11]{MaMaMiTu-trihedral}, this follows from \autoref{L:KLaction}.
\end{proof}

We now define 
finite dimensional quotients of $\nhecke$, which are compatible 
with the cell structure in \autoref{P:Cells}.

\begin{Definition}\label{D:NHedral-level}
For our fixed $\level$, let $\killideal{\level}$ 
be the two-sided ideal in $\nhecke$ generated by
\begin{gather*}
\{\RKL{\bm}{i}\mid i\in \nodes, \bm\in X^{+},\bsum\bm = \level+1\}.
\end{gather*}
We define the \emph{Nhedral Hecke algebra of level $\level$} as
\begin{gather*}
\nhecke[\level]=\nhecke/\killideal{\level}
\end{gather*}
and we call $\killideal{\level}$ the \emph{vanishing ideal} of level $\level$.
\end{Definition}

\begin{Proposition}\label{P:KLbasisLevel}
The set 
\begin{gather*}
\basisC[\level]=
\{1\}\cup
\{\RKL{\bm}{i}\mid i\in \nodes, \bm\in X^{+}(\level)\}
\stackrel{(*)}{=}
\Cbasis[\level]=
\{1\}\cup
\{\LKL{\bm}{i}\mid i\in \nodes, \bm\in X^{+}(\level)\}
\end{gather*}
(the equality (*) follows from \autoref{L:samebasis_c})
is a basis of $\nhecke[\level]$. 
Thus, we have 
$\dim_{\Cv}\nhecke[\level]=1+\rank\tnumbers{\rank}{\level}$.
\end{Proposition}

\begin{proof}
Similarly as \cite[Proof of Proposition 3.14]{MaMaMiTu-trihedral}.
\end{proof}

\begin{Proposition}
The nontrivial cells for the algebra $\nhecke[\level]$ are
\begin{gather*}
\lcell_{i}=\left\{\RKL{\bm}{i}|\bm\in X^{+}(\level)\right\},
\quad\quad
{}_{i}\rcell=\left\{\LKL{\bm}{i}|\bm\in X^{+}(\level)\right\},
\quad\quad \hbox{ for $i\in \nodes$,}
\\
\tcell=
\left\{\RKL{\bm}{i}|\bm\in X^{+}(\level),i\in \nodes\right\}
=
\left\{\LKL{\bm}{i}|\bm\in X^{+}(\level),i\in \nodes\right\},
\end{gather*}
where $\lcell_{i}$, ${}_{i}\rcell$ and $\tcell$ 
are left, right and two-sided cells
respectively. In particular, each left and right cell is of size $\tnumbers{\rank}{\level}$, and $\tcell$ is of size $\rank\tnumbers{\rank}{\level}$.
\end{Proposition}

\begin{proof}
This follows from the previous results.
\end{proof}

\begin{Example}
Left and right cells correspond to the cut-off of the positive Weyl chamber of type $A_{\rank-1}$, similarly as in \cite[Example 3.16]{MaMaMiTu-trihedral}. 
\end{Example}

\subsection{Nhedral complex representations}\label{SS:NhedralComplex}

We now classify all simple representations of $\nhecke[\level]$ on $\Cv$-vector spaces, and this classification implies that 
$\nhecke[\level]$ is semisimple.

To this end, for $\bl=(\lambda_{i})_{i\in I}\in(\Cv)^{\nodes}$ define
\begin{gather*}
\M_{\bl}
\colon\Cv\langle\kl_{i}|i\in\nodes\rangle\to\Cv,
\kl_{i}\mapsto\lambda_{i}, 
\end{gather*}
which determines a one dimensional
$\nhecke[\level]$-representation in the following cases. \changed{For $x\in\Cv$, w}e let $(x,i)\in(\Cv)^{\nodes}$ the element with the $i$th entry $x$ and zero otherwise.

\begin{Proposition}\label{P:OneDimensional}
The following table 
\begin{gather}\label{eq:the-one-dims}
\begin{tikzpicture}[baseline=(current bounding box.center),yscale=0.6]
\matrix (m) [matrix of math nodes, row sep={.85cm,between origins}, column
sep={4.0cm,between origins}, text height=1.5ex, text depth=0.25ex, ampersand replacement=\&] {
\level\equiv 0\bmod\rank \&  \level\not\equiv 0\bmod \rank \\
\changedd{\M_{(\vnum{\rank}!,i)}}\text{ for }i\in\nodes\changedd{\text{, and } \M_{0,\dots,0}}\& \M_{0,\dots,0}  \\
\text{$\rank+1$ in total} \& \text{only one} \\};
\draw[densely dashed] ($(m-1-1.south west)+ (-1.2,0)$) to (m-1-2.south east);
\draw[densely dashed] ($(m-1-2.north west) + (-0.3,0)$) to ($(m-1-2.north west) + (-0.3,-3.75)$);
\draw[densely dashed] ($(m-1-1.south west)+ (-1.2,-1.8)$) to ($(m-1-2.south east)+(0,-1.8)$);
\end{tikzpicture}
\end{gather}
gives a complete and irredundant list of one dimensional 
$\nhecke[\level]$-representations.
\end{Proposition}

\begin{proof}
That $\M_{\changed{0,\ldots,0}}$ is well-defined and simple is immediate. The remaining parts follow from 
\autoref{T:simples} below.
\end{proof}

We now define some representations of dimension $\rank$, which are parametrized by the points in $\vanpar{\level}$. Given $\bs\in E$, we define
\begin{gather*}
\matrep{0}{\bs}
=\scalar
\begin{psmallmatrix}
\vnum{\rank} & \Zfun_{1}(\bs) & \Zfun_{2}(\bs) & \dots & \Zfun_{\rank-1}(\bs) \\
0 & 0 & 0 & \dots & 0\\
\vdots & \vdots & \vdots & \ddots & \vdots \\
0 & 0 & 0 & \dots & 0
\end{psmallmatrix},
\quad
\matrep{1}{\bs}
=\scalar
\begin{psmallmatrix}
0 & 0 & 0 & \dots & 0\\
\Zfun_{\rank-1}(\bs) & \vnum{\rank} & \Zfun_{1}(\bs) & \dots & \Zfun_{\rank-2}(\bs) \\
0 & 0 & 0 & \dots & 0\\
\vdots & \vdots & \vdots & \ddots & \vdots \\
0 & 0 & 0 & \dots & 0
\end{psmallmatrix},
\dots,\\
\matrep{\rank-1}{\bs}
=\scalar
\begin{psmallmatrix}
0 & 0 & \dots & 0 & 0\\
\vdots & \vdots & \ddots & \vdots & \vdots \\
0 & 0 & \dots & 0& 0\\
\Zfun_{1}(\bs) & \Zfun_{2}(\bs) & \dots & \Zfun_{\rank-1}(\bs) & \vnum{\rank}
\end{psmallmatrix}.
\end{gather*}

\begin{Lemma}\label{L:Tinfinityrep}
The assignment $\kl_{i}\mapsto \matrep{i}{\bs}$ is a well-defined representation of $\nhecke[\level]$ if and only if $\bs\in\vanpar{\level}$.
\end{Lemma}

\begin{proof}
The defining relations \autoref{Eq:fundamental-relations} of $\nhecke$ are easy to check, since the matrices have only one nonzero row. It factors through $\vanideal{\level}$ if and only if $(\Zfun_{1}(\bs),\dots,\Zfun_{\rank}(\bs))\in\vanset{\level}$ as in \cite[Lemma 3.18]{MaMaMiTu-trihedral}.
\end{proof}

For $\sigma\in\vanpar{\level}$, denote by $\M(\bs)$ the above defined $\nhecke[\level]$-representation.

\begin{Theorem}\label{T:simples}
We have the following.
\begin{enumerate}
\item For $\bs\in\vanpar{\level}$, the representation $\M(\bs)$ decomposes as the direct sum of $m$ equidimensional nonisomorphic simple representations, where $m$ is the order of the stabilizer of $\bs$ for the $\Z/\rank\Z$-action.
\item For $\bs,\bs'\in\vanpar{\level}$, we have
\begin{gather*}
\dim\ \Hom_{\nhecke[\level]}\big(\M(\bs),\M(\bs')\big) = 
\begin{cases}
0 & \text{if } \bs \text{ and } \bs' \text{ are in different orbits},\\
m & \text{otherwise}.
\end{cases}
\end{gather*}
\item The simple representations appearing in (a), when $\bs$ varying over the orbits of $\vanpar{\level}$, and the representation $\M_{(0,\dots,0)}$ form an complete irredundant list of simple representations.
\item The algebra $\nhecke[\level]$ is semisimple.
\end{enumerate}
\end{Theorem}

\begin{proof}
\changedd{\textit{Part (a).} This is proven in \autoref{Pf:simples-a}.}

\changedd{\textit{Part (b).} This is proven in \autoref{Pf:simples-b}.}

\changedd{\textit{Parts (c) and (d).} These are proven in \autoref{Pf:simples-cd}.}
\end{proof}

\subsection{Nhedral integral representations}\label{SS:NhedralIntegral}

Let $\Nv=\N[\vpar,\vpar^{-1}]$. Recall that an $\Nv$-algebra is a $\Cv$-algebra with a fixed basis such that the structure constants for the multiplication are in $\Nv$. An $\Nv$-representation of such an algebra is a representation over $\Cv$ with a fixed basis such that the structure constants for the action are in $\Nv$.

\begin{Lemma}(This assumes \autoref{QSH})\label{L:integrality}
The Nhedral Hecke algebras $\nhecke$ and $\nhecke[\level]$ are $\Nv$-algebras with respect to the KL basis of \autoref{P:KLbasisInfinity} and \autoref{P:KLbasisLevel}.
\end{Lemma}

We now study $\Nv$-representations of the Nhedral Hecke algebras. To this end, we use the following type of colored graphs. Recall $\nodes=\{0,\dots,\rank-1\}$. 

\begin{Definition}
An unoriented graph $\Gamma=(\vertices,\edges)$ is Ncolored if $\vertices$ is partitioned into $\rank$ disjoint sets $\vertices=\coprod_{i\in \nodes}\vertices_{i}$ such that neighboring vertices have different colors. The vertices in $\vertices_{i}$ are said to be of color $i$.
\end{Definition}

The adjacency matrix of such a graph can be brought in the form
\begin{gather*}
A(\Gamma)=
\begin{pNiceMatrix}[first-row,first-col,hvlines]
& 0 & 1 & \cdots & \rank-1 \\
0 & 0 & Z_{0}^{1} & \dots & Z_{0}^{\rank-1}  \\
1 & Z_{1}^{0} & 0 & \dots & Z_{1}^{\rank-1}\\
\vdots & \vdots & \vdots  & \ddots  & \vdots\\
\rank-1 & Z_{\rank-1}^{0} & Z_{\rank-1}^{1} & \dots & 0 
\end{pNiceMatrix}
\quad\text{ with }\quad (Z_{i}^{j})^{\mathrm{T}} = Z_{j}^{i}.
\end{gather*}

Given an Ncolored graph $\Gamma=(\vertices,\edges)$, define
\begin{gather*}
\M(\Gamma)
\colon\Cv\langle\kl_{i}|i\in\nodes\rangle\to\End_{\Cv}(\Cv\vertices),
\kl_{i}\mapsto\matrep{i}{\Gamma}, 
\end{gather*}
where the matrices are
\begin{gather*}
\matrep{0}{\Gamma}
=\scalar
\begin{psmallmatrix}
\vnum{\rank}\idmor & Z_{0}^{1} & Z_{0}^{2} & \dots & Z_{0}^{\rank-1} \\
0 & 0 & 0 & \dots & 0\\
\vdots & \vdots & \vdots & \ddots & \vdots \\
0 & 0 & 0 & \dots & 0
\end{psmallmatrix},
\quad
\matrep{1}{\Gamma}
=\scalar
\begin{psmallmatrix}
0 & 0 & 0 & \dots & 0\\
Z_{1}^{0} & \vnum{\rank}\idmor & Z_{1}^{2} & \dots & Z_{1}^{\rank-1} \\
0 & 0 & 0 & \dots & 0\\
\vdots & \vdots & \vdots & \ddots & \vdots \\
0 & 0 & 0 & \dots & 0
\end{psmallmatrix},
\dots,\\
\matrep{\rank-1}{\Gamma}
=\scalar
\begin{psmallmatrix}
0 & 0 & \dots & 0 & 0\\
\vdots & \vdots & \ddots & \vdots & \vdots \\
0 & 0 & \dots & 0& 0\\
Z_{\rank-1}^{0} & Z_{\rank-1}^{1} & \dots & Z_{\rank-1}^{\rank-2} & \vnum{\rank}\idmor
\end{psmallmatrix}.
\end{gather*}

Let $Z_{i}^{i}=\idmor$ for $i\in\nodes$.

\begin{Lemma}\label{L:integralTinfinty}
The assignment $\theta_{i}\mapsto \matrep{i}{\Gamma}$ is a well-defined representation of $\nhecke$ if and only if for all $i,j,\mathtt{k}\in\nodes$, $Z_{i+j+k}^{i+j}Z_{\mathtt{i+j}}^{i}=Z_{i+j+k}^{i+k}Z_{i+k}^{i}$.
\end{Lemma}

\begin{proof}
This follows from an easy direct computation (since only one row of the action matrices is nonzero).
\end{proof}

For $i\in\{1,\dots,\rank-1\}$ define the following oriented subgraphs $\Gamma_{i}$ of $\Gamma$ with adjacency matrix $A(\Gamma_{i})$ obtained from $A(\Gamma)$ by setting all blocks to zero except the blocks $Z_{j}^{i+j}$. Note that the condition in \autoref{L:integralTinfinty} is equivalent to the pairwise commutation of the matrices $A(\Gamma_{1}),\dots,A(\Gamma_{\rank-1})$.

\begin{Example}
For $\rank=4$ we have:
\begin{gather*}
A(\Gamma)=\begin{pNiceMatrix}[hvlines]
0 & Z_{0}^{1} & Z_{0}^{2} & Z_{0}^{3}  \\
Z_{1}^{0} & 0 & Z_{1}^{2} & Z_{1}^{3}\\
Z_{2}^{0} & Z_{2}^{1} & 0 & Z_{2}^{3}\\
Z_{3}^{0} & Z_{3}^{1} & Z_{3}^{2} & 0 
\end{pNiceMatrix},\;
A(\Gamma_{1})=A(\Gamma_{3})^{\mathrm{T}} =
\begin{pNiceMatrix}[hvlines]
0 & Z_{0}^{1} & 0 & 0\\
0 & 0 & Z_{1}^{2} & 0\\
0 & 0 & 0 & Z_{2}^{3}\\
Z_{3}^{0} & 0 & 0 & 0 
\end{pNiceMatrix},\;
A(\Gamma_{2})=\begin{pNiceMatrix}[hvlines]
0 & 0 & Z_{0}^{2} & 0  \\
0 & 0 & 0 & Z_{1}^{3}\\
Z_{2}^{0} & 0 & 0 & 0\\
0 & Z_{3}^{1} & 0 & 0 
\end{pNiceMatrix}.
\end{gather*}
If $\Gamma$ is the graph of the type as in \autoref{E:tetrahedron}, then the $\Gamma_{i}$ are obtained by putting an orientation on some edges, and removing other edges. Generically, the picture is: 
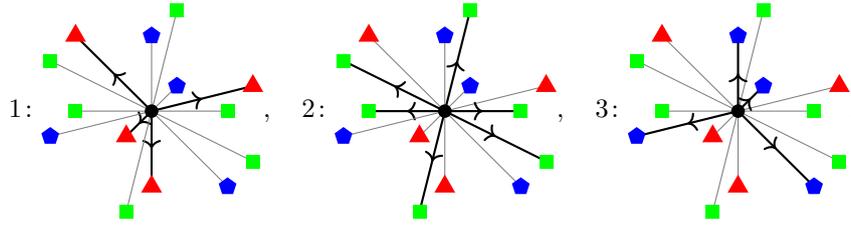
\begin{gather*}
\input{figures/genericsubgraph}
\end{gather*}
These three graphs correspond to tensoring by $\Ll_{\omega_{i}}$. Note that the generic picture for $1$ and $3$ have four edges, while $2$ has six edges, matching the dimensions of the $\Ll_{\omega_{i}}$.
\end{Example}

Let us now assume that $\Gamma$ is such that 
the condition in \autoref{L:integralTinfinty} holds.

\begin{Lemma}\label{L:integralTlevel}
The $\nhecke$-representation $\M(\Gamma)$ \changed{descends} to a $\nhecke[\level]$-representation
if and only if
\begin{gather*}
\pxy{\bm}\big(A(\Gamma_{1}),\dots,
A(\Gamma_{\rank-1})
\big)=0\text{ for all }\bm\text{ with }\bsum\bm=\level+1.
\end{gather*}
\end{Lemma}

\begin{proof}
By construction.
\end{proof}

Some examples of Ncolored graphs $\Gamma$ satisfying \autoref{L:integralTlevel} are obtained through the fusion rules of $\slncatr$. Here is a (non-exhaustive) list:

\begin{Definition}
We define the \emph{graph of type A} of rank $\rank$ and level $\level$ as the graph $\Gamma$ with set of vertices $X^{+}(\level)$, colored using the color $\chi_{c}$, and the vertices $\bm$ and $\bk$ are adjacent if and only if $\Ll_{\bm}$ is a summand of $\bigoplus_{i=1}^{\rank-1}\Ll_{\omega_{i}}\otimes\Ll_{\bk}$. The graph $\Gamma_{i}$ is then the fusion graph of the object $\Ll_{\omega_{i}}$.

Moreover, we define \emph{graphs of type D} in \autoref{Pf:typeD-graphs}. 

For $\rank=4$,
we have additional graphs that we will not describe here but rather refer to \cite[Figures 3 and 4, denoted $2A^{c}_{\level}$, $2(A^{c}_{\level}/2)$ and $E$]{Oc-classification-sun}.
These are called \emph{graphs of conjugate type $A$} (these two are infinite family), 
and \emph{graphs of type E} (there are six of these).
\end{Definition}

\begin{Remark}
There appears to be a small typo in \cite[Figures 3 and 4]{Oc-classification-sun} for the graphs labeled $2A^{c}_{\level}$ and $2(A^{c}_{\level}/2)$: all double edges should be colored blue instead of red (in 2024, we got a colored version of Ocneanu's paper from \url{https://cel.hal.science/cel-00374414/document}).
\end{Remark}

\begin{Example}[Type A]
For $\rank=4$ and $\level=4$, the graph of type $A$ is the tetrahedron of \autoref{E:tetrahedron} and the plot of the eigenvalues of $A(\Gamma_{1})$ are:
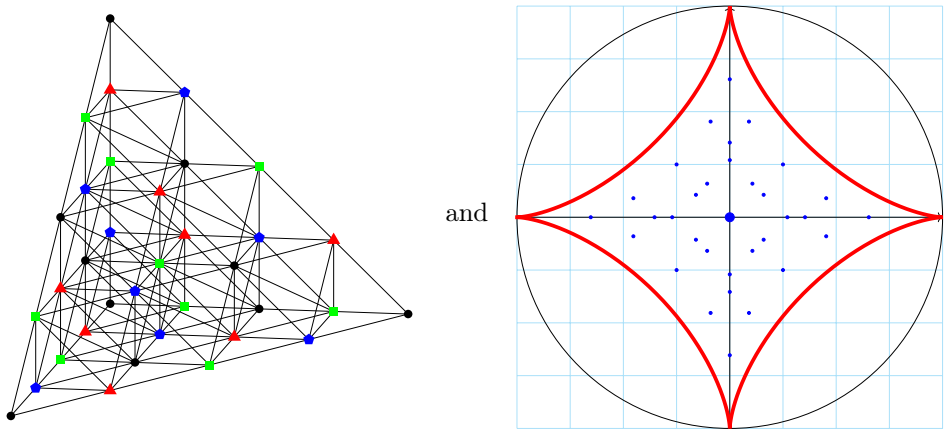
\begin{gather*}
\scalebox{.7}{\input{figures/tetrahedron-4}}
\quad\text{and}\quad
\scalebox{.7}{\input{figures/koorn-A-4-4}}.
\end{gather*}
The joint spectrum of $(A(\Gamma_{1}),A(\Gamma_{2}),A(\Gamma_{3}))$ is the Koornwinder variety $\vanset{\level}$.
\end{Example}

\begin{Example}[Type D]
For $\rank=2$ the type D graphs are type D Dynkin diagrams, for $\rank=3$ see \cite[Appendix 1]{MaMaMiTu-trihedral}, and for $\rank=4$ see \cite[Figures 3 and 4]{Oc-classification-sun}. Here are the graph of type $D$ for $\level=4$ and $\rank=4$ and the plot of the eigenvalues of $A(\Gamma_{1})$ \ochangedd{are}:
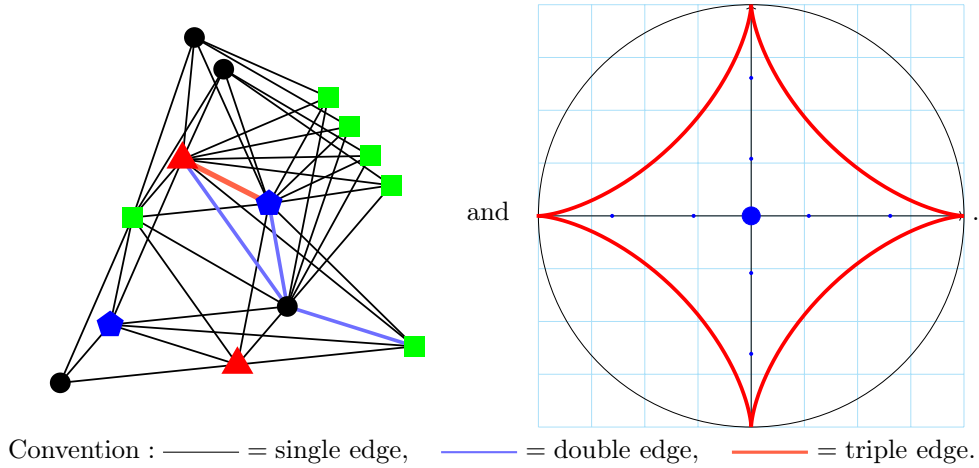
\begin{gather*}
\scalebox{1.7}{\input{figures/graph-D-4-4}}
\quad\text{and}\quad\scalebox{0.7}{\input{figures/koorn-D-4-4}}.
\\
\text{Convention}:
\begin{tikzpicture}[anchorbase]
\draw[black] (0,0) to (1,0);
\end{tikzpicture}
=\text{single edge}
,\quad
\begin{tikzpicture}[anchorbase]
\draw[lightneon, thick] (0,0) to (1,0);
\end{tikzpicture}
=\text{double edge}
,\quad
\begin{tikzpicture}[anchorbase]
\draw[tomato, very thick] (0,0) to (1,0);
\end{tikzpicture}
=\text{triple edge}.
\end{gather*}
Here, similarly as below, the blue slightly thicker edges are double edges and the red thicker edge is a triple edge. This graph is not included in \cite{Oc-classification-sun} and the joint spectrum of $(A(\Gamma_{1}),A(\Gamma_{2}),A(\Gamma_{3}))$ is a subset (with extra multiplicities) of the Koornwinder variety: there are $8$ points with multiplicity $1$ and $3$ points with multiplicity $2$.
\end{Example}

\begin{Example}[Type A conjugate]
For $\level=4$, the graph of type $2A_{\level}^{c}$ and the plot \changed{of the} eigenvalues of $A(\Gamma_{1})$ are:
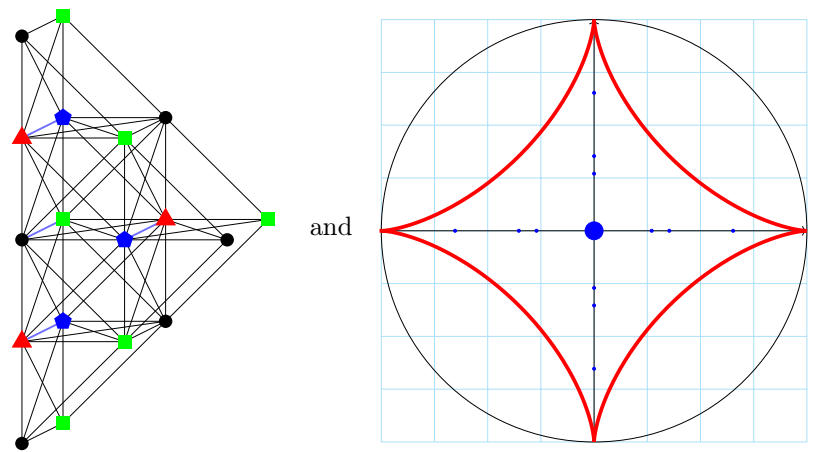
\begin{gather*}
\scalebox{.9}{\input{figures/graph-2Ac-4}}
\quad\text{and}\quad
\scalebox{0.7}{\input{figures/koorn-2Ac-4}}.
\end{gather*}
The joint spectrum of $(A(\Gamma_{1}),A(\Gamma_{2}),A(\Gamma_{3}))$ is a subset (with extra multiplicities) of the Koornwinder variety: there are $12$ points with multiplicity $1$ and $3$ points with multiplicity $2$.

For $\level=4$, the graph of type $2(A_{\level}^{c}/2)$ and the plot eigenvalues of $A(\Gamma_{1})$ are:
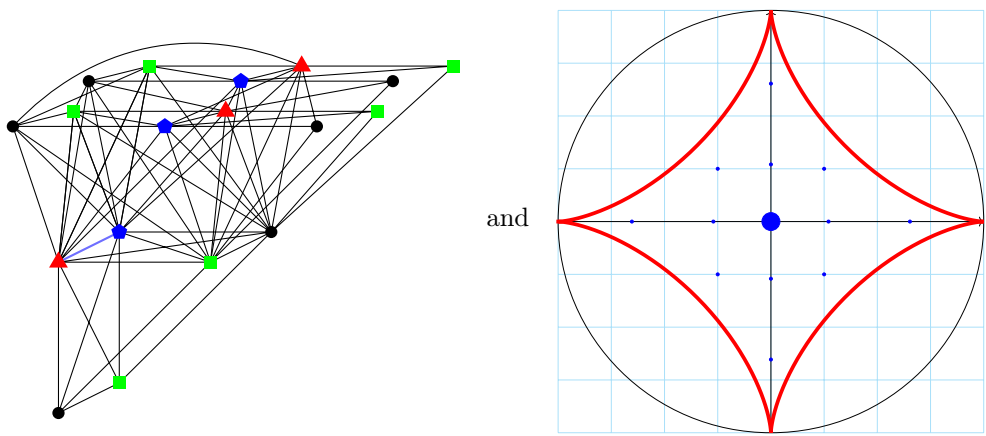
\begin{gather*}
\input{figures/graph-2Ac-2-4}
\quad\text{and}\quad
\scalebox{0.7}{\input{figures/koorn-2Ac-2-4}}.
\end{gather*}
The joint spectrum of $(A(\Gamma_{1}),A(\Gamma_{2}),A(\Gamma_{3}))$ is a subset (with extra multiplicities) of the Koornwinder variety: there are $12$ points with multiplicity $1$ and $3$ points with multiplicity $2$.
\end{Example}

\begin{Example}[Type E]
The exceptional graph $E_{4}$ and the plot \changed{of the} eigenvalues of $A(\Gamma_{1})$ are:
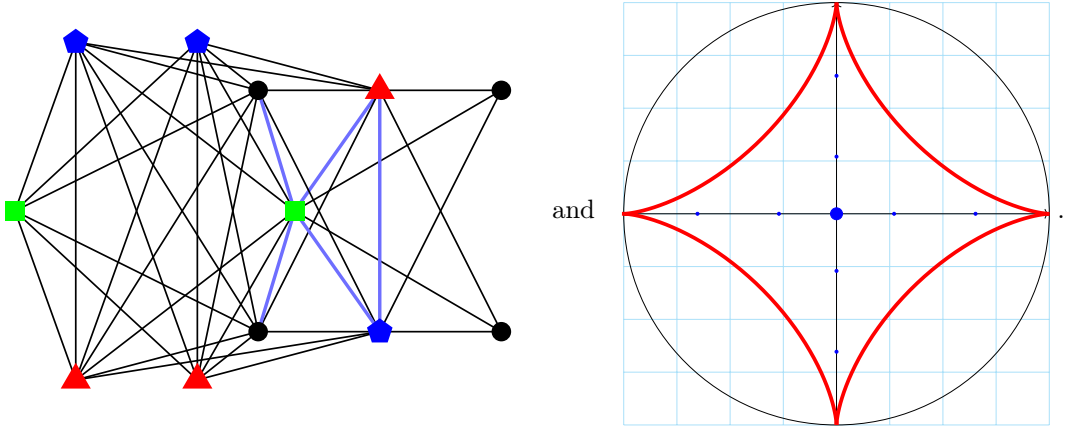
\begin{gather*}
\scalebox{1.6}{\input{figures/graph-E-4}}
\quad\text{and}\quad
\scalebox{0.7}{\input{figures/koorn-E-4}}.
\end{gather*}
The joint spectrum of $(A(\Gamma_{1}),A(\Gamma_{2}),A(\Gamma_{3}))$ is a subset (with extra multiplicities) of the Koornwinder variety: there are $8$ points with multiplicity $1$ and $1$ point with multiplicity $4$.
\end{Example}

Recall, from e.g. \cite[Definition 5.6]{MaMaMiTu-trihedral}, the notion of a transitive $\Zv$-representation.
These representations are the simple objects in the category of $\Zv$-representations, and we call them simple $\Zv$-representations.

\begin{Theorem}\label{T:integralTlevelA}
We have the following (potentially conjectural as indicated below). All of the following give well-defined simple $\Zv$-representations of $\nhecke[\level]$. 
\begin{enumerate}

\item $\M(\Gamma)$ for $\Gamma$ of type $A$ of rank $\rank$ and level $\level$. (Proven.)

\item $\M(\Gamma)$ for $\Gamma$ of type $D$ of rank $\rank$ and level $\level$. (Conjectural; verified in small cases.)

\item $\M(\Gamma)$ for $\Gamma$ of conjugate type $A$ of rank $\rank$ and level $\level$. (Conjectural; verified in small cases.)

\item $\M(\Gamma)$ for $\Gamma$ of type $E$. (Proven via computer, see \cite{LaTuVa-Nhedral-code}.)

\end{enumerate}
\end{Theorem}

\begin{proof}
\changedd{\textit{Part (a).} This is proven in \autoref{Pf:integralTlevelA-a}.}
\end{proof}

\begin{Remark}
One could ask the following problem: could one classify the graphs such that the condition in \autoref{L:integralTlevel} is satisfied? This question has been addressed in the literature for special cases, see e.g. \cite{Zu-gen-dynkin-diagrams}, \cite{Oc-classification-sun}, \cite{MaMaMiTu-trihedral}. We do not know the answer to this in general, not even for $\rank=3$.
In particular, there will be more simple $\Zv$-representations \changed{than} in \autoref{T:integralTlevelA} in general.
\end{Remark}

\subsection{Categorification}

\changed{There should be a categorification of the above story, but we decided not to include it here since that several notions we would need are not in the literature while writing this paper. Instead, we list here what one probably needs to change when compared to \cite{MaMaMiTu-trihedral}.}

\begin{enumerate}[label=(\roman*)]

\item The diagrammatic 2-category in \cite[Section 4.1]{MaMaMiTu-trihedral} should be replaced by its affine type $A_{\rank-1}$ 
analog. While writing this paper, at least to the best of our knowledge, there is no diagrammatic presentation of the 
relevant 2-category but one rather has to work with algebraic singular Soergel bimodules as in \cite{Wi-sing-soergel}. 
For the quantization of this category one needs to use the quantum Cartan matrix as in \cite{El-q-satake}.

\item \emph{Nhedral Soergel bimodules} of level $\infty$ \changed{probably can} then be defined similarly as in \cite[Section 4.2]{MaMaMiTu-trihedral}, replacing `secondary color' therein 
by a subset of $I$ of size $\rank-1$. The main statement \changed{should be} then the analog of \cite[Proposition 4.31]{MaMaMiTu-trihedral}, which \changed{should hold} verbatim and \changed{would provide} a 
categorification of $\nhecke$.

\item Under \autoref{QSH}, one can then \changed{likely} copy \cite[Section 4.3]{MaMaMiTu-trihedral} using $\rank$-colored webs (using e.g. the webs from \cite{CaKaMo-webs-skew-howe}) and \changed{would get} the 
analog of \cite[Proposition 4.48]{MaMaMiTu-trihedral}, categorifying $\nhecke[\level]$ for a fixed $\level$. Hereby, the Nhedral version of the clasps \changed{probably can} be defined 
using \cite{El-ladders-clasps} (which has been proven by now, see \cite{MaSp-type-A-webs}, and which might even give bases {\cf} \cite{AnStTu-cellular-tilting}).

\item Finally, \cite[Section 5.2]{MaMaMiTu-trihedral} \changed{likely} can be, mutatis mutandis, used to categorify \autoref{T:integralTlevelA} in the expected way. One \changed{would exploit} here the 
algebra objects that we give in \autoref{appNhedral} and the algebra object technology for Soergel bimodules developed in \cite{MaMaMiTu-algebra-objects}, \cite{MaMaMiTuZh-soergel-2reps}, \cite{MaMaMiTuZh-bireps}. Alternatively, this \changed{might} be done using a quiver similarly to \cite{AnTu-tilting}, but the quiver is already complicated for $\rank=3$, see \cite[Section 5.3]{MaMaMiTu-trihedral}.

\end{enumerate}

\section{Asymptotic \texorpdfstring{$\crg{\levelthing}{\rank}$}{G(M,M,N)}}\label{S:AsymptoticCategory}

Below we will use some standard terminology that can be found e.g. in 
\cite{EtGeNiOs-tensor-categories}. \changed{We encourage the reader to recall the conventions in the table of notations in \autoref{SS:notations}.} \changedd{In particular, we remind the reader that we have fixed an integer $\level$, and that $\levelthing=\level+\rank$, and that $\rooty$ is the primitive $2\levelthing$th root of unity $\exp(\iunit\pi/\levelthing)$.}

\subsection{The asymptotic category and its Drinfeld center}\label{SS:AsymCat}

We define the asymptotic category $\verg{\levelthing}{\rank}$ as a matrix category indexed by $\Z/\rank\Z$, with entries in suitable subcategories of $\slncatr$.

\begin{Definition}\label{D:AsymptoticCat}
The \emph{asymptotic category} $\verg{\levelthing}{\rank}$ is the monoidal additive $\C$-linear category whose:
\begin{enumerate}

\item \changedd{objects are matrices \changed{$(\obstuff{Y}_{ij})_{i,j\in\Z/\rank\Z}$, for \changed{$\obstuff{Y}_{ij}$} an object in $\slncatr_{i-j}$}},

\item morphisms between \changed{$(\obstuff{Y}_{ij})_{i,j\in\Z/\rank\Z}$} and \changed{$(\obstuff{Z}_{ij})_{i,j\in\Z/\rank\Z}$} are matrices $(f_{ij})_{i,j\in\Z/\rank\Z}$ with $f_{ij}$ being a morphism between \changed{$\obstuff{Y}_{ij}$} and \changed{$\obstuff{Z}_{ij}$} in $\slncatr$,

\end{enumerate}
and the tensor product is given by multiplication of matrices. We equip this category with the ribbon structure inherited from $\slncatr$.
\end{Definition}

\begin{Example}\label{E:EisZero}
When $\level=0$, then \changed{$\levelthing=\rank$ and} $\slncatr$ is equivalent to the category of $\C$-vector spaces, and $\verg{\levelthing}{\rank}$ is then equivalent to the direct sum of $\rank$ copies of the category of $\C$-vector spaces.
\end{Example}

\changedd{As a sneak preview, before defining anything, 
we will do the following with $\verg{\levelthing}{\rank}$:}
\begin{enumerate}[label=$\blacktriangleright$]

\item \changedd{In the final step \autoref{SS:BigCat} (following the introduction of some delicate combinatorics), we define a category that is Morita equivalent to $\verg{\levelthing}{\rank}$ via a straightforward matrix construction. We regard this new category as essentially identical to $\verg{\levelthing}{\rank}$. It has rank given by the ``subregular (=the first nontrivial) KL cell'' and corresponds to the asymptotic category associated with the middle dihedral cell, which is the subregular KL cell in this case. In this sense, we interpret $\verg{\levelthing}{\rank}$ as the ``subregular asymptotic category'' for $\crg{\levelthing}{\rank}$.}

\item \changedd{The Drinfeld center of $\verg{\levelthing}{\rank}$ is a modular category. The first question one would ask about a modular category is what are its $S$ and $T$-matrices. We study this Drinfeld center and then indeed compute parts of its $S$ and $T$-matrices in \autoref{T:DrinfeldCenter}. In contrast to the dihedral case, however, we were unable to compute the entire $S$-matrix.}

\item \changedd{We relate in \autoref{T:Comparison} these $S$ and $T$-matrices to the Fourier matrix and the eigenvalues of the Frobenius of a family of unipotent characters for
$\crg{\levelthing}{\rank}$. This is again in analogy with the asymptotic category for the subregular KL cell in the dihedral group.}

\end{enumerate}

\begin{Lemma}\label{L:AsymptoticIndecomposable}
The asymptotic category is a multifusion category of rank $\rank\tnumbers{\rank}{\level}$ and is indecomposable if and only if $\level\neq 0$.
\end{Lemma}

\begin{proof}
The first statement follows from \autoref{L:VerlindeAgain}, applied row-by-row. For $\level=0$ see \autoref{E:EisZero}. \changedd{Now, suppose $\level\neq 0$ and that the asymptotic category splits as $\catstuff{C}_{1}\oplus\catstuff{C}_2$ with $\catstuff{C}_1$ and $\catstuff{C}_2$ two nonzero multifusion categories. Then, there exist $i,j\in\Z/\rank\Z$ such that $\munit_i$ is a simple object of $\catstuff{C}_1$ and $\munit_{j}$ is a simple object of $\catstuff{C}_2$, where $\munit_k$ denotes the matrix with the unit object of $\slncatr$ at place $(k,k)$ and zero elsewhere. We consider any simple object $\obstuff{X}$ of $\verg{\levelthing}{\rank}$ with $(i,j)$ entry simple and other entries $0$. Such an object exists since $\level\neq 0$ so that the grading of $\slncatr$ by colors is faithful. The definition of the tensor product in the asymptotic category implies that $\munit_{i}\otimes X \otimes \munit_{j} = X$, but $\munit_{i}\otimes X = 0$ or $X\otimes \munit_{j}=0$ since $X$ lies in $\catstuff{C}_1$ or $\catstuff{C}_2$.}
\end{proof}

We now describe the Drinfeld center of the asymptotic category when $\level\geq 0$. The first case is easy: 

\begin{Example}
If $e=0$, then the Drinfeld center is simply $\rank$ copies of the category of vector spaces, with a trivial braiding.
\end{Example}

We now assume that $\level>0$.
Let us denote by $J_{3}(k)=k^{3}\prod_{p|k}(1-\frac{1}{p^{3}})$ Jordan's totient function. See \cite[A059376]{Oeis} for explicit values.
\changedd{We use the concept of \emph{modular closure} below, which is a special case of de-equivariantization from 
\cite[Section 8.23]{EtGeNiOs-tensor-categories}. This is the ``universal way'' of turning a non-modular fusion category into a modular one. Since the definition is quite involved, we will not recall it and we refer the reader to e.g. \cite{Br-cat-modulaires} or \cite{Mu-galois-theory} for more details.}

\begin{Theorem}\label{T:DrinfeldCenter}
We have the following.
\begin{enumerate}

\item $\dcenter{\verg{\levelthing}{\rank}}$ is equivalent as a ribbon category to the modular closure of $\dcenter{\slncatr}_{0}$. In particular, it is a modular category.

\item Write $\stabilizer{{\bk,\bm}}=\gcd(\stabilizer{\bk},\stabilizer{\bm})$. Simple objects of $\dcenter{\verg{\levelthing}{\rank}}$ are indexed by the set 
\begin{gather*}
\big\{(\bm,\bk,i)|
(\bm,\bk)\in X^{+}(\level)^{2}/(\Z/\rank\Z)\text{ with }\chi_{c}(\Ll_{\bm})=\chi_{c}(\Ll_{\bk})\text{, and } i\in\Z/\stabilizer{\bm,\bk}\Z\big\}. 
\end{gather*}
Here $\Z/\rank\Z$ acts diagonally on $X^{+}(\level)^{2}=X^{+}(\level)\times X^{+}(\level)$ \changed{via $(\bm,\bk)\mapsto (\bm^{\turnsymbol},\bk^{\turnsymbol})$}.

\item The rank of $\dcenter{\verg{\levelthing}{\rank}}$ is 
\begin{gather*}
\mathrm{rk}\,\dcenter{\verg{\levelthing}{\rank}}
=\frac{1}{\levelthing^{2}}\sum_{k\mid\gcd(\rank,\levelthing)}J_{3}(k)\binom{\levelthing/k}{\rank/k}^{2}.
\end{gather*}

\item For $\rank$ fixed and $\levelthing\to\infty$ we have
\begin{gather*}
\mathrm{rk}\,\dcenter{\verg{\levelthing}{\rank}}
\sim\tfrac{1}{(\rank!)^{2}}\cdot\levelthing^{2\rank-2}.
\end{gather*}

\item The $S$-matrix $S^{\ba}$ of $\dcenter{\verg{\levelthing}{\rank}}$ satisfies:
\begin{gather*}
\sum_{i=1}^{\stabilizer{\bk,\bm}}S^{\ba}_{(\bm,\bk,i),(\bm^{\prime},\bk^{\prime},i^{\prime})} = \tfrac{1}{\stabilizer{\bk,\bm}}S_{\bm,\bm^{\prime}}^{\sln[]}\overline{S}_{\bk,\bk^{\prime}}^{\sln[]}.
\end{gather*}
Moreover, if $\rank$ is prime, there exists a unique $(\bm,\bk)\in X^{+}(\level)^{2}/(\Z/\rank\Z)$ with $\stabilizer{\bk,\bm}\neq 1$ if and only if $\level\equiv 0\bmod\rank$ and we have
\begin{gather*}
S^{\ba}_{(\bm,\bk,i),(\bm,\bk,j)} 
=
\frac{1}{\rank^{2}}
\begin{cases}
S_{\bm,\bm}^{\sln[]}\overline{S}_{\bk,\bk}^{\sln[]} + (\rank{-}1)\changedd{\theta_{\bm}^{3}\overline{\theta_{\bk}}^{3}}\dim_{c}\slncatr & \text{if } i=j,
\\
S_{\bm,\bm}^{\sln[]}\overline{S}_{\bk,\bk}^{\sln[]} - \changedd{\theta_{\bm}^{3}\overline{\theta_{\bk}}^{3}}\dim_{c}\slncatr & \text{if } i\neq j,
\end{cases}
\end{gather*}
where \changed{$\theta$ is the ribbon element and} $\dim_{c}\slncatr$ is the categorical dimension. (See, for example, \cite[(3.3.9)]{BaKi-lecture-tensor} for the value of $\dim_{c}\slncatr$.) If $\rank$ is not prime, see \autoref{R:NotPrime}.

\item The $T$-matrix of $\dcenter{\verg{\levelthing}{\rank}}$ is the matrix with entries $T^{\ba}_{(\bm,\bk,i),(\bm^{\prime},\bk^{\prime},i^{\prime})} = \delta_{\bm,\bm^{\prime}}\delta_{\bk,\bk^{\prime}}\delta_{i,i^{\prime}}\theta_{\bm}^{-1}\theta_{\bk}$.

\end{enumerate}
\end{Theorem}

\begin{Example}\label{E:SMatrix}
Let $\rank=3$, $\level=3$ and $\levelthing=6$. The category $\slncatr$ has $10$ simple objects and $\dcenter{\verg{\levelthing}{\rank}}_{0}$ has then $34$ simple objects. Among these $34$ \changed{objects}, there are $11$ orbits of size $3$ and one of size $1$ under the action of $\Z/\rank\Z$. Therefore, $\dcenter{\verg{\levelthing}{\rank}}$ has $14$ simple objects, since each orbit of size $3$ will give one simple object, and the orbit of size $1$ will split into $3$ simple objects. 

The $S$-matrix of $\dcenter{\verg{\levelthing}{\rank}}$ is the $14$-by-$14$ matrix
\begin{gather*}
S^{\ba}=
\left(
\scalebox{0.6}{\begin{tikzpicture}[anchorbase]
\matrix[matrix of nodes,
nodes={text centered, minimum size=8mm},
column sep=-0.1cm,
row sep=-0.1cm,
cells={nodes={text centered, minimum size=8mm}},
ampersand replacement=\&] (m) {
\node[fill=tomato!50] {$1$}; \& \node[fill=tomato!50] {$1$}; \& \node[fill=tomato!50] {$3$}; \& \node[fill=tomato!50] {$1$}; \& $4$ \& $4$ \& $4$ \& $4$ \& $4$ \& $4$ \& $3$ \& $3$ \& $3$ \& $3$ \\
\node[fill=tomato!50] {$1$}; \& \node[fill=tomato!50] {$1$}; \& \node[fill=tomato!50] {$3$}; \& \node[fill=tomato!50] {$1$}; \& $4\junit$ \& $4\junit$ \& $4\junit$ \& $4\junit^{2}$ \& $4\junit^{2}$ \& $4\junit^{2}$ \& $3$ \& $3$ \& $3$ \& $3$ \\
\node[fill=tomato!50] {$3$}; \& \node[fill=tomato!50] {$3$}; \& \node[fill=tomato!50] {$-3$}; \& \node[fill=tomato!50] {$3$}; \& $0$ \& $0$ \& $0$ \& $0$ \& $0$ \& $0$ \& $9$ \& $-3$ \& $-3$ \& $-3$ \\
\node[fill=tomato!50] {$1$}; \& \node[fill=tomato!50] {$1$}; \& \node[fill=tomato!50] {$3$}; \& \node[fill=tomato!50] {$1$}; \& $4\junit^{2}$ \& $4\junit^{2}$ \& $4\junit^{2}$ \& $4\junit$ \& $4\junit$ \& $4\junit$ \& $3$ \& $3$ \& $3$ \& $3$ \\
$4$ \& $4\junit$ \& $0$ \& $4\junit^{2}$ \& $4$ \& $4\junit$ \& $4\junit^{2}$ \& $4$ \& $4\junit^{2}$ \& $4\junit$ \& $0$ \& $0$ \& $0$ \& $0$ \\
$4$ \& $4\junit$ \& $0$ \& $4\junit^{2}$ \& $4\junit$ \& $4\junit^{2}$ \& $4$ \& $4\junit^{2}$ \& $4\junit$ \& $4$ \& $0$ \& $0$ \& $0$ \& $0$ \\
$4$ \& $4\junit$ \& $0$ \& $4\junit^{2}$ \& $4\junit^{2}$ \& $4$ \& $4\junit$ \& $4\junit$ \& $4$ \& $4\junit^{2}$ \& $0$ \& $0$ \& $0$ \& $0$ \\
$4$ \& $4\junit^{2}$ \& $0$ \& $4\junit$ \& $4$ \& $4\junit^{2}$ \& $4\junit$ \& $4$ \& $4\junit$ \& $4\junit^{2}$ \& $0$ \& $0$ \& $0$ \& $0$ \\
$4$ \& $4\junit^{2}$ \& $0$ \& $4\junit$ \& $4\junit^{2}$ \& $4\junit$ \& $4$ \& $4\junit$ \& $4\junit^{2}$ \& $4$ \& $0$ \& $0$ \& $0$ \& $0$ \\
$4$ \& $4\junit^{2}$ \& $0$ \& $4\junit$ \& $4\junit$ \& $4$ \& $4\junit^{2}$ \& $4\junit^{2}$ \& $4$ \& $4\junit$ \& $0$ \& $0$ \& $0$ \& $0$; \\
$3$ \& $3$ \& $9$ \& $3$ \& $0$ \& $0$ \& $0$ \& $0$ \& $0$ \& $0$ \& $-3$ \& $-3$ \& $-3$ \& $-3$ \\
$3$ \& $3$ \& $-3$ \& $3$ \& $0$ \& $0$ \& $0$ \& $0$ \& $0$ \& $0$ \& $-3$ \& \node[fill=blue!50] {$9$}; \& \node[fill=blue!50] {$-3$}; \& \node[fill=blue!50] {$-3$}; \\
$3$ \& $3$ \& $-3$ \& $3$ \& $0$ \& $0$ \& $0$ \& $0$ \& $0$ \& $0$ \& $-3$ \& \node[fill=blue!50] {$-3$}; \& \node[fill=blue!50] {$9$}; \& \node[fill=blue!50] {$-3$}; \\
$3$ \& $3$ \& $-3$ \& $3$ \& $0$ \& $0$ \& $0$ \& $0$ \& $0$ \& $0$ \& $-3$ \& \node[fill=blue!50] {$-3$}; \& \node[fill=blue!50] {$-3$}; \& \node[fill=blue!50] {$9$}; \\
};
\end{tikzpicture}}
\right)
,
\end{gather*}
where $\junit=\rooty^{4}$ \changed{to ease notation}. The top colored block corresponds to the $S$-matrix of $\slncatr_{0}$ and the bottom colored block corresponds to the object of $\dcenter{\verg{\levelthing}{\rank}}_{0}$ that splits in three.
\end{Example}

\begin{Remark}
Let $\mathfrak{A}_{4}$ denote the alternating group of order twelve.
The \changed{$14$-by-$14$} matrix in \autoref{E:SMatrix} is the same as the $S$-matrix of \ochanged{the} Drinfeld center of $\mathfrak{A}_{4}$-graded vector spaces (the latter, by duality, is the same as the Drinfeld center of the category of $\mathfrak{A}_{4}$-representations).
However, the appearance of nonnegative integers in the first row and column of the matrix is a small number coincidence and we do not expect any nice description of the $S$-matrix in terms of a group in general\changed{, since the first column corresponds to the quantum dimension of simple objects}.
\end{Remark}

\begin{proof}[Proof of \autoref{T:DrinfeldCenter}]
\textit{Part (a).}
Since the asymptotic category $\verg{\levelthing}{\rank}$ is indecomposable by \autoref{L:AsymptoticIndecomposable}, \cite[Theorem 2.5.1]{KoZh-center-functor} shows that its Drinfeld center is equivalent to the Drinfeld center of any of its diagonal components $\slncatr_{0}$. We then observe that $\dcenter{\slncatr_{0}}$ is monoidally equivalent to the modular closure of $\dcenter{\slncatr}_{0}$, which follows immediately from \cite[Proposition 8.23.11]{EtGeNiOs-tensor-categories}.

\textit{Part (b).}
Let $(\placeholder)^{\rev}$ denote the reverse category as in \cite[Definition 8.1.4]{EtGeNiOs-tensor-categories}. Since $\slncatr$ is modular, its Drinfeld center is equivalent to $\slncatr\boxtimes\slncatr^{\rev}$, see, for example, 
\cite[Proposition 8.6.3]{EtGeNiOs-tensor-categories}. The degree zero part of $\slncatr\boxtimes\slncatr^{\rev}$ is $\bigoplus_{i\in\Z/\rank\Z}\slncatr_{i}\boxtimes(\slncatr^{\rev})_{i}$, and its symmetric center is generated by the object $\Ll_{\level\omega_1}\boxtimes\Ll_{\level\omega_1}$\changed{, by \cite[Lemma 2.2]{La-crossed-S-matrix} and \autoref{L:GradingBraiding}}. Since $\Ll_{\level\omega_1}\otimes \Ll_{\bm} = \Ll_{\bm^{\turnsymbol}}$, the result follows from \cite[Remarques 4.5, 1)]{Br-cat-modulaires} or \cite[Corollary 5.3]{Mu-galois-theory}.

\textit{Part (c).} This is done in \autoref{Pf:DrinfeldCenter-c}.

\textit{Part (d).} For fixed $\rank$, $\mathrm{rk}\,\dcenter{\verg{\levelthing}{\rank}}$ is polynomial in $\levelthing$ and the leading term is obtained for $k=1$. The result follows from $J_{3}(1)=1$ and $\binom{\levelthing}{\rank}\sim\tfrac{\levelthing^{\rank}}{\rank!}$ when $\levelthing\to\infty$.

\textit{Part (e).} This is done in \autoref{Pf:DrinfeldCenter-e}.

\textit{Part (f).}
This follows immediately from \cite[Proposition 4.2]{Mu-galois-theory}.
\end{proof}

\begin{Remark}\label{R:NotPrime}
When $\rank$ is not prime, the calculation of the $S$-matrix of $\dcenter{\verg{\levelthing}{\rank}}$ is more intricate. For $\gcd(\level,\rank) \neq 1$, \autoref{L:ObjectsStabilizer} shows that several objects of $\dcenter{\slncatr}_{0}$ can split in the modular closure. 
\end{Remark}

\subsection{A family of unipotent characters}\label{SS:FamilyUnipotent}

Malle \cite{Ma-unipotente-grade} has introduced the notion of unipotent characters for the complex reflection group $\crg{\levelthing}{\rank}$. These characters are partitioned into \emph{families}, similarly to the unipotent characters of a finite group of Lie type as e.g. in \cite[Chapter 4]{Lu-characters-reductive-groups}. We consider a certain family $\family$ of unipotent characters which, in terms of the combinatorics developed by Malle, are indexed by so-called \emph{$\levelthing$-symbols} with entries in the multiset $\{0^{\level},1^{\rank}\}$, but we do not use this explicitly.

\begin{Remark}\label{R:AValue}
For $\rank=2$ the family $\family$ corresponds to the subregular KL cell of $\mathbf{a}$-value $1$\changed{, which is the only non trivial cell}. In general, the families we consider have $\mathbf{a}$-value $\rank(\rank-1)/2$\changed{, and they play a role \changedd{analogous} to the subregular cell}. \ochangedd{This our reason to use the $\mathbf{a}$-value as the shift in \autoref{R:shift}.} \changed{There are many other non trivial cells that we do not consider here.}
\end{Remark}

\changed{The combinatorics of the family $\family$ is a bit demanding, and some details are postponed to \autoref{Pf:M-symbols}. For now, we just define what is necessary to state our main theorem in this section.}

\changed{We now follow and borrow the notations of \cite[Section 6C]{Ma-unipotente-grade}.}
Let $Y=\{1,\dots,\levelthing\}$ and $\pi\colon Y \to\N$ defined by $\pi(y)=1$ if $1\leq y\leq\rank$ and $\pi(y)=0$ otherwise. \changedd{(Our $\pi$ only takes values in $\{0,1\}$, the reason why we choose $\N$ in its definition is explained in \autoref{R:Symbols}.)} 
We consider the set
\begin{gather*}
\Psi(Y,\pi)=\big\{f\colon Y\to\{0,\ldots,\levelthing-1\}| 
\sum_{y\in Y}f(y)\equiv{\textstyle\binom{\levelthing}{2}}\bmod\levelthing, 
f_{|\pi^{-1}(i)}\text{ is strictly increasing for all }i\in\N\big\}.
\end{gather*}
\changed{There exists an action of $\Z/\levelthing\Z$ on $\Psi(Y,\pi)$ \changedd{defined in \autoref{Pf:M-symbols}}, and the elements of $\family$ are indexed by orbits under this action. We denote equivalence classes by $[f]$. An equivalence class $[f]$ will index several unipotent characters of $\family$, as many as the size of the stabilizer of $f$ under the action of $\Z/\levelthing\Z$, which we will denote by $([f],\changedd{s})$ for $\changedd{s}$ in a finite indexing set.}

\begin{Remark}\label{R:Symbols}
\changedd{The set $\Psi(Y,\pi)$ can be described with Malle's $\levelthing$-symbols \cite[Section 6A]{Ma-unipotente-grade}, which are a generalization of Lusztig's symbols used in the classification of unipotent characters of type D \cite[Section 4.6]{Lu-characters-reductive-groups}. The $\levelthing$-symbol associated to a function $f\in\Psi(Y,\pi)$ has its $j$th row filled with the elements of $\pi(f^{-1}(j))$. The symbols have their entries in $\pi(Y)$, which, for other families, could be \changeddd{some subset of $\N$ other} than $\{0,1\}$.}

\changedd{Explicitly, with $\level=2$ and $\rank=2$, so that $\levelthing=4$, the function $f\in\Psi(Y,\pi)$ given by $f(1)=0,f(2)=2,f(3)=1$ and $f(4)=3$ is the $4$-symbol
\[
\begin{pmatrix}
1 \\
0 \\
1 \\
0
\end{pmatrix}.
\]
The action of $\Z/\levelthing\Z$ on $\Psi(Y,\pi)$ is the cyclic shift on rows of $\levelthing$-symbols. The above $4$-symbol has then a stabilizer of order $2$.}
\end{Remark}

To each family of unipotent characters, Malle \cite{Ma-unipotente-grade} has attached a \emph{Fourier matrix} and eigenvalues of the \emph{Frobenius}, generalizing the constructions from the groups of Lie type as, for example, in 
\cite{Lu-characters-reductive-groups}. The eigenvalue of the Frobenius of the unipotent character of $\family$ indexed by $([f],\changedd{s})$ is
\begin{gather}\label{Eq:Frobenius}
\frob([f],\changedd{s})=\rooty^{-\alpha(f)}
\quad\text{with}\quad
\alpha(f)=\frac{\levelthing(1-\levelthing^{2})}{6}-\sum_{y\in Y}(f(y)^{2}+\levelthing f(y)).
\end{gather}
For the Fourier matrix, we will adopt the conventions of \cite[Section 5.4]{La-thesis} which differ from \cite[Section 6C]{Ma-unipotente-grade} by some signs. For $f\in\Psi(Y,\pi)$, we set $\varepsilon(f)=(-1)^{c(f)+\gamma(f)}$ where
\begin{gather*}
c(f) = \left\lvert\{(y,y^{\prime})\in Y^{2} \mid y<y^{\prime}, f(y)<f(y^{\prime})\}\right\rvert
\quad\text{and}\quad
\gamma(f) = \frac{\levelthing-1}{\levelthing}\Big(\binom{\levelthing}{2} - \sum_{y\in Y}f(y)\Big).
\end{gather*}
Let $\mathcal{S}=(\rooty^{-2ij})_{0\leq i,j < \levelthing}$ and define $\tau = (-1)^{\binom{\levelthing}{2}}\det(\mathcal{S})$. We view $\mathcal{S}$ as the matrix of an endomorphism of $V=\bigoplus_{i=0}^{\levelthing-1}\mathbb{C}v_{i}$. Given $f\colon Y \rightarrow \{0,\ldots,\levelthing-1\}$ strictly increasing on $\pi^{-1}(1)=\{1,\ldots,\rank\}$ and on $\pi^{-1}(0) = \{\rank+1,\ldots,\levelthing\}$, we set
\begin{gather*}
v_{f} = (v_{f(1)}\wedge\cdots\wedge v_{f(\rank)})\otimes (v_{f(\rank+1)}\wedge\cdots\wedge v_{f(\levelthing)}) \in \Lambda^{\rank}V\otimes \Lambda^{\level}V.
\end{gather*}
Such vectors form a basis of $\Lambda^{\rank}V\otimes \Lambda^{\level}V$ and we denote by $(\Lambda^{\rank}\mathcal{S}\otimes \Lambda^{\level}\mathcal{S})_{f,g}$ the entries of the matrix of the endomorphism $\Lambda^{\rank}\mathcal{S}\otimes \Lambda^{\level}\mathcal{S}$ in this basis.

We define the \emph{pre-Fourier matrix} $\prefourier$ of the family $\family$, which \ochangeddd{is} indexed by orbits of $\Psi(Y,\pi)$ under the action of $\Z/\levelthing\Z$, by 
\begin{gather*}
\prefourier_{[f],[g]} = \frac{(-1)^{\levelthing-1}\levelthing}{\tau}\varepsilon(f)\varepsilon(g)\overline{(\Lambda^{\rank}\mathcal{S}\otimes \Lambda^{\level}\mathcal{S})_{f,g}},
\end{gather*}
for all $f,g\in\Psi(Y,\pi)$\changed{, where the bar is the complex conjugate}. 

\begin{Remark}
\changedd{Neither the eigenvalues of the Frobenius nor the entries of pre-Fourier matrix depend on $\changedd{s}$ in the pair $([f],\changedd{s})$, but depend only on $[f]$.} 

\changedd{In contrast to the pre-Fourier matrix, the entries of the Fourier matrix itself of the family $\family$ depend on \changeddd{$s$} in the pair $([f],\changedd{s})$. The Fourier matrix is not entirely determined, see \cite[Remark 5.4.35]{La-thesis}.}
\end{Remark}

Let $\bar{\iota}$ be the bijection, \ochanged{as} in \autoref{Pf:iota}, \changedd{between the $\Z/\levelthing\Z$-orbits of $\Psi(Y,\pi)$ and the $\Z/\rank\Z$-orbits of $\{(\bm,\bk)\in (X^{+}(e))^{2}\mid \chi_{c}(\Ll_{\bm})=\chi_{c}(\Ll_{\bk})\}$}. Recall the $T$-matrix $T^{\ba}$ of the modular category $\dcenter{\verg{\levelthing}{\rank}}$ \changed{in \autoref{T:DrinfeldCenter}.(f)}. We denote the $S$-matrix of the degree zero part of the Drinfeld center of 
$\slncatr$ by $S^{0}$.

\changed{The next theorem compares the Fourier matrix of $\family$ and the $S$-matrix of the modular category $\dcenter{\verg{\levelthing}{\rank}}$. The proof relies on technical calculations, since these two matrices are quite different in nature.}

\begin{Theorem}\label{T:Comparison}
We have the following.
\begin{enumerate}
\item The cardinal of $\family$ is equal to the rank of $\dcenter{\verg{\levelthing}{\rank}}$.
\item For all $f,g\in \Psi(Y,\pi)$, 
\begin{gather*}
\prefourier_{[f],[g]} = \text{factor}(f,g) \cdot
S_{\bar{\iota}(f),\bar{\iota}(g)}^{0}
,\quad\quad
\text{factor}(f,g)=
\tfrac{(-1)^{\sum_{i=1}^{\rank}(f(i)+g(i))}\varepsilon(f)\varepsilon(g)}{\sqrt{\dim_{\changed{c}}(\dcenter{\verg{\levelthing}{\rank}})}}.
\end{gather*}
Thus, up to a diagonal change of basis, the pre-Fourier matrix and the $S$-matrix of $\dcenter{\slncatr}_{0}$ coincide.
\item For all $f\in \Psi(Y,\pi)$, 
\begin{gather*}
\frob([f],\changedd{s}) = T_{\bar{\iota}(f),\bar{\iota}(f)}^{\ba}.
\end{gather*}
\end{enumerate}
\end{Theorem}

\begin{proof}
\changedd{\textit{Part (a).} This is proven in \autoref{Pf:Comparison-a}.}

\changedd{\textit{Part (b).} This is proven in \autoref{Pf:Comparison-b}.}

\changedd{\textit{Part (c).} This is proven in \autoref{Pf:Comparison-c}.}
\end{proof}

\begin{Remark}
By \autoref{T:DrinfeldCenter}.(e), the $S$-matrix of $\dcenter{\verg{\levelthing}{\rank}}$ satisfies \cite[Conjecture 5.4.34]{La-thesis}, and we think it hence deserves to be called the Fourier matrix of the family $\family$.
\end{Remark}

\subsection{The big asymptotic category related to a Calogero--Moser cell}\label{SS:BigCat}

We say that a unipotent character represented by $([f],\changedd{s})$ is in the principal series if $\lvert f^{-1}(\ochanged{j}) \rvert = 1$ for every $0\leq \ochanged{j} \leq \ochanged{\levelthing}-1$. We denote by $\family_{0}$ the unipotent characters of $\family$ lying in the principal series. They are in bijection with a subset of irreducible complex representations of $\crg{\levelthing}{\rank}$, which are represented by $\levelthing$-partitions of $\rank$ with each entry being $(1)$ or empty, see \cite[Section 6A]{Ma-unipotente-grade} for more details. By a (slight) abuse of notations, we will also denote by $\family_{0}$ the set of irreducible complex representations of $\crg{\levelthing}{\rank}$ corresponding to the unipotent characters in the principal series. If $\levelthing > \rank$, a result of Bellamy \cite{Be-calogero-moser} states that $\family_{0}$ is a CM family of irreducible complex representations of $\crg{\levelthing}{\rank}$ for a specific choice of parameters, the so-called \emph{spetsial} ones.

\begin{Remark}
The spetsial parameters correspond to equal parameters for Iwahori--Hecke algebras.
\end{Remark}

In \cite{BoRo-calogero-moser}, Bonnaf\'e--Rouquier construct a partition of a finite complex reflection group into left, right and two-sided CM cells, and conjecture that these cells coincide with the KL cells provided the complex reflection group is a Coxeter group. They also construct a bijection between the CM families and the two-sided CM cells. Let us denote by $\cmcell$ the two-sided CM cell associated with the family $\family_{0}$.

We now define what we call the \emph{big asymptotic category}. 

\begin{Definition}\label{D:BigCat}
Let $\verG{\levelthing}{\rank}$ be the matrix category over $\verg{\levelthing}{\rank}$ of size $(\rank{-}1)!$.
\end{Definition}

Since $\verg{\levelthing}{\rank}$ is a ribbon category, so it the big asymptotic category $\verG{\levelthing}{\rank}$.

\begin{Theorem}\label{T:CMCells}
We have the following.
\begin{enumerate}
\item $\verg{\levelthing}{\rank}$ and $\verG{\levelthing}{\rank}$ are Morita equivalent (they have the same $2$-representation theory).
\item The Drinfeld centers of $\verg{\levelthing}{\rank}$ and $\verG{\levelthing}{\rank}$ are equivalent as modular categories.
\item $\mathrm{rk}\ \verG{\levelthing}{\rank} = \tfrac{(\rank!)^{2}}{\levelthing}\binom{\levelthing}{\rank} = \lvert \cmcell \rvert$.
\item $\ggroupc{\verG{\levelthing}{\rank}} \simeq \bigoplus_{V\in\family_{0}} \mathrm{Mat}_{\dim_{\C}V}(\C)$.
\end{enumerate}
\end{Theorem}

\begin{proof}
\textit{Part (a).} 
This is \cite[Proposition 2.27]{MaMaMiTuZh-bireps}.

\textit{Part (b).}
This follows from \cite[Theorem 2.5.1, (7)]{KoZh-center-functor}.

\textit{Part (c).} The rank of $\verG{\levelthing}{\rank}$ is $((\rank{-}1)!)^{2}\mathrm{rk}\ \verg{\levelthing}{\rank} = ((\rank{-}1)!)^{2}\rank\tnumbers{\rank}{\level}$. The second equality follows from \cite[Theorem 12.2.7, (c)]{BoRo-calogero-moser} and taking dimensions in \textit{(d)}.

\textit{Part (d).} \ochangedd{This is proven in} \autoref{Pf:CMCells-d}.
\end{proof}

\section{Appendix: Some more technical \changed{proofs}}\label{app}

\subsection{Missing proofs for \autoref{S:slnthings}}\label{appslnthings}

\subsubsection{\changedd{Proof of \autoref{L:KoornwinderVarietyZ}}}\label{Pf:KoornwinderVarietyZ}
Given $\bm\in X^{+}$, we define $\Efun_{\bm}\colon E \to \C$ by 
\begin{gather*}
\Efun_{\bm}(\bs) = \sum_{w\in\sym_{\rank}}(-1)^{l(w)}\exp\left(\iunit (w(\bm),\bs)\right).
\end{gather*}

The functions $\Zfun_{i}$ and $\Efun_{\bm}$ are invariant by $\bs\mapsto \bs+2\pi\alpha_{i}$ and therefore define functions on the torus $T^{\rank-1}=E/2\pi Y$, where $Y$ is the root lattice. Moreover, the functions $\Zfun_{i}$ are invariant and the functions $\Efun_{\bm}$ are antiinvariant under the action of the symmetric group. The fundamental domain of the quotient $T^{\rank-1}/\sym_{\rank}$ is equal to
\begin{gather*}
D=\left\{\sum_{i=1}^{\rank-1}\lambda_{i}\alpha_{i} \Big| 2\lambda_{i}\geq \lambda_{i-1}+\lambda_{i+1} \text{ for all } 1\leq i <N, \lambda_{1}+\lambda_{\rank-1}\leq 2\pi\right\}.
\end{gather*}
Weyl's denominator formula implies that the zeroes of $\Efun_{\rho}$ are on the boundary of $D$. Therefore the functions $\Zfun_{i}$ and $\Efun_{\bm+\rho}/\Efun_{\rho}$ are defined on the interior of $D$.

Using Weyl's character formula we have, for $\bs$ in the interior of $D$, we have
\begin{gather*}
\pxy{\bm}\big(\Zfun_{1}(\bs),\dots,\Zfun_{\rank-1}(\bs)\big) = \Efun_{\bm+\rho}(\bs)/\Efun_{\rho}(\bs).
\end{gather*}
If moreover $\bm\in X^{+}$ is such that $\sum\bm=e+1$, then \cite[Equation (3.3.11)]{BaKi-lecture-tensor} implies that, for $\bs\in\vanpar{\level}$, we have $\pxy{\bm}(\Zfun_{1}(\bs),\dots,\Zfun_{\rank-1}(\bs))=0$.\null\hfill$\square$

\subsubsection{\changedd{Proof of \autoref{T:KoornwinderVariety} \textit{Part (b).}}}\label{Pf:KoornwinderVariety-b}
\changed{The dimension of the quotient by $\vanideal{\level}$ being of dimension $\tnumbers{\rank}{\level}$, Hilbert's Nullstellensatz (via e.g. \cite[Corollary I.7.4]{Fu-algebraic-curves}) implies that $\vanset{\level}$ is discrete with at most $\tnumbers{\rank}{\level}$ points.} For the converse, we use the following important lemma:

\begin{Lemma}\label{L:ZInjective}
The map $\bs\mapsto(\Zfun_{1}(\bs),\dots,\Zfun_{\rank-1}(\bs))$ is an injection on the interior of the fundamental domain $D$.
\end{Lemma}

\begin{proof}
In \cite[Equation 5.9]{Be-Chebyshev-Laplace-Beltrami}, Beerends calculated the value of the Jacobian $\left\lvert \frac{\partial(\Zfun_{1},\dots,\Zfun_{\rank-1})}{\partial(\sigma_{1},\ldots,\sigma_{\rank-1})}\right\rvert$:
\begin{gather*}
\left\lvert \frac{\partial(\Zfun_{1},\dots,\Zfun_{\rank-1})}{\partial(\sigma_{1},\ldots,\sigma_{\rank-1})}\right\rvert = \prod_{1\leq i < j \leq \rank} \lvert e^{\iunit(\sigma_{i}-\sigma_{i-1})} - e^{\iunit(\sigma_{j}-\sigma_{j-1})}\rvert = 2^{\binom{\rank}{2}}\prod_{1\leq i < j \leq \rank} \left\lvert\sin\left(\frac{\sigma_{i}-\sigma_{i+1}-\sigma_{j}+\sigma_{j+1}}{2}\right)\right\rvert.
\end{gather*}
Therefore, the Jacobian vanishes only on the reflecting hyperplanes of the representation $E$ of the symmetric group, and we are done.
\end{proof}

Since $\vanpar{\level}$ lies in the interior of $D$, \autoref{L:KoornwinderVarietyZ} implies that we have found $\tnumbers{\rank}{\level}$ points in $\vanset{\level}$.\null\hfill$\square$

\subsubsection{\changedd{Proof of \autoref{T:KoornwinderVariety} \textit{Part (c).}}}\label{Pf:KoornwinderVariety-c} We need the following two lemmas.

\begin{Lemma}\label{L:Newton}
Given $\bs=\sum_{j=1}^{\rank-1}\sigma_{j}\alpha_{j}\in E$, the function $\Zfun_{i}(\bs)$ is the evaluation of the $i$th elementary symmetric function at $(e^{\iunit\sigma_{1}},e^{\iunit(\sigma_{2}-\sigma_{1})},\dots,e^{\iunit(\sigma_{\rank-1}-\sigma_{\rank-2})},e^{-\iunit\sigma_{\rank-1}})$.
\end{Lemma}

\begin{proof}
The function $\Zfun_{i}$ is the evaluation of the character of the fundamental $\sln$-representation $\Ll_{\omega_{i}}$ at the diagonal matrix $\mathrm{diag}(\iunit\sigma_{1},\iunit(\sigma_{2}-\sigma_{1}),\dots,\iunit(\sigma_{\rank-1}-\sigma_{\rank-2}),-\iunit\sigma_{\rank-1})$. The claim follows from the usual correspondence between the character of $\Ll_{\omega_{i}}$ and the $i$th elementary symmetric function.
\end{proof}

Using Newton's identities, we can hence express $\Zfun_{i}(\bs)$ in terms of the power sums $\Zfun_{1}(k\bs)$.

\begin{Lemma}\label{L:OrbitsKoornwinder}
For $\bs\in\vanpar{\level}$, we have $\Zfun_{i}(\bs^{\turnsymbol}) = \rankroot^{-i}\Zfun_{i}(\bs)$.
\end{Lemma}

\begin{proof}
For $\bk=\sum_{i=1}^{\rank-1}k_{i}\omega_{i}$, we write $\lambda_{i}(\bk)$ its coordinates in the basis $(\alpha_{1},\ldots,\alpha_{\rank-1})$. Using \cite[Plate I]{Bo-chapters-4-6}, we explicitly have
\begin{gather*}
\lambda_{i}(\bk) = \frac{1}{\rank}\left((\rank-i)\sum_{j=1}^{i-1}jk_{j} + (\rank-i)ik_{i} + i\sum_{j=i+1}^{\rank-1}(\rank-j)k_{j}\right).
\end{gather*}
Then, one easily show that $\rank(\lambda_{i}(\bk)-\lambda_{i-1}(\bk))=-\sum_{j=1}^{\rank-1}jk_{j} + \rank\sum_{j=i}^{\rank-1}k_{j}$. A straightforward calculation also shows that, for $\bk\in X^{+}(\level)$ and $1\leq i \leq \rank-1$, we have
\begin{gather*}
\lambda_{i+1}(\bk^{\turnsymbol}+\rho)-\lambda_{i}(\bk^{\turnsymbol}+\rho) = \lambda_{i}(\bk+\rho)-\lambda_{i-1}(\bk+\rho) - \frac{\rank+\level}{\rank} \quad\text{and}\\
\lambda_{1}(\bk^{\turnsymbol}+\rho) = -\lambda_{\rank-1}(\bk+\rho) - \frac{\rank+\level}{\rank} + \level+\rank.
\end{gather*}
The result then follows from \autoref{L:Newton}.
\end{proof}

\autoref{L:OrbitsKoornwinder} and \autoref{T:KoornwinderVariety}.(a) \changed{imply} that the map 
\begin{equation*}
\ochangedd{X^{+}(\level)\to\vanset{\level},\ \bk\mapsto (\Zfun_{i}(2\iunit\pi(\bk+\rho)/(\level+\rank)))_{1\leq i \leq \rank-1}}
\end{equation*}
is a bijection which is $\Z/\rank\Z$-equivariant. We then do the counting in $X^{+}(\level)$ instead of in the Koornwinder variety. 

\begin{Lemma}\label{L:periodicity}
\changedd{Let us temporarily write $X^{+}_{\rank}(\level)$ instead of $X^{+}(\level)$. Let $\lambda=(\lambda_1,\dots,\lambda_{\rank-1})\in X^{+}_{\rank}(\level)$ with stabilizer of order $m$ under the action $\placeholder{}^{\turnsymbol}$. Then $m\mid \level$ and $\lambda'=(\lambda_{1},\dots,\lambda_{\rank/m-1})\in X^{+}_{\rank/m}(\level/m)$ and its stabilizer under the action $\placeholder{}^{\turnsymbol}$ is of order $1$.}

\changedd{This defines a bijection between the elements of $X^{+}_{\rank}(\level)$ with stabilizer of order $m$ and the elements of $X^{+}_{\rank/m}(\level/m)$ with stabilizer of order $1$.}
\end{Lemma}

\begin{proof}
\changedd{Since the stabilizer of $\lambda$ is of order $m$, the extended tuple $(\lambda_{1},\dots,\lambda_{\rank-1},\level-\bsum\lambda)$ is $\rank/m$-periodic, and $\rank/m$ is the smallest period. Therefore $\lambda_{1},\dots,\lambda_{\changeddd{N/m}}$ appear each $m$ times in the extended tuple so that $\lambda_{1}+\dots+\lambda_{\changeddd{N/m}}=\level/m$. This implies that $\lambda'=(\lambda_{1},\dots,\lambda_{\rank/m-1})\in X^{+}_{\rank/m}(\level/m)$. This element has stabilizer of order $1$ since $\rank/m$ is the smallest period of the extended tuple of $\lambda$.}

\changedd{Conversely, if $\lambda'=(\lambda_{1},\dots,\lambda_{\rank/m-1})\in X^{+}_{\rank/m}(\level/m)$ has a stabilizer of order $1$, then $\sum_{j=0}^{m-1}\sum_{i=1}^{\rank/m}\lambda_{i}\omega_{i+j\rank/m}$ (where $\omega_{N}=0$) is in $X^{+}_{\rank}(\level)$ and is of stabilizer of order $m$.}
\end{proof}

Denote by $c(\rank,\level,m)$ the number of weights in $X^{+}(\level)$ with stabilizer of order $m$. \changedd{\autoref{L:periodicity} implies that $m\mid \level$ and that $c(\rank,\level,m)=c(\rank/m,\level/m,1)$.} It then suffices to prove the formula for $m=1$.

We proceed by induction on $\gcd(\rank,\level)$. If $\gcd(\rank,\level)=1$ then all elements of $X^{+}(\level)$ have a stabilizer of order $1$. Therefore, $c(\rank,\level,1)=\tnumbers{\rank}{\level} = \binom{\rank-1}{\levelthing-1} = \tfrac{\levelthing}{\rank}\binom{\rank}{\levelthing}$.

Now suppose that $\gcd(\rank,\level)>1$. We have
\begin{gather*}
c(\rank,\level,1) 
= 
\frac{\rank}{\levelthing}\binom{\levelthing}{\rank} - \sum_{\substack{k\mid\gcd(\levelthing,\rank)\\k\neq 1}}c(\rank,\level,k)
=
\frac{\rank}{\levelthing}\binom{\levelthing}{\rank} - \sum_{\substack{k\mid\gcd(\levelthing,\rank)\\k\neq 1}}c(\rank/k,\levelthing/k,1).
\end{gather*}

Since $k\neq 1$ in the last sum, we can apply the induction hypothesis and we find that
\begin{gather*}
c(\rank,\level,1) 
=
\frac{\rank}{\levelthing}\binom{\levelthing}{\rank} - \sum_{\substack{k\mid\gcd(\levelthing,\rank)\\k\neq 1}}\frac{\rank}{\levelthing} \sum_{k'\mid \gcd\left(\rank/k,\levelthing/k\right)}\mu(k')\binom{\levelthing/kk'}{\rank/kk'}.
\end{gather*}

But
\begin{gather*}
\sum_{\substack{k\mid\gcd(\levelthing,\rank)}}\frac{\rank}{\levelthing} \sum_{k'\mid \gcd\left(\rank/k,\levelthing/k\right)}\mu(k')\binom{\levelthing/kk'}{\rank/kk'}
=
\frac{\rank}{\levelthing}\sum_{k\mid \gcd(\rank,\levelthing)}\sum_{k'\mid k}\mu\left(\frac{k}{k'}\right)\binom{\levelthing/k}{\rank/k}
\end{gather*}
and we conclude using the fact \ochanged{that} the M{\"o}bius function is the convolution inverse of the \ochangedd{identity}.\null\hfill$\square$

\subsection{Missing proofs for \autoref{S:Nhedral} and type D graphs}\label{appNhedral}

\subsubsection{\changedd{Proof of \autoref{L:Embedding}}}\label{Pf:Embedding}
We have two things to check: that the map is well-defined and that it is injective. The first is significantly more difficult than in \cite{MaMaMiTu-trihedral}.

\textit{Well-defined.} To see that the $\kltwo_{w_{i}}$ satisfy the defining relations of 
$\nhecke$ in \autoref{Eq:fundamental-relations} we first observe that $\kltwo_{w_{i}}\kltwo_{w_{i}}=\vnum{\rank}!\cdot\kltwo_{w_{i}}$ as follows from, for example, 
\cite[(2.8)]{El-thick-soergel-typea} or \cite[Proposition 10.5.2.(d)]{Bo-cells} since $\# W_{\mathtt{i}}=\rank!$. The second relation in \autoref{Eq:fundamental-relations} is much more difficult to prove and we use \autoref{QSH} to do so: under \autoref{QSH} going from 
$\kl_{k}$ to $\kl_{k+i}$ corresponds to tensoring with $\Ll_{\omega_{i}}$ and from $\kl_{k+i}$ 
to $\kl_{k+i+j}$ corresponds to tensoring with $\Ll_{\omega_{j}}$, and similarly, but reversed, for the other side of the equation. In particular, this equation in $\hecke$ holds since we have $\Ll_{\omega_{j}}\otimes\Ll_{\omega_{i}}\cong\Ll_{\omega_{i}}\otimes\Ll_{\omega_{j}}$.

\textit{Injective.} The same argument as in \cite[Proof of Lemma 3.2]{MaMaMiTu-trihedral} works.
\null\hfill$\square$

\subsubsection{\changedd{Proof of \autoref{T:simples} \textit{Part (a).}}}\label{Pf:simples-a}
We start with a lemma.

\begin{Lemma}\label{L:stabilizer}
Let $\bs\in\vanpar{\level}$. The stabilizer of $\bs$ is of order $\gcd\big(\{j | \Zfun_{j}(\bs)\neq 0\}\cup\{\rank\}\big)$. 
\end{Lemma}

\begin{proof}
Let $m$ be such that $Z_{j}(\bs) = 0$ if $m\nmid j$. By \autoref{L:OrbitsKoornwinder}, we have $\Zfun_{i}(\bs^{\turnsymbol})=\rankroot^{-i}\Zfun_{i}(\bs)$ for all $1\leq i \leq \rank$. This implies that $\Zfun_{i}\big(\bs^{\turnsymbol \rank/m}\big)=\Zfun_{i}(\bs)$, either because $\Zfun_{i}(\bs)=0$, if $m\nmid j$, or because $\rankroot^{-i \rank/m}=1$, otherwise. By \autoref{L:ZInjective}, we have $\bs^{\turnsymbol \rank/m} = \bs$.
\end{proof}

Since $\Zfun_{j}(\bs) = 0$ if $j\not\equiv 0 \bmod m$, the representation $\M(\bs)$ decomposes as $\Ll(\bs)_{0}\oplus\dots\oplus \Ll(\bs)_{m-1}$ where $\theta_{i},\theta_{i+m},\dots,\theta_{i+(\rank/m-1)m}$ act\ochangedd{s} on $\Ll(\bs)_{i}$ by the matrices
\begin{gather*}
\scalar
\begin{psmallmatrix}
\vnum{\rank} & \Zfun_{m}(\bs) & \Zfun_{2m}(\bs) & \dots & \Zfun_{(\rank/m-1)m}(\bs) \\
0 & 0 & 0 & \dots & 0\\
\vdots & \vdots & \vdots & \ddots & \vdots \\
0 & 0 & 0 & \dots & 0
\end{psmallmatrix},
\quad
\scalar
\begin{psmallmatrix}
0 & 0 & 0 & \dots & 0\\
\Zfun_{(\rank/m-1)m}(\bs) & \vnum{\rank} & \Zfun_{m}(\bs) & \dots & \Zfun_{(\rank/m-2)m}(\bs) \\
0 & 0 & 0 & \dots & 0\\
\vdots & \vdots & \vdots & \ddots & \vdots \\
0 & 0 & 0 & \dots & 0
\end{psmallmatrix},
\dots,\\
\scalar
\begin{psmallmatrix}
0 & 0 & \dots & 0 & 0\\
\vdots & \vdots & \ddots & \vdots & \vdots \\
0 & 0 & \dots & 0& 0\\
\Zfun_{m}(\bs) & \Zfun_{2m}(\bs) & \dots & \Zfun_{(\rank/m-1)}(\bs) & \vnum{\rank}
\end{psmallmatrix}
\end{gather*}
and by zero otherwise.

\begin{Lemma}\label{L:endoring}
The endomorphism ring of $\Ll(\bs)_{i}$ is one dimensional.
\end{Lemma}

\begin{proof}
Let $f$ be a matrix intertwiner for $\Ll(\bs)_{i}$. Since $\vnum{\rank}\neq 0$, we immediately obtain that $f$ must be a diagonal matrix $\mathrm{diag}(f_{0},\dots,f_{\rank/m-1})$. The relation $f\theta_{i+km} = \theta_{i+km}f$ implies that $f_{k}\Zfun_{rm}(\bs) = \Zfun_{rm}(\bs)f_{k+r}$ so that $f_{k}=f_{k+r}$ for all $0\leq k < \rank/m$ and $0 < r < \rank/m$ such that $\Zfun_{rm}(\bs)\neq 0$, where the indices are taken modulo $\rank/m$. By \autoref{L:stabilizer}, the set $\{1\leq r < \rank/m | \Zfun_{rm}(\bs)\neq 0\}$ generates $\Z/(\rank/m)\Z$. Therefore, we obtain that $f$ is diagonal.
\end{proof}

Since $\vnum{\rank}\neq 0$, the explicit form of the matrices above implies that the representations $\Ll(\bs)_{i}$ are pairwise nonisomorphic. They are moreover simple because their endomorphism rings are one dimensional by \autoref{L:endoring}.
\null\hfill$\square$

\subsubsection{\changedd{Proof of \autoref{T:simples} \textit{Part (b).}}}\label{Pf:simples-b} This follows from (a) and the following lemma.

\begin{Lemma}
The representations $\M(\bs)$ and $\M(\bs')$ are isomorphic if and only if $\bs$ and $\bs'$ are in the same orbit.
\end{Lemma}

\begin{proof}
We have two directions.

\textit{Case $\Leftarrow$.}
Let $j$ be such that $\bs'=\bs^{\turnsymbol j}$. By \autoref{L:OrbitsKoornwinder}, this implies that $\Zfun_{i}(\bs') = \rankroot^{-ij}\Zfun_{i}(\bs)$. An easy calculation shows that, for all $i\in \nodes$, we have
\begin{gather*}
\matrep{i}{\bs'} = \mathrm{diag}(1,\rankroot^{j},\dots,\rankroot^{(\rank-1)j})
\matrep{i}{\bs}
\mathrm{diag}(1,\rankroot^{j},\dots,\rankroot^{(\rank-1)j})^{-1}.
\end{gather*}
Therefore the representations $\M(\bs)$ and $\M(\bs')$ are isomorphic.

\textit{Case $\Rightarrow$.} As in the proof of \autoref{L:endoring}, an invertible matrix intertwiner $\M(\bs) \to \M(\bs')$ is necessarily a diagonal matrix $\mathrm{diag}(f_{0},\dots,f_{\rank-1})$. The relation $\matrep{i}{\bs}f = f\matrep{i}{\bs'}$ implies that $f_{i}\Zfun_{k}(\bs) = \Zfun_{k}(\bs')f_{i+k}$ for all $0\leq i < \rank$ and $1\leq k < \rank$. Since $f$ is invertible, we obtain that $\Zfun_{k}(\bs') = 0$ if and only if $\Zfun_{k}(\bs) = 0$. By \autoref{L:stabilizer}, the stabilizers of $\bs$ and $\bs'$ have the same order $m$.

If for all $1\leq k < \rank$ we have $\Zfun_{k}(\bs)=\Zfun_{k}(\bs')=0$, then, by \autoref{L:ZInjective}, we have $\bs=\bs'$. We then suppose that all Z-functions do not vanish at $\bs$ (and consequently at $\bs'$). By \autoref{L:stabilizer}, $m$ is a generator of the subgroup of $\Z/\rank\Z$ generated by the $1\leq j < \rank$ with $\Zfun_{j}(\bs)\neq 0$. Therefore, using the relations $f_{i}\Zfun_{k}(\bs) = \Zfun_{k}(\bs')f_{i+k}$ with $k$ such that $\Zfun_{k}(\bs)\neq 0 \neq \Zfun_{k}(\bs')$, one may show that the ratios $f_{i}/f_{i+m}$ are all equal. As \begin{gather*}
1=\frac{f_{0}}{f_{m}} \cdot \frac{f_{m}}{f_{2m}} \dots \frac{f_{(\rank/m-1)m}}{f_{0}} = \left(\frac{f_{0}}{f_{m}}\right)^{\rank/m},
\end{gather*}
the ratio $f_{0}/f_{m}$ is equal to $\rankroot^{-mj}$ for some $j$.

The relation $f_{0}\Zfun_{k}(\bs) = \Zfun_{k}(\bs')f_{k}$ gives then $\Zfun_{k}(\bs') = \rankroot^{-jk}\Zfun_{k}(\bs)$. Indeed, either $k$ is a multiple of $m$ and this follows from $f_{0} = \rankroot^{-mjr}f_{rm}$, either $k$ is not a multiple of $m$ and $\Zfun_{k}(\bs)=\Zfun_{k}(\bs')=0$. We finally obtain that $\Zfun_{k}(\bs') = \Zfun_{k}(\bs^{\turnsymbol j})$ for all $1\leq k < \rank$\changed{,} and then $\bs'=\bs^{\turnsymbol j}$ by \autoref{L:ZInjective}.
\null\hfill$\square$

\subsubsection{\changedd{Proof of \autoref{T:simples} \textit{Parts (c) and (d).}}}\label{Pf:simples-cd} The simple representations we found in (a) satisfy
\begin{gather*}
\dim_{\Cv}(\M_{0})^2 + \sum_{\bs \in \text{orbits}}\sum_{i=0}^{m-1}\dim_{\Cv}(\Ll(\bs)_{i})^2 = \dim_{\Cv}(\nhecke[\level]).
\end{gather*}
This follows from \autoref{T:KoornwinderVariety}.(c). Indeed, each orbit in $\vanset{\level}$ with stabilizer of order $m$ defines $m$ representations of dimension $\rank/m$ and there are $\frac{m}{\rank}\frac{\rank}{\levelthing}\sum_{k | \gcd(\rank/m, \levelthing/m)}\mu(k)\binom{\levelthing/mk}{\rank/mk}$ such orbits. Hence
\begin{gather*}
\sum_{\bs \in \text{orbits}}\sum_{i=0}^{m-1}\dim_{\Cv}(\Ll(\bs)_{i})^2 = \sum_{m | \gcd(\rank,\level)} m\frac{m}{\rank}\frac{\rank^2}{m^2}\frac{\rank}{\levelthing}\sum_{k | \gcd(\rank/m, \levelthing/m)}\mu(k)\binom{\levelthing/mk}{\rank/mk} 
\\
= \frac{\rank^2}{\levelthing} \sum_{m | \gcd(\rank,\level)}\sum_{k | \gcd(\rank/m, \levelthing/m)}\mu(k)\binom{\levelthing/mk}{\rank/mk} \stackrel{(\star)}{=} \frac{\rank^2}{\levelthing}\binom{\levelthing}{\rank} = \rank\tnumbers{\level}{\rank},
\end{gather*}
where the equality $(\star)$ follows from the convolution property of the M{\"o}bius function.
The result now follows \changed{since the sum of squares of dimensions of simple representations is equal to the dimension of the algebra, which implies semisimplicity}.
\end{proof}

\subsubsection{\changedd{Graphs of type D}}\label{Pf:typeD-graphs}

We also construct other examples of Ncolored \changed{graphs} through the fusion rules of module categories over $\slncatr$\changed{. The corresponding module categories are} constructed as categories of modules over an algebra in $\slncatr$. This is similar to the orbifold procedure of \cite{Ko-orbifold}.

Given a rank $\rank$\changed{, a level} $\level$ \text{and $1\leq i < \rank$}, the pointed subcategory generated by $\Ll_{\level\omega_{i}}$ has the fusion rules of $\Z/\rank\Z$. However, this subcategory is not a symmetric subcategory in general. Recall, for a finite group $G$ and a central element $z\in G$ with $z\neq 1$ and $z^2=1$, the braided fusion category $\catstuff{Rep}(G,z)$ as in \cite[Example 9.9.1.(3)]{EtGeNiOs-tensor-categories}. Let $g=\gcd(\rank,\level)$ and $p=N/g$. 

\begin{Lemma}
The pointed subcategory $\point$ generated by $\Ll_{\level\omega_{p}}$ is a symmetric subcategory of $\slncatr$. Moreover $\point$ is braided equivalent to $\catstuff{Rep}(\Z/g\Z,z)$ if ($\level/g$ and $p$ are odd and $g$ is even) and to $\catstuff{Rep}(\Z/g\Z)$ otherwise.
\end{Lemma}

\begin{proof}
The twist of the simple object $\Ll_{\level\omega_{i}}$ is given by $\rooty^{(\level\omega_{i},\level\omega_{i}+2\rho)} = \exp\big(\iunit\pi i (\rank-i)/\rank\big)$. Therefore, since $\Ll_{\level\omega_{pi}}\otimes\Ll_{\level\omega_{pj}}\simeq \Ll_{\level\omega_{p(i+j)}}$, we have $\br_{\Ll_{\level\omega_{pi}},\Ll_{\level\omega_{pj}}}=\theta_{\level\omega_{p(i+j)}}\theta_{\level\omega_{pi}}^{-1}\theta_{\level\omega_{pj}}^{-1} = \idmor$, and the subcategory $\point$ is symmetric.

By \cite[Corollary 9.9.25]{EtGeNiOs-tensor-categories}, it remains to compute the braiding of the generating object $\Ll_{\level\omega_{p}}$ with itself. Since the object $\Ll_{\level\omega_{p}}\otimes \Ll_{\level\omega_{p}}$ is simple, the braiding $\br_{\Ll_{\level\omega_{p}},\Ll_{\level\omega_{p}}}$ is a scalar multiple of the identity. By \cite[Exercise 8.10.15]{EtGeNiOs-tensor-categories}, since the invertible object $\Ll_{\level\omega_{p}}$ is of quantum dimension $1$, we have $\br_{\Ll_{\level\omega_{p}},\Ll_{\level\omega_{p}}} = \theta_{\level\omega_{p}}$. One may check that $\theta_{\level\omega_{p}} = (-1)^{\level p(g-1)/g}\idmor$, which is equal to $-1$ if and only if $\level/g$, $p$ and $g-1$ are all odd.
\end{proof}

\changed{We describe the classical construction of an algebra \changedd{in $\catstuff{Rep}(\Z/g\Z)$}, the case of $\catstuff{Rep}(\Z/g\Z,z)$ being similar. In the category $\catstuff{Rep}(\Z/g\Z)$ we consider the algebra of complex valued functions, which is an algebra in \changedd{this category}. \changedd{In the category $\catstuff{Rep}(\Z/g\Z,z)$, there is a similar algebra object, but the multiplication is twisted by a sign, through the action of $z$.} Fourier transform for abelian finite groups implies that these algebras are, as representations, isomorphic to the direct sum of simple representations.}

\changed{The category $\point$ \changedd{is} equivalent to $\catstuff{Rep}(\Z/g\Z)$ or $\catstuff{Rep}(\Z/g\Z,z)$\changedd{. W}e let $\alg$ be the algebra in $\point$ corresponding to the above algebra in $\catstuff{Rep}(\Z/g\Z)$ or $\catstuff{Rep}(\Z/g\Z,z)$.} 

\begin{Lemma}
In $\slncatr$, we have $\alg\simeq \Ll_{0}\oplus\Ll_{\level\omega_{p}}\oplus\dots\oplus\Ll_{\level\omega{(g-1)p}}$.
\end{Lemma}

\begin{proof}
\changed{The result  follows from the decomposition of the algebra in $\catstuff{Rep}(\Z/g\Z)$ and $\catstuff{Rep}(\Z/g\Z,z)$, since the $g$ simple objects of $\catstuff{Rep}(\Z/g\Z)$ or $\catstuff{Rep}(\Z/g\Z,z)$ are identified with the simple objects $\Ll_{0},\Ll_{\level\omega_{p}},\dots,\Ll_{\level\omega{(g-1)p}}$.}
\end{proof}

\begin{Remark}\label{R:webs}
For $\rank=3$, in \cite[Proposition 5.4.3]{MaMaMiTu-trihedral} the calculation of the algebra object $\alg$ was done using symmetric webs \cite{RoTu-symmetric-howe}, \cite{RoWa-sym-homology}, \cite{LaTu-gln-webs}. We instead use a corollary of Deligne\changedd{'s} theorem on Tannakian categories, see \cite[\changed{Corollary 9.9.25}]{EtGeNiOs-tensor-categories}.
\end{Remark}

We consider the category $\catstuff{M}$ of (right) modules over $\alg$ in $\slncatr$. In the \changed{case} $\point$ is equivalent to $\catstuff{Rep}(\Z/g\Z)$, the category $\catstuff{M}$ is the de-equivariantization of $\slncatr$, \changed{and} is a braided $\Z/g\Z$-crossed category.

\begin{Lemma}
The category $\catstuff{M}$ is a finite semisimple module category over $\slncatr$.
\end{Lemma}

\begin{proof}
The only nontrivial statement is the semisimplicity, which follows from the separability of $\alg$, by \cite[Proposition 7.8.30]{EtGeNiOs-tensor-categories}. By \autoref{R:webs}, separability is equivalent to the digon scalar to be invertible, which is immediate from the definition of symmetric webs. \changed{Here, we assume familiarity with webs.}
\end{proof}

There is $\Z/\rank\Z$-grading on $\catstuff{M}$. The simple objects of $\catstuff{M}$ are all obtained as summands of free objects $\alg\otimes \Ll_{\bm}$.  Since all summands of $\alg$ are of color $0$, $\alg\otimes \Ll_{\bm}$ is of color $\chi_{c}(\Ll_{\bm})$. This grading is moreover compatible with the grading of $\slncatr$.

We define the \changed{\emph{graph of type D}} of rank $\rank$ and level $\level$ as the graph $\Gamma$ with set of vertices $\simples{\M}$, colored with the above $\Z/\rank\Z$-grading, and the vertices $\obstuff{X}$ and $\obstuff{Y}$ are adjacent if and only if $\obstuff{Y}$ is a summand of $\bigoplus_{i=1}^{\rank-1}\Ll_{\omega_{i}}\otimes \obstuff{X}$. The graph $\Gamma_{i}$ is then the fusion graph of the object $\Ll_{\omega_{i}}$.

\subsubsection{\changedd{Proof of \autoref{T:integralTlevelA} \textit{Part (a)}}}\label{Pf:integralTlevelA-a}
For type $A$, this is explained in \cite[Section 1.3, Example 2]{Zu-gen-dynkin-diagrams}. We give a self-contained proof.

\begin{Lemma}
For $\Gamma$ of type $A$ of rank $\rank$ and level $\level$, the matrices $(A(\Gamma_{1}),\dots,A(\Gamma_{\rank-1}))$ can be simultaneously diagonalized and their joint spectrum is the Koornwinder variety $\vanset{\level}$.
\end{Lemma}

\begin{proof}
The matrix $A(\Gamma_{i})$ is the fusion matrix of $\Ll_{\omega_{i}}$ in the category $\slncatr$. The $S$-matrix diagonalizes the fusion rules \cite[Corollary 8.14.5]{EtGeNiOs-tensor-categories} and the joint spectrum of $(A(\Gamma_{1}),\dots,A(\Gamma_{\rank-1}))$ is:
\begin{gather*}
\Big\{\big(S^{\sln[]}_{\bk,\omega_{1}}/S^{\sln[]}_{\bk,0},\dots,S^{\sln[]}_{\bk,\omega_{\rank-1}}/S^{\sln[]}_{\bk,0}\big) \Big| \bk\in X^{+}(\level)\Big\}.
\end{gather*}
A direct application of Weyl character formula shows that
\begin{gather*}
\frac{S^{\sln[]}_{\bk,\omega_{i}}}{S^{\sln[]}_{\bk,0}}
=
\sum_{j}\rooty^{2(\bk+\rho,w_{j}^i)}
=
\Zfun_{i}\left(\frac{2\pi}{\level+\rank}(\bk+\rho)\right).
\end{gather*}
We conclude using \autoref{L:KoornwinderVarietyZ} and \autoref{T:KoornwinderVariety}.
\end{proof}

By definition of the Koornwinder variety and \autoref{L:integralTlevel}, $\M(\Gamma)$ is a $\nhecke[\level]$-representation for $\Gamma$ of type $A$.
\null\hfill$\square$

\subsection{Missing proofs for \autoref{S:AsymptoticCategory}}\label{appasymptotic}

\subsubsection{\changedd{Proof of \autoref{T:DrinfeldCenter}
\textit{Part (c).}}}\label{Pf:DrinfeldCenter-c} We will first compute the rank of $\dcenter{\slncatr}_{0}$ and then the number of simple objects with a stabilizer of a given order. We will denote by $\Ll_{\bm,\bk}$ the simple object $\Ll_{\bm}\boxtimes\Ll_{\bk}$. In these lemmas we use Euler's totient function $\varphi$ and M{\"o}bius function $\mu$.

\begin{Lemma}\label{L:SimpleCenterZero}
The rank of $\dcenter{\slncatr}_{0}$ is
\begin{gather*}
\frac{\rank}{\levelthing^{2}}\sum_{k\mid \gcd(\levelthing,\rank)}\varphi(k)\binom{\levelthing/k}{\rank/k}^{2}.
\end{gather*}
\end{Lemma}

\begin{proof}
The simple objects are indexed by pairs $(\bm,\bk)\in \ochangedd{(}X^{+}(\level)\ochangedd{)^{2}}$ such that $\chi_{c}(\Ll_{\bm}) = \chi_{c}(\Ll_{\bk})$. Let $V$ be a complex vector space with basis $(v_{0},\ldots,v_{\rank{-}1})$ equipped with the endomorphism $\Omega$ defined by $\Omega(v_{i}) = \omega^{i}v_{i}$, where $\omega$ is a primitive $\rank$th root of unity. This endows $V$ with an action of $\Z/\rank\Z$\changed{,} and then also the symmetric algebra $S(V\times V^{*}) = S(V)\otimes S(V^{*})$. \ochangedd{If we denote by $(w_{0},\ldots,w_{\rank-1})$ the dual basis of $(v_{0},\ldots,v_{\rank-1})$, then $(v_{\bm,\bk})_{(\bm,\bk)\in (X^{+}(\level))^{2}}$, where
\begin{equation*}
v_{\bm,\bk} = v_0^{\level-\bsum\bm}v_1^{m_1}\dots v_{\rank-1}^{\bm_{\rank-1}}\otimes w_{0}^{\level-\bsum\bk}w_{1}^{\bk_{1}}\dots w_{\rank-1}^{\bk_{\rank-1}},
\end{equation*}
is a basis of $S^{\level}(V)\otimes S^{\level}(V^*)$, on which $\Omega$ acts diagonally: $\Omega(v_{\bm,\bk}) = \omega^{\chi_{c}(\Ll_{\bm})-\chi_{c}(\Ll_{\bk})}v_{\bm,\bk}$. Therefore} the rank of $\dcenter{\slncatr}_{0}$ is equal to the dimension of $(S^{\level}(V)\otimes S^{\level}(V^*))^{\Z/\rank\Z}$. Using a multigraded version of Molien's formula, one finds that the rank of $\dcenter{\slncatr}_{0}$ is the coefficient of $u^{\level}v^{\level}$ in
\begin{gather*}
\frac{1}{\rank}\sum_{i\in \Z/\rank\Z}\frac{1}{\det_{V}(1-\Omega^{i}u)\det_{V^*}(1-\Omega^{i}v)}.
\end{gather*}

Since $\Omega$ is diagonal in the basis $(v_{0},\ldots,v_{\rank{-}1})$, we have $\det_{V}(1-\Omega^{i}u) = \prod_{j=0}^{\rank{-}1}(1-\omega^{ij}u) = (1-u^{k_{i}})^{n/k_{i}}$, where $k_{i}$ is the order of the root $\omega^{i}$. Similarly, $\det_{V^{*}}(1-\Omega^{i}v) = (1-v^{k_{i}})^{n/k_{i}}$. Since there are exactly $\varphi(k)$ powers of $\omega$ of order $k$, one has
\begin{align*}
\frac{1}{\rank}\sum_{i\in \Z/\rank\Z}\frac{1}{\det_{V}(1-\Omega^{i}u)\det_{V^*}(1-\Omega^{i}v)} 
&=
\frac{1}{\rank}\sum_{k\mid\rank}\frac{\varphi(k)}{(1-u^{k})^{\rank/k}(1-v^{k})^{\rank/k}}\\
&=
\frac{1}{\rank}\sum_{k\mid \rank} \varphi(k) \sum_{i,j \geq 0} \binom{\rank/k+i-1}{\rank/k-1}\binom{\rank/k+j-1}{\rank/k-1}u^{ki}v^{kj}.
\end{align*}

Only the terms with $k\mid\level$ will contribute to the coefficient of $u^{\level}v^{\level}$ and we finally obtain that the rank of $\dcenter{\slncatr}_{0}$ is
\begin{gather*}
\frac{1}{\rank}\sum_{k\mid\gcd(\level,\rank)}\varphi(k)\binom{\rank/k+\level/k-1}{\rank/k-1}^{2}
=
\frac{\rank}{\levelthing^{2}}\sum_{k\mid \gcd(\levelthing,\rank)}\varphi(k)\binom{\levelthing/k}{\rank/k}^{2}\changedd{.}
\end{gather*}
the last equality is obtained using $\level+\rank = \levelthing$ \changedd{and $\binom{a-1}{b-1}=\frac{b}{a}\binom{a}{b}$}.
\end{proof}

\begin{Lemma}\label{L:ObjectsStabilizer}
Let $m\in\Z$. 
\changedd{If $m\nmid\gcd(\rank,\level)$, then there is no simple objects of $\dcenter{\slncatr}_{0}$ with stabilizer of order $m$ under $\Ll_{\level\omega_{1},\level\omega_{1}}\otimes\placeholder$.} 
If $m\mid\gcd(\rank,\level)$, then the number of simple objects of $\dcenter{\slncatr}_{0}$ with stabilizer of order $m$ under $\Ll_{\level\omega_{1},\level\omega_{1}}\otimes\placeholder$ is
\begin{gather*}
\ochangedd{c_{\catstuff{Z}}}(\rank,\levelthing,m)=\frac{m\rank}{\levelthing^{2}}\sum_{k\mid \gcd\left(\rank/m,\levelthing/m\right)}\mu(k)\binom{\levelthing/mk}{\rank/mk}^{2}\changedd{.} 
\end{gather*}
\end{Lemma}

\begin{proof}
The count is similar to the count of \autoref{T:KoornwinderVariety}.(c). The base case is dealt by \autoref{L:SimpleCenterZero}, and, at the end of the inductive step, we use the fact Euler's totient function is the convolution of the M{\"o}bius function and of the identity function.
\end{proof}

Each orbit of $X^{+}(\level)^{2}$ \changed{with} stabilizer of order $m$ will \changed{contribute} $m$ simple objects in the modular closure of $\dcenter{\slncatr}_{0}$, and there are $\tfrac{m}{\rank}\ochangedd{c_{\catstuff{Z}}}(\rank,\levelthing,m)$ such orbits. Hence, the rank of $\dcenter{\slncatr_{0}}$ is
\begin{align*}
\sum_{m\mid\gcd(\levelthing,\rank)}\frac{m^{2}}{\rank}\ochangedd{c_{\catstuff{Z}}}(\rank,\levelthing,m) 
&= 
\frac{1}{\levelthing^{2}}\sum_{m\mid\gcd(\levelthing,\rank)}\sum_{k\mid\gcd(\levelthing/m,\rank/m)}m^{3}\mu(k)\binom{\levelthing/mk}{\rank/mk}^{2}  
\\
&=
\frac{1}{\levelthing^{2}}\sum_{m\mid\gcd(\levelthing,\rank)}\sum_{k\mid m} k^3\mu\left(\frac{m}{k}\right)\binom{\levelthing/m}{\rank/m}^{2}.
\end{align*}
We obtain the formula in \changed{\autoref{T:DrinfeldCenter}}(c) since $J_{3}$ is the convolution of the M{\"o}bius function $\mu$ and the cube function.\null\hfill$\square$

\subsubsection{\changedd{Proof of \autoref{T:DrinfeldCenter}
\textit{Part (e).}}}\label{Pf:DrinfeldCenter-e}  We will use the description of the modular closure of \ochangedd{the category} $\catstuff{R}=\left(\slncatr\boxtimes\slncatr^{\rev}\right)_{0}$ as in \cite[Definition 3.12]{Mu-galois-theory}. Recall that it is obtained as the idempotent completion of the category where the objects are the same as in $\catstuff{R}$, but morphisms between $X$ and $Y$ are $\bigoplus_{i\in\Z/\rank\Z}\Hom_{\catstuff{R}}(X,\Ll_{\level\omega_{i},\level\omega_{i}}\otimes Y)$. \changed{Let us} denote by $\gamma$ an \changed{endomorphism of $\Ll_{\bm,\bk}$ in the modular closure} different from the identity and of order $\stabilizer{\bm,\bk}$\changed{. Then} the primitive idempotents of $\End(L_{\bm,\bk})$ in the modular closure are $p_{j}=\tfrac{1}{\stabilizer{\bm,\bk}}\sum_{k=1}^{\stabilizer{\bm,\bk}}\xi^{jk}\gamma^{k}$ for $j\in\Z/\rank\Z$, where $\xi$ is a primitive $\stabilizer{\bm,\bk}$th root of unity. Using similar notations for the object $\Ll_{\bm^{\prime},\bk^{\prime}}$, we have
\begin{align*}
S_{(\bm,\bk,i),(\bm^{\prime},\bk^{\prime},j)} 
&= \Tr(p_{i}\otimes p^{\prime}_{j} \circ \beta_{\Ll_{\bm,\bk},\Ll_{\bm,\bk}}\circ \beta_{\Ll_{\bm,\bk},\Ll_{\bm,\bk}}) \\
&= \frac{1}{\stabilizer{\bm,\bk}}\sum_{k=1}^{\stabilizer{\bm,\bk}}\xi^{ik}\Tr(\gamma^{k}\otimes p^{\prime}_{j^{\prime}}\circ\beta_{\Ll_{\bm^{\prime},\bk^{\prime}},\Ll_{\bm,\bk}}\circ \beta_{\Ll_{\bm,\bk},\Ll_{\bm^{\prime},\bk^{\prime}}}).
\end{align*}
Therefore, summing over $i$, we obtain
\begin{gather*}
\sum_{i=1}^{\stabilizer{\bk,\bm}}S_{(\bm,\bk,i),(\bm^{\prime},\bk^{\prime},i^{\prime})} 
=
\Tr(\idmor\otimes p^{\prime}_{j^{\prime}}\circ\beta_{\Ll_{\bm^{\prime},\bk^{\prime}},\Ll_{\bm,\bk}}\circ \beta_{\Ll_{\bm,\bk},\Ll_{\bm^{\prime},\bk^{\prime}}})
=
\frac{1}{\stabilizer{\bm,\bk}}\Tr(\beta_{\Ll_{\bm^{\prime},\bk^{\prime}},\Ll_{\bm,\bk}}\circ \beta_{\Ll_{\bm,\bk},\Ll_{\bm^{\prime},\bk^{\prime}}}),
\end{gather*}
since the other terms in the sum defining $p^{\prime}_{j^{\prime}}$ are of trace $0$, as they are traces of morphisms between different objects in $\catstuff{R}$. 

Now suppose that $\rank$ is prime. Then \autoref{L:ObjectsStabilizer} shows that there exists a unique pair $(\bm,\bk)\in X^{+}(\level)^{2}$ with $\chi_{c}(\Ll_{\bm})=\chi_{c}(\Ll_{\bk})$ and $\stabilizer{\bm,\bk}\neq 1$ if and only if $\level\equiv 0\mod\rank$, and one easily checks that $\bm=\bk=(\level/\rank,\ldots,\level/\rank)$ provides such a pair. A similar argument as above shows that $S_{(\bm,\bk,i),(\bm,\bk,j)}=S_{(\bm,\bk,i^{\prime}),(\bm,\bk,j^{\prime})}$ if $i-j\equiv i^{\prime}-j^{\prime} \mod \rank$. 

We now use the relation $\tau^{-}T^{-1}ST^{-1}=STS$, where $\tau^{-}$ is the Gauss sum as in \cite[Definition 8.15.1]{EtGeNiOs-tensor-categories}. This relation is satisfied since the modular closure is a modular category. The entry $((\bm,\bk,0),(\bm,\bk,i))$ of this relation gives
\begin{align*}
\tau^{-}\theta_{\bm,\bk}^{2}S_{(\bm,\bk,0),(\bm,\bk,i)} =& \sum_{(\bm^{\prime},\bk^{\prime})\neq (\bm,\bk)}\theta_{\bm^{\prime},\bk^{\prime}}^{-1}S_{(\bm,\bk,0),(\bm^{\prime},\bk^{\prime})}S_{(\bm^{\prime},\bk^{\prime}),(\bm,\bk,0)} \\
&+ \theta_{\bm,\bk}^{-1}\sum_{j=1}^{\rank}S_{(\bm,\bk,0),(\bm,\bk,j)}S_{(\bm,\bk,j),(\bm,\bk,i)}.
\end{align*}
Let $k\not\equiv 0 \mod \rank$. Multiplying by $\xi^{ik}$ and summing then gives
\begin{align*}
\tau^{-}\theta_{\bm\bk}^{3}\sum_{i=1}^{\rank}\xi^{ik}S_{(\bm,\bk,0),(\bm,\bk,i)} 
&= \sum_{i,j=1}^{\rank} \xi^{ik} S_{(\bm,\bk,0),(\bm,\bk,j)}S_{(\bm,\bk,j),(\bm,\bk,i)} \\
&= \sum_{i,j=1}^{\rank} \xi^{ik} S_{(\bm,\bk,0),(\bm,\bk,j)}S_{(\bm,\bk,0),(\bm,\bk,j-i)} 
= \big(\sum_{i=1}^{\rank}\xi^{ik}S_{(\bm,\bk,0),(\bm,\bk,i)}\big)^{2}.
\end{align*}
But $\sum_{i=1}^{\rank}\xi^{ik}S_{(\bm,\bk,0),(\bm,\bk,i)}\neq 0$ since a similar calculation shows that
\begin{gather*}
\sum_{i=1}^{\rank}\xi^{ik}S_{(\bm,\bk,0),(\bm,\bk,i)} = \sum_{i=1}^{\rank}(S^{2})_{(\bm,\bk,0),(\bm,\bk,i)}
\end{gather*}
and that $S^{2}$ is, up to a nonzero constant, the permutation matrix given by the duality on simple objects.

We finally obtain that
\begin{gather*}
S_{(\bm,\bk,i),(\bm,\bk,j)} 
=\frac{1}{\rank^{2}}
\begin{cases}
S_{(\bm,\bk),(\bm,\bk)} + \rank(\rank{-}1)\theta_{\bm,\bk}^{3}\tau^{-} & \text{if }i=j,
\\
S_{(\bm,\bk),(\bm,\bk)} - \rank\theta_{\bm,\bk}^{3}\tau^{-} & \text{if }i\neq j.
\end{cases}
\end{gather*}
It \changedd{is} a calculation to show that the Gauss sum $\tau^{-}$ is equal to \changed{$\dim_{c}\,\slncatr/\rank$}. \changed{First, one  shows that the Gauss sum $\tau^{-}$ is equal to the Gauss sum $\tau^{-}(\catstuff{R})$ of $\catstuff{R}$. Then, using the grading, one shows that $\rank\tau^{-}(\catstuff{R}) = \tau^{-}\left(\slncatr\boxtimes\slncatr^{\rev}\right)$. Finally, one concludes using \cite[Proposition 8.15.4]{EtGeNiOs-tensor-categories}.}\null\hfill$\square$

\subsubsection{\changedd{Definition of the action}}\label{Pf:M-symbols}

Now we come to \autoref{SS:FamilyUnipotent}, so the reader might want to recall the various terminology. We define an action of $\Z/\levelthing\Z$ on $\Psi(Y,\pi)$ as follows. Given $f\in \Psi(Y,\pi)$ and $k\in \Z/\levelthing\Z$, we define $f+k$ as the unique function $Y\rightarrow \{0,\ldots,\levelthing-1\}$ such that $(f+k)(y) \equiv f(y)+k \ \bmod\levelthing$ for all $y\in Y$\changed{. Since $f+k$ is not necessarily increasing on $\pi^{-1}(i)$ for all $i\in\N$, we denote by} $f^{+k}$ the unique function in $\Psi(Y,\pi)$ such that $f^{+k}(\pi^{-1}(i)) = (f+k)(\pi^{-1}(i))$ for all $i\in\N$. Given $f\in\Psi(Y,\pi)$, we denote by $s(f)$ the cardinal of the stabilizer of $f$ for \changed{the} action of $\Z/\levelthing\Z$ \changed{induced by $f\mapsto f^{+1}$}. The unipotent characters lying in the family $\family$ are then indexed by pairs $([f],\changedd{s})$ where $[f]$ is an equivalence class of $\Psi(Y,\pi)$ under the action of $\Z/\levelthing\Z$ and $1\leq \changedd{s} \leq s(f)$. 

\changed{For $f\in \Psi(Y,\pi)$, let $\bar{f} \colon \{1,\dots,\rank\} \to \{0,\dots,\levelthing-1\}$ be the only increasing function such that $\{\bar{f}(1),\ldots,\bar{f}(\rank)\}$ is the complement of $f(\pi^{-1}\changedd{(0)})$ in $\{0,\ldots,\levelthing-1\}$.} \ochangedd{Then we have $\sum_{i=1}^{\rank}(f(i)-\bar{f}(i))\equiv 0 \mod \levelthing$.}

\begin{Example}\label{E:symbol-1}
\changedd{We take $\level=2$ and $\rank=3$, so that $\levelthing=5$, and $f\in\Psi(Y,\pi)$ given by $f(1)=0$, $f(2)=1$, $f(3)=2$, $f(4)=0$ and $f(5)=2$\changeddd{. This corresponds to Malle's $5$-symbol,}
\begin{equation*}
\begin{pmatrix}
0 & 1\\ 1 & \\ 0 & 1 \\ \placeholder & \\ \placeholder & 
\end{pmatrix}
\end{equation*}
\changeddd{see \autoref{R:Symbols} for the construction. The $1$'s in this symbol appear in the rows $f(1),f(2)$ and $f(3)=f(\rank)$, whereas the $0$'s appear in the rows $f(4)=f(\rank+1)$ and $f(5) = f(\rank+\level) = f(\levelthing)$, and $\bar{f}$ records the rows which do not have $0$'s.}
}

\changeddd{Finally, the functions $f^{+1}$ and $f^{+3}$ satisfy $f^{+1}(1)=1, f^{+1}(2)=2, f^{+1}(3)=3, f^{+1}(4)=1$ and $f^{+1}(5)=3$ and $f^{+3}(1)=0,f^{+3}(2)=3,f^{+3}(3)=4,f^{+3}(4)=0$ and $f^{+3}(5)=3$. They respectively correspond to the symbols}
\begin{equation*}
\changeddd{\begin{pmatrix}
\placeholder & \\ 0 & 1 \\ 1 & \\ 0 & 1 \\ \placeholder
\end{pmatrix}
\quad\text{and}\quad
\begin{pmatrix}
0 & 1 \\ \placeholder & \\ \placeholder & \\ 0 & 1 \\ 1 & 
\end{pmatrix}.}
\end{equation*}
\changeddd{The symbols corresponding to $f^{+1}$ and $f^{+3}$ are obtained from the symbol corresponding to $f$ via a cyclic shift.}
\end{Example}

\subsubsection{\changedd{Proof of \autoref{T:Comparison}
\textit{Part (a).}}} \label{Pf:Comparison-a} The computation of the cardinal of $\family$ is similar to \autoref{T:DrinfeldCenter}(c). 

\changedd{\begin{Lemma}\label{L:card-psi}
The cardinal of $\Psi(Y,\pi)$ is $\frac{1}{\levelthing}\sum_{k\mid \gcd(\rank,\levelthing)} \binom{\levelthing/k}{\rank/k}^{2}$.
\end{Lemma}}

\changedd{\begin{proof}
The proof is similar to \autoref{L:SimpleCenterZero}. We consider $V$ a $\mathbb{C}$-vector space of dimension $\levelthing$ equipped with a basis $(v_{0},\ldots,v_{\levelthing-1})$. Let $\omega$ be a primitive $\levelthing$-th root of unity, and define $\Omega$ as the endomorphism of $V$ such that $\Omega(v_{i})=\omega^{i}v_{i}$. This endows $V$ with an action of $\Z/\levelthing\Z$, and also the exterior power $\Lambda(V\times V^{*}) = \Lambda(V)\otimes \Lambda(V^{*})$. We denote by $(w_{0},\ldots,w_{\levelthing-1})$ the dual basis of $(v_{0},\ldots,v_{\levelthing-1})$. Given a function $f\colon Y \to \{0,\ldots,\levelthing-1\}$ strictly increasing on $\pi^{-1}(1)=\{1,\ldots,\rank\}$ and on $\pi^{-1}(0) = \{\rank+1,\ldots,\levelthing\}$, we set
\begin{equation*}
v_{f} = (v_{f(1)}\wedge\dots\wedge v_{f(\rank)})\otimes(w_{\bar{f}(1)}\wedge\dots\wedge w_{\bar{f}(\rank)})\in \Lambda^{\rank}(V)\otimes\Lambda^{\rank}(V^{*}).
\end{equation*}
Such vectors form a basis of $\Lambda^{\rank}(V)\otimes\Lambda^{\rank}(V^{*})$ on which $\Omega$ acts diagonally: $\Omega(v_{f}) = \omega^{\sum_{i=1}^{\rank}(f(i)-\bar{f}(i))}v_{f}$. Therefore, the cardinality of $\Psi(Y,\pi)$ is equal to the dimension of $(\Lambda^{\rank}(V)\otimes\Lambda^{\rank}(V^{*}))^{\Z/\levelthing\Z}$, and, via a Molien-type argument, is the $u^{\rank}v^{\rank}$ coefficient in
\begin{equation*}
\frac{1}{\levelthing}\sum_{i\in \Z/\levelthing\Z}\det\nolimits_{V}(1+\Omega^{i}u)\det\nolimits_{V^{*}}(1+\Omega^{i}v).
\end{equation*}
Using arguments similar to \autoref{L:SimpleCenterZero}, we find that the $u^{\rank}v^{\rank}$ coefficient of the above polynomial is $\frac{1}{\levelthing}\sum_{k\mid\gcd(\rank,\levelthing)}\varphi(k)\binom{\levelthing/k}{\rank/k}^{2}$.
\end{proof}}

\changedd{\begin{Lemma}\label{L:NumberOrbits}
Let $m\in\Z$. If $m\nmid\gcd(\rank,\levelthing)$, then there is no $f\in \Psi(Y,\pi)$ with $s(f)=m$.
If $m\mid\gcd(\rank,\levelthing)$, then the number of $f\in\Psi(Y,\pi)$ with $s(f)=m$ is
\begin{gather*}
c_{\Psi}(\rank,\levelthing,m)=\frac{m}{\levelthing}\sum_{k\mid \gcd\left(\rank/m,\levelthing/m\right)}\mu(k)\binom{\levelthing/mk}{\rank/mk}^{2}. 
\end{gather*}
\end{Lemma}}

\changedd{\begin{proof}
As for \autoref{L:ObjectsStabilizer}, the proof is similar to proof of \autoref{T:KoornwinderVariety}.(c), the base case being dealt by \autoref{L:card-psi}.
\end{proof}}

\changedd{Each orbit of $\Psi(Y,\pi)$ with stabilizer of order $m$ will contribute $m$ unipotent characters in $\family$, and there are $\frac{m}{\levelthing}c_{\Psi}(\rank,\levelthing,m)$ such orbits. Hence, the cardinality of $\family$ is
\begin{align*}
\sum_{m\mid\gcd(\levelthing,\rank)}\frac{m^{2}}{\levelthing}c_{\Psi}(\rank,\levelthing,m) &= \frac{1}{\levelthing^{2}}\sum_{m\mid\gcd(\levelthing,\rank)}\sum_{k\mid \gcd\left(\rank/m,\levelthing/m\right)}m^{3}\mu(k)\binom{\levelthing/mk}{\rank/mk}^{2}\\
&= \frac{1}{\levelthing^{2}}\sum_{m\mid\gcd(\levelthing,\rank)}\sum_{k\mid m} k^3\mu\left(\frac{m}{k}\right)\binom{\levelthing/m}{\rank/m}^{2}.
\end{align*}
Since $J_{3}$ is the convolution of the M{\"o}bius function $\mu$ and the cube function, we obtain 
\begin{equation*}
\lvert\family\rvert = \frac{1}{\levelthing}\sum_{k\mid\gcd(\rank,\levelthing)}J_{3}(k)\binom{\levelthing/k}{\rank/k}^{2}.
\end{equation*}
By \autoref{T:DrinfeldCenter}.(c), this cardinality is the rank of $\dcenter{\verg{\levelthing}{\rank}}$.}
\null\hfill$\square$

\subsubsection{\changedd{ The map \texorpdfstring{$\bar{\iota}$}{iota}}}\label{Pf:iota}

We define a map $\iota\colon \Psi(Y,\pi) \to (X^{+}(\level))^{2}, f\mapsto (\bm_{f},\bk_{f})$ where
\begin{gather*}
\bm_{f}=(f(2)-f(1)-1,\ldots,f(\rank)-f(\rank-1)-1)^{\turnsymbol r_{f}}\text{ and } \bk_{f}=(\bar{f}(2)-\bar{f}(1)-1,\ldots,\bar{f}(\rank)-\bar{f}(\rank-1)-1),
\end{gather*}
with $r_{f}\in\Z$ being the unique integer such that $\sum_{i=1}^{\rank}(f(i)-\bar{f}(i))=r_{f}\ochangedd{\levelthing}$. \ochangedd{The appearance of a rotation by $r_f$ is necessary to make \autoref{L:bij_same_color} true. Note that, going back to the definition of $\placeholder{}^{\turnsymbol}$, using the weight notation (with the convention that $\omega_{0}=0$), and considering all indices and arguments modulo $\rank$, we have
\begin{equation}
\label{eq:m_f-expanded}
\bm_{f} = \sum_{i=1}^{\rank-1}\big(f(i+1-r_f)-f(i-r_{f})-1\big)\omega_{i} + \levelthing \omega_{r_{f}}.
\end{equation}}

\begin{Example}\label{E:symbol-2}
\changeddd{We continue \autoref{E:symbol-1}, and consider the same function $f\in \Psi(Y,\pi)$. We have $r_f=-1$, and the corresponding weights are $\bm_{f}=(0,0)^{\turnsymbol(-1)}=(0,2)$ and $\bk_{f}=(1,0)$. Both are of color $1\in \Z/\rank\Z$, see \autoref{L:bij_same_color}.}

\changeddd{The reader may verify that $r_{f^{+1}}=0$ and $r_{f^{+3}}=0$, that $\bm_{f^{+1}} = (\bm_{f})^{\turnsymbol}$ and $\bk_{f^{+1}} = (\bk_{f})^{\turnsymbol}$, and that $\bm_{f^{+3}} = (\bm_{f})^{\turnsymbol 2}$ and $\bk_{f^{+3}} = (\bk_{f})^{\turnsymbol 2}$. This behavior of $\iota$ under the action of $\Z/\levelthing\Z$ is explained in the proof of \autoref{L:iota_bijection}}
\end{Example}

We now show that $\iota$ induces an isomorphism between $\family$ and $\simples{\dcenter{\verg{\levelthing}{\rank}}}$.

\begin{Lemma}\label{L:bij_same_color}
For any $f\in \Psi(Y,\pi)$, we have $\chi_{c}(\Ll_{\bm_{f}}) = \chi_{c}(\Ll_{\bk_{f}})\in \Z/\rank\Z$.
\end{Lemma}

\begin{proof}
\changedd{A straightforward calculation shows that $\chi_{c}(\Ll_{\bk_{f}}) = -\sum_{i=1}^{\rank}\overline{f}(i)-\binom{\rank}{2}$. Since $\placeholder{}^{\turnsymbol}$ adds $\level$ to the color, we obtain that $\chi_{c}(\Ll_{\bm_{f}}) = -\sum_{i=1}^{\rank}f(i)-\binom{\rank}{2}+r_{f}\level$. Therefore $\chi_{c}(\Ll_{\bk_{f}})-\chi_{c}(\Ll_{\bm_{f}}) = r_f(\levelthing - \level) = 0$ since $\levelthing=\level+\rank$.}
\end{proof}

\begin{Lemma}\label{L:iota_bijection}
The map $\iota$ induces a bijection \ochangedd{$\bar{\iota}$} between the orbits of $\Psi(Y,\pi)$ under the action of $\Z/\levelthing\Z$ and the orbits of $\{(\bm,\bk)\in (X^{+}(e))^{2}\mid \chi_{c}(\Ll_{\bm})=\chi_{c}(\Ll_{\bk})\}$ under the action of $\Z/\rank\Z$.
\end{Lemma}

\begin{proof}
\changedd{We first study the behavior of $\iota$ under the action of $\Z/\levelthing\Z$. For $f\in \Psi(Y,\pi)$, we show that 
\begin{equation*}
\iota(f^{+1}) = 
\begin{cases}
\iota(f) & \text{if } \bar{f}(\rank)\neq \levelthing-1,\\
\iota(f)^{\turnsymbol} & \text{if } \bar{f}(\rank) = \levelthing-1.
\end{cases}
\end{equation*}
We distinguish four different cases depending on whether $f(\rank)$ and $\bar{f}(\rank)$ are equal to $\levelthing-1$. If $f(\rank)=\levelthing-1$ and $\bar{f}(\changedd{\rank})\neq \levelthing - 1$, then $f^{\changedd{+1}}(1) = 0$, $f^{\changedd{+1}}(i) = f(i-1)+1$ for all $1< i \leq \rank$ and $\overline{f^{\changedd{+1}}}(i) = \bar{f}(i)+1$ for all $1\leq i \leq \rank$. Therefore $r_{f^{\changedd{+1}}}=r_{f}-1$ and $\iota(f^{\changedd{+1}})=\iota(f)$. The three other cases are similar.} \changedd{Therefore, $\iota$ induces a well-defined map $\bar{\iota}$ from the orbits of of $\Psi(Y,\pi)$ under the action of $\Z/\levelthing\Z$ to the orbits of $\{(\bm,\bk)\in (X^{+}(e))^{2}\mid \chi_{c}(\Ll_{\bm})=\chi_{c}(\Ll_{\bk})\}$ under the action of $\Z/\rank\Z$.}

We now prove that if $f,g\in \Psi(Y,\pi)$ are such that $\iota(f)$ and $\iota(g)$ are in the same orbit of $(X^{+}(e))^{2}$ under the action of $\Z/\rank\Z$, then $f$ and $g$ are in the same orbit of $\Psi(Y,\pi)$ under the action of $\Z/\levelthing\Z$. \changedd{This would imply that the above induced map $\bar{\iota}$ is injective.} Let $1\leq l \leq \rank$ such that $\iota(f)=\iota(g)^{\turnsymbol l}$, and $0\leq k < \levelthing$ such that $\bar{g}(\rank-l+1) = \levelthing-k$. Then \changedd{one may check that} $\iota(g^{+k})=\iota(g)^{\turnsymbol l}$ and we may and will suppose that $\iota(f)=\iota(g)$. \ochangedd{The above behavior of $\iota$ also shows that $\iota(f^{+(\levelthing-\bar{f}(\rank)-1)})=\iota(f)$, and similarly for $g$. Therefore we also may and will suppose that $\bar{f}(\rank) = \bar{g}(\rank)$}.

With these extra assumptions, it remains to show that $f=g$. Since $\bar{f}(\rank)=\bar{g}(\rank)$ and $\bk_{f}=\bk_{g}$, we easily deduce that $\bar{f}(i) = \bar{g}(i)$ for all $1\leq i \leq \rank$. Let $0\leq r'_{f},r'_{g} < \rank$ such that $r'_{f} \ochangedd{\equiv} r_{f} \mod \rank$ and $r'_{g}\ochangedd{\equiv} r_{g} \mod \rank$. Since $\bm_{f}=\bm_{g}$ we have:
\begin{align*}
0 
= \sum_{i=1}^{\rank-1}i(\bm_{f,i}-\bm_{g,i}) 
&= \sum_{i=1}^{\rank}(g(i)-f(i))+\rank(f(N-r'_{f})-g(N-r'_{g}))+\levelthing(r'_{f}-r'_{g})\\
&= \levelthing(r'_{f}-r_{f}+r_{g}-r'_{g}) + \rank(f(N-r'_{f})-g(N-r'_{g})).
\end{align*}
We hence deduce that $\levelthing$ divides $f(N-r'_{f})-g(N-r'_{g})$, which is between $-\levelthing+1$ and $\levelthing-1$ so that $f(N-r'_{f})=g(N-r'_{g})$. The equality $\bm_{f}=\bm_{g}$ then implies that $f(i-r'_{f})=g(i-r'_{g})$ for all $i$ and that $r_{f}=r_{g}$ since $\levelthing(r_{f}-r_{g}) = \sum_{i=1}^{\rank}(f(i)-g(i))$. We obtain that $f=g$ as expected.

Therefore \changedd{$\bar{\iota}$ is an injective map} between the $\Z/\levelthing\Z$ orbits of $\Psi(Y,\pi)$ and the $\Z/\rank\Z$ orbits of $\{(\bm,\bk)\in (X^{+}(e))^{2}\mid \chi_{c}(\Ll_{\bm})=\chi_{c}(\Ll_{\bk})\}$. Since the number of such orbits is the same we deduce that this injection is a bijection.
\end{proof}

\subsubsection{\changedd{Proof of \autoref{T:Comparison}
\textit{Part (b).}}} \label{Pf:Comparison-b} 

Recall that the symmetric group $\mathfrak{S}_{\rank}$ acts on the weights. We will denote by $\bullet$ the action shifted by $\rho$. \changedd{We need the following two technical lemmas.}

\begin{Lemma}\label{L:ActionWeyl}
Let $f\in \Psi(Y,\pi)$ and $w\in \mathfrak{S}_{\rank}$. Then
\begin{gather*}
\ochangedd{(w\bullet \bm_{f})_{i} \equiv  f\big(w^{-1}(i+1)-r_{f}\big)-f\big(w^{-1}(i)-r_{f}\big)-1 \mod \levelthing,}
\\
\ochangedd{(w\bullet \bk_{f})_{i} \equiv \bar{f}\big(w^{-1}(i+1)\big)-\bar{f}\big(w^{-1}(i)\big)-1 \mod \levelthing.}
\end{gather*}
\end{Lemma}

\begin{proof}
It suffices to treat the case of a simple reflection $s_{j}$ through the hyperplane orthogonal to $\alpha_{j}$. \changedd{We only do the case of $\bm_{f}$, the case of $\bk_{f}$ being similar. By definition, we have $s_{j}\bullet \bm_{f} = \bm_{f}-\langle \bm_{f}+\rho,\alpha^{\vee}_{j}\rangle \alpha_{j}$. Using \eqref{eq:m_f-expanded}, we have $\langle \bm_{f}+\rho,\alpha^{\vee}_{j}\rangle \equiv f(j+1-r_{f})-f(j-r_{f}) \mod \rank$. Since $\alpha_{j} = -\omega_{j-1}+2\omega_{j}-\omega_{j+1}$, we verify that the statement of the lemma holds for $w=s_{j}$.}
\end{proof}

\begin{Lemma}\label{L:DifferenceRoots}
Let $f\in \Psi(Y,\pi)$ and $w,w'\in \mathfrak{S}_{\rank}$. Then $w\bullet \bm_{f} - w'\bullet \bk_{f} = \sum_{i=1}^{\rank-1}\nu_{i}\alpha_{i}$ with 
\begin{gather*}
\nu_i \ochangedd{\equiv}  \sum_{j=1}^{i}\Big(\bar{f}\big(w'^{-1}(j)\big)-f\big(w^{-1}(f)-r_{f}\big)\Big)\mod \levelthing.
\end{gather*}
\end{Lemma}

\begin{proof}
By definition of the fundamental roots, we have $\nu_{i} = \langle w\bullet\bm_{f}-w'\bullet\bk_{f},\omega_{i}\rangle$. \changedd{Using \autoref{eq:m_f-expanded}}, we have
\begin{align*}
\nu_{i} 
&= \levelthing\langle\omega_{r_{f}},\omega_{i}\rangle + \sum_{j=1}^{\rank-1}\Big(f\big(w^{-1}(j+1)-r_{f}\big) - f\big(w^{-1}(j)-r_{f}\big)-\bar{f}\big(w'^{-1}(j+1)\big)+\bar{f}\big(w'^{-1}(j)\big)\Big)\langle\omega_{j},\omega_{i}\rangle\\
&= \levelthing\langle\omega_{r_{f}},\omega_{i}\rangle + \sum_{j=1}^{\rank}\Big(f\big(w^{-1}(j)-r_{f}\big)-\bar{f}\big(w'^{-1}(j)\big)\Big)\langle\omega_{j-1} - \omega_{j}, \omega_{i}\rangle.
\end{align*}

Since $\langle \omega_{i}, \omega_{j} \rangle = \min(i,j)-\frac{ij}{\rank}$, we deduce that
\begin{align*}
\nu_{i} 
&= \levelthing\langle\omega_{r_{f}},\omega_{i}\rangle + \frac{i}{\rank}\sum_{j=1}^{\rank}\Big(f\big(w^{-1}(j)-r_{f}\big)-\bar{f}\big(w'^{-1}(j)\big)\Big) - \sum_{j=1}^{j}\Big(f\big(w^{-1}(j)-r_{f}\big)-\bar{f}\big(w'^{-1}(j)\big)\Big)\\
&=\levelthing\min(r_{f},i) + \sum_{j=1}^{j}\Big(\bar{f}\big(w'^{-1}(j)\big)-f\big(w^{-1}(j)-r_{f}\big)\Big),
\end{align*}
the last equality following from the definition of $r_{f}$.
\end{proof}

We \changedd{now} give another expression for the pre-Fourier matrix of the family $\family$, which will be more suitable for the comparison with the modular data arising from the asymptotic category. \changedd{Recall the matrix $\mathcal{S}=(\rooty^{-2ij})_{0\leq i,j < \levelthing}$.}

\begin{Lemma}\label{L:Adjugate}
Let $0\leq i_1 < \ldots < i_{\level} \leq \levelthing-1$ and $0\leq j_1 < \ldots < j_{\level} \leq \levelthing-1$. Let $\mathbf{i}$ be the $\level$-tuple $(i_1,\ldots,i_{\level})$, $\mathbf{j}$ be the $\level$-tuple $(j_1,\ldots,j_{\level})$ and ${}^c\mathbf{i}$ be the strictly increasing $n$-tuple obtained from the complement of $\mathbf{i}$ in $\{0,\ldots,\levelthing-1\}$ and similarly for ${}^c\mathbf{j}$. Then
\[    
\left(\Lambda^{\level}\mathcal{S}\right)_{\mathbf{i},\mathbf{j}}=(-1)^{\sum_{k=1}^{\level}({}^ci_k+{}^cj_k)}\frac{\det(\mathcal{S})}{\levelthing^{\rank}}\left(\Lambda^{\rank}\overline{\mathcal{S}}\right)_{{}^c\mathbf{i},{}^c\mathbf{j}}.
\]
\end{Lemma} 

\begin{proof}
Since $\frac{1}{\levelthing}\bar{\mathcal{S}}$ is the inverse of $\mathcal{S}$, the inverse of $\Lambda^{\rank}\mathcal{S}$ is $\frac{1}{\levelthing}\Lambda^{\rank}\bar{\mathcal{S}}$. But the inverse of $\Lambda^{\rank}\mathcal{S}$ can also be expressed in \changedd{terms} of the $\rank$th adjugate matrix of $\mathcal{S}$. The lemma follows then from the explicit form of the adjugate matrix.
\end{proof}

As $\tau = (-1)^{\binom{\levelthing}{2}}\det(\mathcal{S})$, we then deduce that, given $f,g\in \Psi(Y,\pi)$, the pre-Fourier matrix of $\family$ have the following expression:
\begin{gather*}
\tilde{\mathbb{S}}_{[f],[g]} = \frac{(-1)^{\sum_{i=1}^{\rank}(\bar{f}(i)+\bar{g}(i))}\varepsilon(f)\varepsilon(g)}{\levelthing^{\rank-1}}\sum_{w,w'\in\mathfrak{S}_{\rank}}(-1)^{l(w)+l(w')}\prod_{i=1}^{\rank}\rooty^{2(f(w(i))g(i)-\bar{f}(w'(i))\bar{g}(i))}.
\end{gather*}

Given $f,g\in \Psi(Y,\pi)$, using the explicit formula for $S$-matrix of $\dcenter{\ochangedd{\slncatr}}_{0}$ \changedd{(see \autoref{L:Smat-sln} for the $S$-matrix of $\slncatr$)}, we have
\begin{gather*}
S^{\ochangedd{0}}_{\iota(f),\iota(g)} \changedd{= S^{\sln[]}_{\bm_{f},\bm_{g}}\overline{S^{\sln[]}_{\bk_{f},\bk_{g}}}} = \frac{\sum_{w,w'\in\mathfrak{S}_{\rank}}(-1)^{l(w)+l(w')}\rooty^{2\langle\bm_{f}+\rho,w(\bm_{g}+\rho)\rangle-2\langle\bk_{f}+\rho,w'(\bk_{g}+\rho)\rangle}}{\sum_{w,w'\in\mathfrak{S}_{\rank}}(-1)^{l(w)+l(w')}\rooty^{2\langle\rho,w(\rho)\rangle-2\langle\rho,w'(\rho)\rangle}}.
\end{gather*}
We are then concerned with the value, modulo $\levelthing$, of $\langle\bm_{f}+\rho,w(\bm_{g}+\rho)\rangle-\langle\bk_{f}+\rho,w'(\bk_{g}+\rho)\rangle$. Using the shifted action of the symmetric group, we have
\begin{gather*}
\langle\bm_{f}+\rho,w(\bm_{g}+\rho)\rangle-\langle\bk_{f}+\rho,w'(\bk_{g}+\rho)\rangle = \langle\bm_{f}-\bk_{f},w'\bullet \bk_{g}+\rho\rangle + \langle w\bullet \bm_{g}-w'\bullet\bk_{g},w'\bullet\bk_{g}+\rho\rangle \changedd{.}
\end{gather*}
Using \autoref{L:ActionWeyl} and \autoref{L:DifferenceRoots}, we find
\begin{gather*}
\langle\bm_{f}-\bk_{f},w'\bullet \bk_{g}+\rho\rangle
\ochangedd{\equiv} 
\sum_{j=1}^{\rank-1}\Big(\bar{g}\big(w'^{-1}(\rank)\big)-\bar{g}\big(w'^{-1}(j)\big)\Big)\big(\bar{f}(j)-f(j-r_{f})\big) 
\mod \levelthing
\end{gather*}
and
\begin{gather*}
\langle w\bullet \bm_{g}-w'\bullet\bk_{g},w'\bullet\bk_{g}+\rho\rangle
\ochangedd{\equiv} 
\sum_{j=1}^{\rank-1}\big(f(n-r_{f})-f(j-r_{f})\big)\Big(\bar{g}\big(w'^{-1}(j)\big)-g\big(w^{-1}(j)-k_{g}\big)\ochangedd{\Big)}
\mod \levelthing.
\end{gather*}
Since $f\in \Psi(Y,\pi)$ we have $\sum_{i=1}^{\rank}(f(i)-\bar{f}(i))\ochangedd{\equiv} 0\mod \levelthing$ and similarly for $g$. Therefore, we obtain
\begin{gather*}
\langle\bm_{f}+\rho,w(\bm_{g}+\rho)\rangle-\langle\bk_{f}+\rho,w'(\bk_{g}+\rho)\rangle
\ochangedd{\equiv} 
\sum_{j=1}^{\rank}\Big(f(j-r_{f})-g\big(w^{-1}(j)-r_{g}\big)-\bar{f}(j)\bar{g}\big(w'^{-1}(j)\big)\Big)
\mod \levelthing
\end{gather*}
and since $w\mapsto (j\mapsto w(j+r_{g})-r_{g})$ is a signature preserving bijection of $\mathfrak{S}_{\rank}$ we have
\begin{gather*}
S_{\iota(f),\iota(g)} = \frac{\sum_{w,w'\in\mathfrak{S}_{\rank}}(-1)^{l(w)+l(w')}\prod_{i=1}^{\rank}\rooty^{f(w(j))g(j)-\bar{f}(w'(j))\bar{g}(j)}}{\sum_{w,w'\in\mathfrak{S}_{\rank}}(-1)^{l(w)+l(w')}\rooty^{2\langle\rho,w(\rho)\rangle-2\langle\rho,w'(\rho)\rangle}}.
\end{gather*}

\changedd{We finally} renormalize the $S$-matrix by the positive square root of the categorical dimension of $\dcenter{\verg{\levelthing}{\rank}}$\changedd{, in order to compare with the pre-Fourier matrix of family $\family$}. By \cite[Corollary 8.23.12]{EtGeNiOs-tensor-categories}, we have $\dim(\dcenter{\slncatr})=\rank^{2}\dim(\dcenter{\verg{\levelthing}{\rank}})$. Since \ochangedd{the category} $\slncatr$ is modular, we deduce that $\dim(\dcenter{\slncatr}) = \dim(\slncatr)^{2}$. Using \cite[Theorem 3.3.20]{BaKi-lecture-tensor}, we find that the positive square root of $\changedd{\dim(}\dcenter{\verg{\levelthing}{\rank}}\changedd{)}$ is then $\frac{\dim(\slncatr)}{\rank}=\levelthing^{\rank-1}\lvert \sum_{w\in \mathfrak{S}_{\rank}}(-1)^{l(w)}\rooty^{2\langle\rho,w(\rho)\rangle} \rvert^{-2}$.
\null\hfill$\square$

\subsubsection{\changedd{Proof of \autoref{T:Comparison}
\textit{Part (c).}}} \label{Pf:Comparison-c} We first rewrite slightly the eigenvalue of the Frobenius.

\begin{Lemma}\label{L:FrobeniusSimplified}
Given $f\in \Psi(Y,\pi)$ and $1\leq \changedd{s} \leq s(f)$, we have $\frob([f],\changedd{s}) = \rooty^{-\alpha(f)}$ with
\begin{gather*}
\alpha(f) \equiv -2\sum_{i=1}^{\rank}\big(\bar{f}(i)-f(i-r_{f})\big)\left(\sum_{j=1}^{i-1}\bar{f}(j)-\sum_{j=1}^{i}f(j-r_{f})\right) \mod 2\levelthing.
\end{gather*}
\end{Lemma}

\begin{proof}
The value of $\alpha(f)$ is given in \autoref{Eq:Frobenius}. First, using the definition of $\bar{f}$, we have
\begin{gather*}
\sum_{y\in Y}\big(f(y)^{2}+\levelthing f(y)\big) = \sum_{i=1}^{\rank}\big(f(i)^{2}+\levelthing f(i) - \bar{f}(i)^{2}-\levelthing \bar{f}(i)\big) + \sum_{j=1}^{\levelthing}(j^{2}+\levelthing j)
\end{gather*}
so that $\alpha(f) \equiv \sum_{i=1}^{\rank}(\bar{f}(i)^{2}+\levelthing \bar{f}(i) - f(i)^{2}-\levelthing f(i)) \mod 2\levelthing$. Now, we rewrite $\bar{f}(\rank)$ using the fact that $f\in\Psi(Y,\pi)$ and we obtain
\begin{multline*}
\bar{f}(\rank)^{2}+\levelthing \bar{f}(\rank)
= 
\left(f(\rank-r_{f})+\sum_{i=1}^{\rank-1}\big(f(i-r_{f})-f(i)\big)-er_{f}\right)^{2}\\
+\levelthing\left(f(\rank-r_{f})+\sum_{i=1}^{\rank-1}\big(f(i-r_{f})-f(i)\big)-er_{f}\right).
\end{multline*}
Expanding the square, we obtain
\begin{multline*}
\bar{f}(\rank)^{2}+\levelthing \bar{f}(\rank)
\equiv 
f(\rank-r_{f})^{2} + \left(\sum_{i=1}^{\rank-1}\big(f(i-r_{f})-f(i)\big)\right)^{2}+2f(\rank-r_{f})\sum_{i=1}^{\rank-1}\big(f(i-r_{f})-f(i)\big)\\
+ \levelthing\left(f(\rank-r_{f})+\sum_{i=1}^{\rank-1}\big(f(i-r_{f})-f(i)\big)\right) \mod 2\levelthing.
\end{multline*}
We then check that this implies that
\begin{gather*}
\alpha(f) \equiv 2\sum_{i=1}^{\rank-1}\big(f(i-r_{f})-\bar{f}(i)\big)\left(\sum_{j=i+1}^{\rank}f(j-r_{f})-\sum_{j=i}^{\rank-1}\bar{f}(j)\right) \mod 2\levelthing,
\end{gather*}
and we conclude the proof using once again that $\sum_{i=1}^{\rank}(\bar{f}(i)-f(i))\equiv 0 \mod \levelthing$.
\end{proof}

Given $g\in\Psi(Y,\pi)$, the value for the ribbon in $\dcenter{\verg{\levelthing}{\rank}}$ is $\theta_{\iota(f)} = \rooty^{\langle\bm_{f},\bm_{f}+2\rho\rangle-\langle\bk_{f},\bk_{f}+2\rho\rangle}$. We notice the equality $\langle\bm_{f},\bm_{f}+2\rho\rangle-\langle\bk_{f},\bk_{f}+2\rho\rangle = \langle\bm_{f}-\bk_{f},\bm_{f}-\bk_{f}\rangle + 2\langle\bm_{f}-\bk_{f},\bk_{f}+\rho\rangle$. On the one hand, using \autoref{L:DifferenceRoots}, $\langle\alpha_{i},\alpha_{j}\rangle=a_{i,j}$ and $\sum_{i=1}^{\rank}(f(i)-\bar{f}(i))\equiv 0 \mod\levelthing$ we have
\begin{align*}
\langle\bm_{f}-\bk_{f},\bm_{f}-\bk_{f}\rangle 
\equiv 2\sum_{i=1}^{\rank}\left(\sum_{j=1}^{i}\big(\bar{f}(j)-f(j-r_{f})\big)\right)\big(\bar{f}(i)-f(i-r_{f})\big)\mod 2\levelthing.
\end{align*}
On the other hand, using once again \autoref{L:DifferenceRoots}, we have
\begin{gather*}
\langle\bm_{f}-\bk_{f},\bk_{f}+\rho\rangle 
\equiv \sum_{i=1}^{\rank-1}\left(\sum_{j=1}^{i}\big(\bar{f}(j)-f(j-r_{f})\big)\right)\big(\bar{f}(i+1)-\bar{f}(i)\big)
\equiv -\sum_{i=1}^{\rank}\bar{f}(i)\big(\bar{f}(i)-f(i-r_{f})\big)\mod \levelthing
\end{gather*}
Therefore $\theta_{\iota(f)} = \rooty^{2\beta(f)}$, where
\begin{gather*}
\beta(f) = \sum_{i=1}^{\rank}\big(\bar{f}(i)-f(i-r_{f})\big)\left(\sum_{j=1}^{i-1}\bar{f}(j)-\sum_{j=1}^{i}f(j-r_{f})\right).
\end{gather*}
Then \autoref{L:FrobeniusSimplified} shows that the ribbon and the eigenvalue of the Frobenius coincide, and we conclude since the $T$-matrix is the diagonal matrix with entries the inverse of the ribbon.    
\null\hfill$\square$

\subsubsection{\changedd{Proof of \autoref{T:CMCells} \textit{Part (d)}}}\label{Pf:CMCells-d}
Since $\verG{\levelthing}{\rank}$ is a matrix category over $\verg{\levelthing}{\rank}$, we first compute the Artin--Wedderburn decomposition of $\ggroupc{\verg{\levelthing}{\rank}}$, which is a subring of $\mathrm{Mat}_{\rank}(\ggroupc{\slncatr})$, by determining all the primitive central idempotents.

It is easy to see that the center of $\ggroupc{\verg{\levelthing}{\rank}}$ consists of diagonal matrices with entries in the center of $\ggroupc{\slncatr_{0}}$. A central idempotent of $\ggroupc{\verg{\levelthing}{\rank}}$ is then a diagonal matrix with coefficients central idempotents of $\ggroupc{\slncatr_{0}}$.

\begin{Lemma}
The primitive central idempotents of $\ggroupc{\slncatr_{0}}$ are in bijection with the $\Ll_{\level\omega_{1}}\otimes\placeholder$-orbits on $\simples{\slncatr}$.
\end{Lemma}

\begin{proof}
Since $\slncatr_{0}$ is a braided fusion category, its complexified Grothendieck ring is commutative and semisimple. The primitive central idempotents of $\ggroupc{\slncatr_{0}}$ are then in bijection with the characters $\ggroupc{\slncatr_{0}}\to\C$. Given $\Ll_{\bm}\in\simples{\slncatr}$, the linear extension $\chi_{\bm}$ of $[\Ll_{\bk}]\mapsto S^{\sln[]}_{\bm,\bk}/\dim(\Ll_{\bm})$ is a character of $\ggroupc{\slncatr_{0}}$, see \cite[Proposition 8.13.11]{EtGeNiOs-tensor-categories}. Since $\slncatr$ is a modular category, it follows from \cite[Theorem 8.20.7 and Corollary 8.20.11]{EtGeNiOs-tensor-categories} that all this \changedd{exhausts} all characters $\ggroupc{\slncatr_{0}}\to\C$, and that $\chi_{\bm}=\chi_{\bm^{\prime}}$ if and only if $\Ll_{\bm}$ and $\Ll_{\bm^{\prime}}$ are in the same $(\slncatr_{0})^{\prime}$-module component. Here, $(\slncatr_{0})^{\prime}$ denotes the centralizer of $\slncatr_{0}$ in $\slncatr$ as in \cite[Section 8.20]{EtGeNiOs-tensor-categories}. 
\changeddd{Using the tensor product rules in $\slncatr$ (which can be obtained, for example, from \autoref{T:KoornwinderVariety}.(a)), one}
may easily show that $(\slncatr_{0})^{\prime}$ is the fusion subcategory generated by $\Ll_{\level\omega_{1}}$, which implies that $\Ll_{\bm}$ and $\Ll_{\bm^{\prime}}$ are in the same $(\slncatr_{0})^{\prime}$-module component if and only if $\Ll_{\bm}$ and $\Ll_{\bm^{\prime}}$ are in the same orbit under $\Ll_{\level\omega_{1}}\otimes\placeholder$.
\end{proof}

The idempotent of $\ggroupc{\slncat}$ corresponding to the object $\Ll_{\bm}$ is then a scalar multiple of $\regular_{\bm}=\sum_{X\in\simples{\slncatr}}\chi_{\bm}(X)[X^{*}]$. Let us denote by $\regular_{\bm,i}=\sum_{X\in\simples{\slncatr_{i}}}\chi_{\bm}(X)[X^{*}]$\changedd{. Then} the idempotent $\idemp_{\bm}$ of $\ggroupc{\slncat_{0}}$ corresponding to the object $\Ll_{\bm}$ is a scalar multiple of $\regular_{\bm,0}$. One can also check that $\regular_{\bm,0} [\Ll_{\bk}] = \chi_{\bm}([\Ll_{\bk}])\regular_{\bm,-i}$ for $\Ll_{\bk}$ a simple object of color $i$.

\begin{Lemma}
Let $\bm\in X^{+}(\level)$. Define $A_{\bm}=\{i\in\Z/\rank\Z | \exists \Ll_{\bk}, \chi_{\bm}(\Ll_{\bk})\neq 0 \text{ and } \chi_{c}(\Ll_{\bk}) = i\}$. Then $A_{\bm}$ is a subgroup of $\Z/\rank\Z$ and the diagonal matrix $\Idemp_{\bm,i+A_{\bm}}$ supported by the right coset $i+A_{\bm}$ with nonzero entries equal to $\idemp_{\bm}$ is a primitive central idempotent of $\ggroupc{\verg{\levelthing}{\rank}}$.
\end{Lemma}

\begin{proof}
First, it is clear that $A_{\bm}$ is an additive subgroup of $\Z/\rank\Z$. One may easily show that each diagonal matrix supported by a right coset of $A_{\bm}$ with nonzero entries equal to $\idemp_{\bm}$ is a central idempotent of $\ggroupc{\verg{\levelthing}{\rank}}$.

Let $f$ be a primitive central idempotent of $\ggroupc{\verg{\levelthing}{\rank}}$ and choose $i$ such that the diagonal entry $f_{i,i}$ is nonzero. Since $e$ is primitive, there exists $\bm$ such that $f_{i,i} = \idemp_{\bm}$. Suppose that there exists a simple object $\Ll_{\bk}$ of color $j$ such that $\chi_{\bm}(\Ll_{\bk})\neq 0$ and consider $x\in\ggroupc{\verg{\levelthing}{\rank}}$ whose only non zero entry is at the position $(i,i+j)$ and is equal to $[\Ll_{\bk}]$. Then, since $f$ is central, we obtain that $x_{i,i+j}f_{i+j,i+j} = f_{i,i}x_{i,i+j} = \idemp_{\bm}[L_{\bk}]$ which is a scalar multiple of $\regular_{\bm,-j}$. Since $\chi_{\bm}(\Ll_{\bk})\neq 0$, the element $\regular_{\bm,j}$ is nonzero as well as $x_{i,i+j}f_{i+j,i+j}$. In particular, the coefficient $f_{i+j,i+j}$ of $f$ is nonzero and a scalar multiple of $\regular_{\bm,0}$. Therefore $f_{i+j,i+j}=\idemp_{\bm}$. We deduce that $f=\Idemp_{m,i+A_{\bm}}$. 
\end{proof}

In the Artin--Wedderburn decomposition of $\ggroupc{\verg{\levelthing}{\rank}}$, the idempotent $\Idemp_{\bm,i+A_{\bm}}$ will contribute to a matrix algebra of size $\lvert A_{\bm} \rvert$. Let $\Stabilizer{\bm}$ be the stabilizer subgroup of $\Ll_{\bm}$ under the action \changedd{of} $\Ll_{\level\omega_{1}}\otimes\placeholder$. 

\begin{Lemma}
Let $\bm\in X^{+}(\level)$ and $r\in\Z/\rank\Z$. Then $r\in\Stabilizer{\bm}$ if and only if $ri \equiv 0 \bmod \rank$ for all $i\in A_{\bm}$.
\end{Lemma}

\begin{proof}
Using \autoref{L:GradingBraiding}, we have, for any $\bm,\bk\in X^{+}(\level)$ and $r\in\N$:
\begin{gather*}
\chi_{\bm^{\turnsymbol r}}(\Ll_{\bk}) 
= \frac{\dim(\Ll_{\bk})}{\dim(\Ll_{\bm})}\chi_{\bk}(\Ll_{\bm^{\turnsymbol r}}) 
= \frac{\dim(\Ll_{\bk})}{\dim(\Ll_{\bm})}\chi_{\bk}(\Ll_{\bm})\chi_{\bk}(\Ll_{\level\omega_{r}})
= \rooty^{-2r\level\chi_{c}(\Ll_{\bk})/\rank}\chi_{\bm}(\Ll_{\bk}).
\end{gather*}

Suppose that $r\in\Stabilizer{\bm}$ and let $i\in A_{\bm}$. We then choose $\bk\in X^{+}(\level)$ such that $\chi_{c}(\Ll_{\bk})=i$ and $\chi_{\bm}(\Ll_{\bk})\neq 0$. By assumption, $\bm^{\turnsymbol r}=\bm$ and therefore we have $\chi_{\bm}(\Ll_{\bk}) = \chi_{\bm^{\turnsymbol r}}(\Ll_{\bk}) = \rooty^{-2r\level\chi_{c}(\Ll_{\bk})/\rank}\chi_{\bm}(\Ll_{\bk})$. As $\chi_{\bm}(\Ll_{\bk})\neq 0$, we obtain that $r\chi_{\bm}(\Ll_{\bk})\equiv 0 \mod \rank$.

Conversely, suppose that $ri\equiv 0 \bmod \rank$ for all $i\in A_{\bm}$. Since $\slncatr$ is modular, $\chi_{\bm}=\chi_{\bm^{\turnsymbol r}}$ implies that $\Ll_{\bm}\simeq \Ll_{\bm^{\turnsymbol r}}$. Let $\bk\in X^{+}(\level)$. If $\chi_{c}(\Ll_{\bk})\not\in A_{\bm}$ then $\chi_{\bm}(\Ll_{\bk})=0$ and $\chi_{\bm^{\turnsymbol r}}(\Ll_{\bk})=\rooty^{-2r\level\chi_{c}(\Ll_{\bk})/\rank}\chi_{\bm}(\Ll_{\bk})=0$. Otherwise, $r\chi_{c}(\Ll_{\bk})\equiv 0\bmod \rank$ and we have $\chi_{\bm^{\turnsymbol r}}(\Ll_{\bk})=\rooty^{-2r\level\chi_{c}(\Ll_{\bk})/\rank}\chi_{\bm}(\Ll_{\bk})=\chi_{\bm}(\Ll_{\bk})$. Therefore the characters $\chi_{\bm}$ and $\chi_{\bm^{\turnsymbol r}}$ conicide.
\end{proof}

We hence deduce that $\lvert A_{\bm}\rvert$ is equal to the cardinal of the orbit of $\Ll_{\bm}$ under $\Ll_{\level\omega_{1}}\otimes\placeholder$.

\begin{Lemma}\label{L:Wedderburn}
We have
\begin{gather*}
\ggroupc{\verg{\levelthing}{\rank}} \simeq \bigoplus_{m\mid \gcd(\levelthing,\rank)} \mathrm{Mat}_{\rank/m}(\C)^{\oplus n_{m}}, 
\quad\ggroupc{\verG{\levelthing}{\rank}} \simeq \bigoplus_{m\mid \gcd(\levelthing,\rank)} \mathrm{Mat}_{\rank!/m}(\C)^{\oplus n_{m}},
\end{gather*}
where $n_{m} = \frac{m^{2}}{\levelthing}\sum_{k\mid\gcd(\rank/m,\levelthing/m)}\mu(k)\binom{\levelthing/mk}{\rank/mk}$.
\end{Lemma}

\begin{proof}
To obtain the Artin--Wedderburn decomposition of $\verg{\levelthing}{\rank}$, it remains to compute the number of orbits of simple objects of $\slncatr$ under $\Ll_{\level\omega_{1}}\otimes\placeholder$ of a given cardinality $m$, which is done similarly to \autoref{L:ObjectsStabilizer}. Indeed, each such orbit gives rise to $m$ idempotents with a matrix algebra of size $\rank/m$. The Artin--Wedderburn decomposition of $\verG{\levelthing}{\rank}$ follows from the definition of the big asymptotic category.
\end{proof}

\changedd{The next lemma studies} the dimensions of the representations of the Calogero--Moser family $\family_{0}$.

\begin{Lemma}\label{L:principal_series}
The Calogero--Moser family $\family_{0}$ contains only representations of dimension $\rank!/m$ for every $m|\gcd(\levelthing,\rank)$. The number of representations of dimension $\rank!/m$ is $n_{m}$.
\end{Lemma}

\begin{proof}
Recall that in term\ochangedd{s} of $\levelthing$-partitions of $\rank$, the representations of $\family_{0}$ are indexed by orbits of multipartitions with entries equal to $(1)$ or $\emptyset$ under the action of $\Z/\levelthing\Z$ by cyclic shifting. Each orbit of size $m$ parametrizes $\levelthing/m$ representations of dimension $\rank!/m$.

To count the number of representations of a given dimension, one may proceed similarly to the proof of \autoref{T:Comparison}.(a).
\end{proof}

\changedd{Therefore, combining \autoref{L:Wedderburn} and \autoref{L:principal_series}, we obtain
\begin{equation*}
\ggroupc{\verG{\levelthing}{\rank}} \simeq \bigoplus_{m\mid \gcd(\levelthing,\rank)} \mathrm{Mat}_{\rank!/m}(\C)^{\oplus n_{m}}\simeq \bigoplus_{V\in \family_{0}}\mathrm{Mat}_{\dim_{\C}(V)}(\C).
\end{equation*}
The proof is now complete.}
\null\hfill$\square$

\newcommand{\etalchar}[1]{$^{#1}$}

\end{document}

%% file: figures/tetrahedron-1.tex
\tdplotsetmaincoords{2}{0}
\begin{tikzpicture}[tdplot_main_coords,scale=1.4,anchorbase]
  \coordinate (o1) at (4/3,1/3,1/3);
  \coordinate (o2) at (1/3,4/3,1/3);
  \coordinate (o3) at (1/3,1/3,4/3);
  \draw (0,0,0) -- (o1);
  \draw (0,0,0) -- (o2);
  \draw (0,0,0) -- (o3);
  \draw (o1) -- (o2);
  \draw (o1) -- (o3);
  \draw (o2) -- (o3);
  \node[circle,inner sep=0.9*1.8pt,fill=zero] at (0,0,0){};
  \node[triangle,inner sep=0.9*1.6pt,fill=one] at (o1){};
  \node[square,inner sep=0.9*1.8pt,fill=two] at (o2){};
  \node[pentagon,inner sep=0.9*2.0pt,fill=three] at (o3){};
\end{tikzpicture}

%% file: figures/tetrahedron-4.tex
\tdplotsetmaincoords{2}{0}
\begin{tikzpicture}[tdplot_main_coords,square/.style={regular polygon,regular polygon sides=4},triangle/.style={regular polygon,regular polygon sides=3},pentagon/.style={regular polygon,regular polygon sides=5},scale=1.4,anchorbase]
  \coordinate (o1) at (4/3,1/3,1/3);
  \coordinate (o2) at (1/3,4/3,1/3);
  \coordinate (o3) at (1/3,1/3,4/3);
  \draw (0,0,0) -- (o1);
  \draw (0,0,0) -- (o2);
  \draw (0,0,0) -- (o3);
  \draw (o1) -- (o2);
  \draw (o1) -- (o3);
  \draw (o2) -- (o3);
  \draw (o1) -- ($2*(o1)$);
  \draw (o1) -- ($(o1)+(o2)$);
  \draw (o1) -- ($(o1)+(o3)$);
  \draw (o2) -- ($(o1)+(o2)$);
  \draw (o2) -- ($2*(o2)$);
  \draw (o2) -- ($(o2)+(o3)$);
  \draw (o2) -- ($(o1)+(o3)$);
  \draw (o3) -- ($(o1)+(o3)$);
  \draw (o3) -- ($(o2)+(o3)$);
  \draw (o3) -- ($2*(o3)$);
  \draw ($2*(o1)$) -- ($(o1)+(o2)$);
  \draw ($2*(o1)$) -- ($(o1)+(o3)$);
  \draw ($(o1)+(o2)$) -- ($(o1)+(o3)$);
  \draw ($(o1)+(o2)$) -- ($2*(o2)$);
  \draw ($(o1)+(o2)$) -- ($(o2)+(o3)$);
  \draw ($(o1)+(o3)$) -- ($(o2)+(o3)$);
  \draw ($(o1)+(o3)$) -- ($2*(o3)$);
  \draw ($2*(o2)$) -- ($(o2)+(o3)$);
  \draw ($(o2)+(o3)$) -- ($2*(o3)$);
  \draw ($2*(o1)$) -- ($3*(o1)$);
  \draw ($2*(o1)$) -- ($2*(o1)+(o2)$);
  \draw ($2*(o1)$) -- ($2*(o1)+(o3)$);
  \draw ($(o1)+(o2)$) -- ($2*(o1)+(o2)$);
  \draw ($(o1)+(o2)$) -- ($(o1)+2*(o2)$);
  \draw ($(o1)+(o2)$) -- ($(o1)+(o2)+(o3)$);
  \draw ($(o1)+(o2)$) -- ($2*(o1)+(o3)$);
  \draw ($(o1)+(o3)$) -- ($2*(o1)+(o3)$);
  \draw ($(o1)+(o3)$) -- ($(o1)+(o2)+(o3)$);
  \draw ($(o1)+(o3)$) -- ($(o1)+2*(o3)$);
  \draw ($2*(o2)$) -- ($(o1)+2*(o2)$);
  \draw ($2*(o2)$) -- ($3*(o2)$);
  \draw ($2*(o2)$) -- ($2*(o2)+(o3)$);
  \draw ($2*(o2)$) -- ($(o1)+(o2)+(o3)$);
  \draw ($(o2)+(o3)$) -- ($(o1)+(o2)+(o3)$);
  \draw ($(o2)+(o3)$) -- ($2*(o2)+(o3)$);
  \draw ($(o2)+(o3)$) -- ($(o2)+2*(o3)$);
  \draw ($(o2)+(o3)$) -- ($(o1)+2*(o3)$);
  \draw ($2*(o3)$) -- ($(o1)+2*(o3)$);
  \draw ($2*(o3)$) -- ($(o2)+2*(o3)$);
  \draw ($2*(o3)$) -- ($3*(o3)$);
  \draw ($3*(o1)$) -- ($2*(o1)+(o2)$);
  \draw ($3*(o1)$) -- ($2*(o1)+(o3)$);
  \draw ($2*(o1)+(o2)$) -- ($(o1)+2*(o2)$);
  \draw ($2*(o1)+(o2)$) -- ($(o1)+(o2)+(o3)$);
  \draw ($2*(o1)+(o2)$) -- ($2*(o1)+(o3)$);
  \draw ($(o1)+2*(o2)$) -- ($3*(o2)$);
  \draw ($(o1)+2*(o2)$) -- ($2*(o2)+(o3)$);
  \draw ($(o1)+2*(o2)$) -- ($(o1)+(o2)+(o3)$);
  \draw ($3*(o2)$) -- ($2*(o2)+(o3)$);
  \draw ($2*(o1)+(o3)$) -- ($(o1)+2*(o3)$);
  \draw ($2*(o1)+(o3)$) -- ($(o1)+(o2)+(o3)$);
  \draw ($(o1)+(o2)+(o3)$) -- ($(o1)+2*(o3)$);
  \draw ($(o1)+(o2)+(o3)$) -- ($(o2)+2*(o3)$);
  \draw ($(o1)+(o2)+(o3)$) -- ($2*(o2)+(o3)$);
  \draw ($2*(o2)+(o3)$) -- ($(o2)+2*(o3)$);
  \draw ($(o1)+2*(o3)$) -- ($3*(o3)$);
  \draw ($(o1)+2*(o3)$) -- ($(o2)+2*(o3)$);
  \draw ($(o2)+2*(o3)$) -- ($3*(o3)$);
  \draw ($3*(o1)$) -- ($4*(o1)$);
  \draw ($3*(o1)$) -- ($3*(o1)+(o2)$);
  \draw ($3*(o1)$) -- ($3*(o1)+(o3)$);
  \draw ($2*(o1)+(o2)$) -- ($3*(o1)+(o2)$);
  \draw ($2*(o1)+(o2)$) -- ($2*(o1)+2*(o2)$);
  \draw ($2*(o1)+(o2)$) -- ($2*(o1)+(o2)+(o3)$);
  \draw ($2*(o1)+(o2)$) -- ($3*(o1)+(o3)$);
  \draw ($2*(o1)+(o3)$) -- ($3*(o1)+(o3)$);
  \draw ($2*(o1)+(o3)$) -- ($2*(o1)+(o2)+(o3)$);
  \draw ($2*(o1)+(o3)$) -- ($2*(o1)+2*(o3)$);
  \draw ($(o1)+2*(o2)$) -- ($2*(o1)+2*(o2)$);
  \draw ($(o1)+2*(o2)$) -- ($(o1)+3*(o2)$);
  \draw ($(o1)+2*(o2)$) -- ($(o1)+2*(o2)+(o3)$);
  \draw ($(o1)+2*(o2)$) -- ($2*(o1)+(o2)+(o3)$);
  \draw ($(o1)+(o2)+(o3)$) -- ($2*(o1)+(o2)+(o3)$);
  \draw ($(o1)+(o2)+(o3)$) -- ($(o1)+2*(o2)+(o3)$);
  \draw ($(o1)+(o2)+(o3)$) -- ($(o1)+(o2)+2*(o3)$);
  \draw ($(o1)+(o2)+(o3)$) -- ($2*(o1)+2*(o3)$);
  \draw ($(o1)+2*(o3)$) -- ($2*(o1)+2*(o3)$);
  \draw ($(o1)+2*(o3)$) -- ($(o1)+(o2)+2*(o3)$);
  \draw ($(o1)+2*(o3)$) -- ($(o1)+3*(o3)$);
  \draw ($3*(o2)$) -- ($(o1)+3*(o2)$);
  \draw ($3*(o2)$) -- ($4*(o2)$);
  \draw ($3*(o2)$) -- ($3*(o2)+(o3)$);
  \draw ($3*(o2)$) -- ($(o1)+2*(o2)+(o3)$);
  \draw ($2*(o2)+(o3)$) -- ($(o1)+2*(o2)+(o3)$);
  \draw ($2*(o2)+(o3)$) -- ($3*(o2)+(o3)$);
  \draw ($2*(o2)+(o3)$) -- ($2*(o2)+2*(o3)$);
  \draw ($2*(o2)+(o3)$) -- ($(o1)+(o2)+2*(o3)$);
  \draw ($(o2)+2*(o3)$) -- ($(o1)+(o2)+2*(o3)$);
  \draw ($(o2)+2*(o3)$) -- ($2*(o2)+2*(o3)$);
  \draw ($(o2)+2*(o3)$) -- ($(o2)+3*(o3)$);
  \draw ($(o2)+2*(o3)$) -- ($(o1)+3*(o3)$);
  \draw ($3*(o3)$) -- ($(o1)+3*(o3)$);
  \draw ($3*(o3)$) -- ($(o2)+3*(o3)$);
  \draw ($3*(o3)$) -- ($4*(o3)$);
  \draw ($4*(o1)$) -- ($3*(o1)+(o2)$);
  \draw ($4*(o1)$) -- ($3*(o1)+(o3)$);
  \draw ($3*(o1)+(o2)$) -- ($3*(o1)+(o3)$);
  \draw ($3*(o1)+(o2)$) -- ($2*(o1)+(o2)+(o3)$);
  \draw ($3*(o1)+(o2)$) -- ($2*(o1)+2*(o2)$);
  \draw ($3*(o1)+(o3)$) -- ($2*(o1)+2*(o3)$);
  \draw ($3*(o1)+(o3)$) -- ($2*(o1)+(o2)+(o3)$);
  \draw ($2*(o1)+2*(o2)$) -- ($2*(o1)+(o2)+(o3)$);
  \draw ($2*(o1)+2*(o2)$) -- ($(o1)+2*(o2)+(o3)$);
  \draw ($2*(o1)+2*(o2)$) -- ($(o1)+3*(o2)$);
  \draw ($2*(o1)+(o2)+(o3)$) -- ($2*(o1)+2*(o3)$);
  \draw ($2*(o1)+(o2)+(o3)$) -- ($(o1)+(o2)+2*(o3)$);
  \draw ($2*(o1)+(o2)+(o3)$) -- ($(o1)+2*(o2)+(o3)$);
  \draw ($2*(o1)+2*(o3)$) -- ($(o1)+3*(o3)$);
  \draw ($2*(o1)+2*(o3)$) -- ($(o1)+(o2)+2*(o3)$);
  \draw ($(o1)+3*(o2)$) -- ($(o1)+2*(o2)+(o3)$);
  \draw ($(o1)+3*(o2)$) -- ($3*(o2)+(o3)$);
  \draw ($(o1)+3*(o2)$) -- ($4*(o2)$);
  \draw ($(o1)+2*(o2)+(o3)$) -- ($(o1)+(o2)+2*(o3)$);
  \draw ($(o1)+2*(o2)+(o3)$) -- ($2*(o2)+2*(o3)$);
  \draw ($(o1)+2*(o2)+(o3)$) -- ($3*(o2)+(o3)$);
  \draw ($(o1)+(o2)+2*(o3)$) -- ($(o1)+3*(o3)$);
  \draw ($(o1)+(o2)+2*(o3)$) -- ($(o2)+3*(o3)$);
  \draw ($(o1)+(o2)+2*(o3)$) -- ($2*(o2)+2*(o3)$);
  \draw ($(o1)+3*(o3)$) -- ($(o2)+3*(o3)$);
  \draw ($(o1)+3*(o3)$) -- ($4*(o3)$);
  \draw ($4*(o2)$) -- ($3*(o2)+(o3)$);
  \draw ($3*(o2)+(o3)$) -- ($2*(o2)+2*(o3)$);
  \draw ($2*(o2)+2*(o3)$) -- ($(o2)+3*(o3)$);
  \draw ($(o2)+3*(o3)$) -- ($4*(o3)$);
  \node[circle,inner sep=0.9*1.8pt,fill=zero] at (0,0,0){};
  \node[triangle,inner sep=0.9*1.6pt,fill=one] at (o1){};
  \node[square,inner sep=0.9*1.8pt,fill=two] at (o2){};
  \node[pentagon,inner sep=0.9*2.0pt,fill=three] at (o3){};
  \node[square,inner sep=0.9*1.8pt,fill=two] at ($2*(o1)$){};
  \node[pentagon,inner sep=0.9*2.0pt,fill=three] at ($(o1)+(o2)$){};
  \node[circle,inner sep=0.9*1.8pt,fill=zero] at ($(o1)+(o3)$){};
  \node[circle,inner sep=0.9*1.8pt,fill=zero] at ($2*(o2)$){};
  \node[triangle,inner sep=0.9*1.6pt,fill=one] at ($(o2)+(o3)$){};
  \node[square,inner sep=0.9*1.8pt,fill=two] at ($2*(o3)$){};
  \node[pentagon,inner sep=0.9*2.0pt,fill=three] at ($3*(o1)$){};
  \node[circle,inner sep=0.9*1.8pt,fill=zero] at ($2*(o1)+(o2)$){};
  \node[triangle,inner sep=0.9*1.6pt,fill=one] at ($2*(o1)+(o3)$){};
  \node[triangle,inner sep=0.9*1.6pt,fill=one] at ($(o1)+2*(o2)$){};
  \node[square,inner sep=0.9*1.8pt,fill=two] at ($(o1)+(o2)+(o3)$){};
  \node[pentagon,inner sep=0.9*2.0pt,fill=three] at ($(o1)+2*(o3)$){};
  \node[square,inner sep=0.9*1.8pt,fill=two] at ($3*(o2)$){};
  \node[pentagon,inner sep=0.9*2.0pt,fill=three] at ($2*(o2)+(o3)$){};
  \node[circle,inner sep=0.9*1.8pt,fill=zero] at ($(o2)+2*(o3)$){};
  \node[triangle,inner sep=0.9*1.6pt,fill=one] at ($3*(o3)$){};
  \node[circle,inner sep=0.9*1.8pt,fill=zero] at ($4*(o1)$){};
  \node[triangle,inner sep=0.9*1.6pt,fill=one] at ($3*(o1)+(o2)$){};
  \node[square,inner sep=0.9*1.8pt,fill=two] at ($3*(o1)+(o3)$){};
  \node[square,inner sep=0.9*1.8pt,fill=two] at ($2*(o1)+2*(o2)$){};
  \node[pentagon,inner sep=0.9*2.0pt,fill=three] at ($2*(o1)+(o2)+(o3)$){};
  \node[circle,inner sep=0.9*1.8pt,fill=zero] at ($2*(o1)+2*(o3)$){};
  \node[pentagon,inner sep=0.9*2.0pt,fill=three] at ($(o1)+3*(o2)$){};
  \node[circle,inner sep=0.9*1.8pt,fill=zero] at ($(o1)+2*(o2)+(o3)$){};
  \node[triangle,inner sep=0.9*1.6pt,fill=one] at ($(o1)+(o2)+2*(o3)$){};
  \node[square,inner sep=0.9*1.8pt,fill=two] at ($(o1)+3*(o3)$){};
  \node[circle,inner sep=0.9*1.8pt,fill=zero] at ($4*(o2)$){};
  \node[triangle,inner sep=0.9*1.6pt,fill=one] at ($3*(o2)+(o3)$){};
  \node[square,inner sep=0.9*1.8pt,fill=two] at ($2*(o2)+2*(o3)$){};
  \node[pentagon,inner sep=0.9*2.0pt,fill=three] at ($(o2)+3*(o3)$){};
  \node[circle,inner sep=0.9*1.8pt,fill=zero] at ($4*(o3)$){};
\end{tikzpicture}

%% file: figures/genericsubgraph.tex
\tdplotsetmaincoords{0}{0}
1\colon
\begin{tikzpicture}[tdplot_main_coords,anchorbase]
  \coordinate (o1) at (4/3,1/3,1/3);
  \coordinate (o2) at (1/3,4/3,1/3);
  \coordinate (o3) at (1/3,1/3,4/3);

  \draw[thick,->-] (0,0,0) -- (o1);
  \draw[thick,->-] (0,0,0) -- ($-1*(o1)+(o2)$);
  \draw[thick,->-] (0,0,0) -- ($-1*(o2)+(o3)$);
  \draw[thick,->-] (0,0,0) -- ($-1*(o3)$);

  \draw[gray] (0,0,0) -- (o2);
  \draw[gray] (0,0,0) -- ($(o1)-1*(o2)+(o3)$);
  \draw[gray] (0,0,0) -- ($-1*(o1)+(o3)$);
  \draw[gray] (0,0,0) -- ($(o1)-1*(o3)$);
  \draw[gray] (0,0,0) -- ($-1*(o1)+(o2)-1*(o3)$);
  \draw[gray] (0,0,0) -- ($-1*(o2)$);
  
  \draw[gray] (0,0,0) -- (o3);
  \draw[gray] (0,0,0) -- ($(o2)-1*(o3)$);
  \draw[gray] (0,0,0) -- ($(o1)-1*(o2)$);
  \draw[gray] (0,0,0) -- ($-1*(o1)$);
  
  \node[circle,inner sep=1.8pt,fill=zero] at (0,0,0){};

  \node[triangle,inner sep=1.6pt,fill=one] at (o1){};
  \node[triangle,inner sep=1.6pt,fill=one] at ($-1*(o1)+(o2)$){};
  \node[triangle,inner sep=1.6pt,fill=one] at ($-1*(o2)+(o3)$){};
  \node[triangle,inner sep=1.6pt,fill=one] at ($-1*(o3)$){};

  \node[square,inner sep=1.8pt,fill=two] at (o2){};
  \node[square,inner sep=1.8pt,fill=two] at ($(o1)-1*(o2)+(o3)$){};
  \node[square,inner sep=1.8pt,fill=two] at ($-1*(o1)+(o3)$){};
  \node[square,inner sep=1.8pt,fill=two] at ($(o1)-1*(o3)$){};
  \node[square,inner sep=1.8pt,fill=two] at ($-1*(o1)+(o2)-1*(o3)$){};
  \node[square,inner sep=1.8pt,fill=two] at ($-1*(o2)$){};

  \node[pentagon,inner sep=2.0pt,fill=three] at (o3){};
  \node[pentagon,inner sep=2.0pt,fill=three] at ($(o2)-1*(o3)$){};
  \node[pentagon,inner sep=2.0pt,fill=three] at ($(o1)-1*(o2)$){};
  \node[pentagon,inner sep=2.0pt,fill=three] at ($-1*(o1)$){};  
\end{tikzpicture}
,\quad 2\colon
\begin{tikzpicture}[tdplot_main_coords,anchorbase]
  \coordinate (o1) at (4/3,1/3,1/3);
  \coordinate (o2) at (1/3,4/3,1/3);
  \coordinate (o3) at (1/3,1/3,4/3);

  \draw[gray] (0,0,0) -- (o1);
  \draw[gray] (0,0,0) -- ($-1*(o1)+(o2)$);
  \draw[gray] (0,0,0) -- ($-1*(o2)+(o3)$);
  \draw[gray] (0,0,0) -- ($-1*(o3)$);

  \draw[thick,->-] (0,0,0) -- (o2);
  \draw[thick,->-] (0,0,0) -- ($(o1)-1*(o2)+(o3)$);
  \draw[thick,->-] (0,0,0) -- ($-1*(o1)+(o3)$);
  \draw[thick,->-] (0,0,0) -- ($(o1)-1*(o3)$);
  \draw[thick,->-] (0,0,0) -- ($-1*(o1)+(o2)-1*(o3)$);
  \draw[thick,->-] (0,0,0) -- ($-1*(o2)$);

  \draw[gray] (0,0,0) -- (o3);
  \draw[gray] (0,0,0) -- ($(o2)-1*(o3)$);
  \draw[gray] (0,0,0) -- ($(o1)-1*(o2)$);
  \draw[gray] (0,0,0) -- ($-1*(o1)$);
  
  \node[circle,inner sep=1.8pt,fill=zero] at (0,0,0){};
  \node[triangle,inner sep=1.6pt,fill=one] at (o1){};
  \node[triangle,inner sep=1.6pt,fill=one] at ($-1*(o1)+(o2)$){};
  \node[triangle,inner sep=1.6pt,fill=one] at ($-1*(o2)+(o3)$){};
  \node[triangle,inner sep=1.6pt,fill=one] at ($-1*(o3)$){};

  \node[square,inner sep=1.8pt,fill=two] at (o2){};
  \node[square,inner sep=1.8pt,fill=two] at ($(o1)-1*(o2)+(o3)$){};
  \node[square,inner sep=1.8pt,fill=two] at ($-1*(o1)+(o3)$){};
  \node[square,inner sep=1.8pt,fill=two] at ($(o1)-1*(o3)$){};
  \node[square,inner sep=1.8pt,fill=two] at ($-1*(o1)+(o2)-1*(o3)$){};
  \node[square,inner sep=1.8pt,fill=two] at ($-1*(o2)$){};

  \node[pentagon,inner sep=2.0pt,fill=three] at (o3){};
  \node[pentagon,inner sep=2.0pt,fill=three] at ($(o2)-1*(o3)$){};
  \node[pentagon,inner sep=2.0pt,fill=three] at ($(o1)-1*(o2)$){};
  \node[pentagon,inner sep=2.0pt,fill=three] at ($-1*(o1)$){};  
\end{tikzpicture}
,\quad 3\colon
\begin{tikzpicture}[tdplot_main_coords,anchorbase]
  \coordinate (o1) at (4/3,1/3,1/3);
  \coordinate (o2) at (1/3,4/3,1/3);
  \coordinate (o3) at (1/3,1/3,4/3);

  \draw[gray] (0,0,0) -- (o1);
  \draw[gray] (0,0,0) -- ($-1*(o1)+(o2)$);
  \draw[gray] (0,0,0) -- ($-1*(o2)+(o3)$);
  \draw[gray] (0,0,0) -- ($-1*(o3)$);

  \draw[gray] (0,0,0) -- (o2);
  \draw[gray] (0,0,0) -- ($(o1)-1*(o2)+(o3)$);
  \draw[gray] (0,0,0) -- ($-1*(o1)+(o3)$);
  \draw[gray] (0,0,0) -- ($(o1)-1*(o3)$);
  \draw[gray] (0,0,0) -- ($-1*(o1)+(o2)-1*(o3)$);
  \draw[gray] (0,0,0) -- ($-1*(o2)$);
  
  \draw[thick,->-] (0,0,0) -- (o3);
  \draw[thick,->-] (0,0,0) -- ($(o2)-1*(o3)$);
  \draw[thick,->-] (0,0,0) -- ($(o1)-1*(o2)$);
  \draw[thick,->-] (0,0,0) -- ($-1*(o1)$);
  
  \node[circle,inner sep=1.8pt,fill=zero] at (0,0,0){};
  
  \node[triangle,inner sep=1.6pt,fill=one] at (o1){};
  \node[triangle,inner sep=1.6pt,fill=one] at ($-1*(o1)+(o2)$){};
  \node[triangle,inner sep=1.6pt,fill=one] at ($-1*(o2)+(o3)$){};
  \node[triangle,inner sep=1.6pt,fill=one] at ($-1*(o3)$){};

  \node[square,inner sep=1.8pt,fill=two] at (o2){};
  \node[square,inner sep=1.8pt,fill=two] at ($(o1)-1*(o2)+(o3)$){};
  \node[square,inner sep=1.8pt,fill=two] at ($-1*(o1)+(o3)$){};
  \node[square,inner sep=1.8pt,fill=two] at ($(o1)-1*(o3)$){};
  \node[square,inner sep=1.8pt,fill=two] at ($-1*(o1)+(o2)-1*(o3)$){};
  \node[square,inner sep=1.8pt,fill=two] at ($-1*(o2)$){};

  \node[pentagon,inner sep=2.0pt,fill=three] at (o3){};
  \node[pentagon,inner sep=2.0pt,fill=three] at ($(o2)-1*(o3)$){};
  \node[pentagon,inner sep=2.0pt,fill=three] at ($(o1)-1*(o2)$){};
  \node[pentagon,inner sep=2.0pt,fill=three] at ($-1*(o1)$){};
\end{tikzpicture}
.

%% file: figures/koorn-A-4-4.tex
\begin{tikzpicture}[anchorbase]
\def\a{1} \def\b{4} \def\scale{1}
\draw[cyan!30,very thin] ({-\b*\scale},{-\b*\scale}) grid ({\b*\scale},{\b*\scale});
\draw[->] ({-\b*\scale},0) -- ({\b*\scale},0);
\draw[->] (0,{-\b*\scale}) -- (0,{\b*\scale});
\draw (0,0) circle ({\b*\scale});
\draw[line width=2pt,red] plot[samples=100,domain=0:\a*360,smooth,variable=\t] ({\scale*((\b-\a)*cos(\t)+\a*cos((\b-\a)*\t/\a)},{\scale*((\b-\a)*sin(\t)-\a*sin((\b-\a)*\t/\a)});
\draw[blue,fill=blue] (2.613125929752753,0.0) circle (4*0.200000000000000pt);
\draw[blue,fill=blue] (1.8122548927057758,0.36047991100347393) circle (4*0.200000000000000pt);
\draw[blue,fill=blue] (1.0,1.0) circle (4*0.200000000000000pt);
\draw[blue,fill=blue] (0.360479911003474,1.8122548927057758) circle (4*0.200000000000000pt);
\draw[blue,fill=blue] (-0.0,2.613125929752753) circle (4*0.200000000000000pt);
\draw[blue,fill=blue] (1.4142135623730951,0.0) circle (4*0.200000000000000pt);
\draw[blue,fill=blue] (0.6363792902864169,0.42521504738362825) circle (4*0.200000000000000pt);
\draw[blue,fill=blue] (0.0,1.082392200292394) circle (4*0.200000000000000pt);
\draw[blue,fill=blue] (-0.360479911003474,1.8122548927057758) circle (4*0.200000000000000pt);
\draw[blue,fill=blue] (0.0,0.0) circle (4*0.600000000000000pt);
\draw[blue,fill=blue] (-0.6363792902864169,0.42521504738362825) circle (4*0.200000000000000pt);
\draw[blue,fill=blue] (-1.0,1.0) circle (4*0.200000000000000pt);
\draw[blue,fill=blue] (-1.4142135623730951,0.0) circle (4*0.200000000000000pt);
\draw[blue,fill=blue] (-1.8122548927057758,0.36047991100347393) circle (4*0.200000000000000pt);
\draw[blue,fill=blue] (-2.613125929752753,-0.0) circle (4*0.200000000000000pt);
\draw[blue,fill=blue] (1.8122548927057758,-0.36047991100347393) circle (4*0.200000000000000pt);
\draw[blue,fill=blue] (1.082392200292394,-0.0) circle (4*0.200000000000000pt);
\draw[blue,fill=blue] (0.42521504738362814,0.6363792902864169) circle (4*0.200000000000000pt);
\draw[blue,fill=blue] (0.0,1.4142135623730951) circle (4*0.200000000000000pt);
\draw[blue,fill=blue] (0.6363792902864169,-0.42521504738362825) circle (4*0.200000000000000pt);
\draw[blue,fill=blue] (-0.42521504738362814,0.6363792902864169) circle (4*0.200000000000000pt);
\draw[blue,fill=blue] (-0.6363792902864169,-0.42521504738362825) circle (4*0.200000000000000pt);
\draw[blue,fill=blue] (-1.082392200292394,0.0) circle (4*0.200000000000000pt);
\draw[blue,fill=blue] (-1.8122548927057758,-0.36047991100347393) circle (4*0.200000000000000pt);
\draw[blue,fill=blue] (1.0,-1.0) circle (4*0.200000000000000pt);
\draw[blue,fill=blue] (0.42521504738362814,-0.6363792902864169) circle (4*0.200000000000000pt);
\draw[blue,fill=blue] (-0.0,-1.082392200292394) circle (4*0.200000000000000pt);
\draw[blue,fill=blue] (-0.42521504738362814,-0.6363792902864169) circle (4*0.200000000000000pt);
\draw[blue,fill=blue] (-1.0,-1.0) circle (4*0.200000000000000pt);
\draw[blue,fill=blue] (0.360479911003474,-1.8122548927057758) circle (4*0.200000000000000pt);
\draw[blue,fill=blue] (0.0,-1.4142135623730951) circle (4*0.200000000000000pt);
\draw[blue,fill=blue] (-0.360479911003474,-1.8122548927057758) circle (4*0.200000000000000pt);
\draw[blue,fill=blue] (0.0,-2.613125929752753) circle (4*0.200000000000000pt);
\end{tikzpicture}

%% file: figures/graph-D-4-4.tex
\tdplotsetmaincoords{8}{10}
\begin{tikzpicture}[tdplot_main_coords,anchorbase]
  \coordinate (o1) at (4/3,1/3,1/3);
  \coordinate (o2) at (1/3,4/3,1/3);
  \coordinate (o3) at (1/3,1/3,4/3);
  \coordinate (delta) at ($1/5*(o2)-1/5*(o1)$);
  \draw (0,0,0) -- (o1);
  \draw (0,0,0) -- (o2);
  \draw (0,0,0) -- (o3);
  \draw (o1) -- (o2);
  \draw (o1) -- (o3);
  \draw (o1) -- ($77/100*(o1)+(o2)$);
  \draw (o1) -- ($(o1)+(o3)$);
  \draw (o1) -- ($2*(o1)$);
  \draw (o2) -- (o3);
  \draw (o2) -- ($77/100*(o1)+(o2)$);
  \draw (o2) -- ($(o1)+(o3)$);
  \draw (o2) -- ($(o2)+(o3)$);
  \draw (o2) -- ($2*(o2)+35/1000*(o1)-60/100*(delta)$);
  \draw (o2) -- ($2*(o2)+1/2*(delta)$);
  \draw (o3) -- ($(o1)+(o3)$);
  \draw (o3) -- ($(o2)+(o3)$);
  \draw (o3) -- ($2*(o1)$);
  \draw ($2*(o1)$) -- ($77/100*(o1)+(o2)$);
  \draw[lightneon, thick] ($2*(o1)$) -- ($(o1)+(o3)$);
  \draw ($2*(o1)$) -- ($(o2)+(o3)$);
  \draw[lightneon, thick] ($77/100*(o1)+(o2)$) -- ($(o1)+(o3)$);
  \draw ($77/100*(o1)+(o2)$) -- ($2*(o2)+35/1000*(o1)-60/100*(delta)$);
  \draw ($77/100*(o1)+(o2)$) -- ($2*(o2)+1/2*(delta)$);
  \draw[tomato, very thick] ($77/100*(o1)+(o2)$) -- ($(o2)+(o3)$);
  \draw ($77/100*(o1)+(o2)$) -- ($(o1)+(o2)+(o3)-3/2*(delta)$);
  \draw ($77/100*(o1)+(o2)$) -- ($(o1)+(o2)+(o3)-1/2*(delta)$);
  \draw ($77/100*(o1)+(o2)$) -- ($(o1)+(o2)+(o3)+1/2*(delta)$);
  \draw ($77/100*(o1)+(o2)$) -- ($(o1)+(o2)+(o3)+3/2*(delta)$);
  \draw[lightneon, thick] ($(o1)+(o3)$) -- ($(o2)+(o3)$);
  \draw ($(o1)+(o3)$) -- ($(o1)+(o2)+(o3)-3/2*(delta)$);
  \draw ($(o1)+(o3)$) -- ($(o1)+(o2)+(o3)-1/2*(delta)$);
  \draw ($(o1)+(o3)$) -- ($(o1)+(o2)+(o3)+1/2*(delta)$);
  \draw ($(o1)+(o3)$) -- ($(o1)+(o2)+(o3)+3/2*(delta)$);
  \draw ($2*(o2)+1/2*(delta)$) -- ($(o2)+(o3)$);
  \draw ($2*(o2)+1/2*(delta)$) -- ($(o1)+(o2)+(o3)+1/2*(delta)$);
  \draw ($2*(o2)+1/2*(delta)$) -- ($(o1)+(o2)+(o3)+3/2*(delta)$);
  \draw ($2*(o2)+35/1000*(o1)-60/100*(delta)$) -- ($(o2)+(o3)$);
  \draw ($2*(o2)+35/1000*(o1)-60/100*(delta)$) -- ($(o1)+(o2)+(o3)-1/2*(delta)$);
  \draw ($2*(o2)+35/1000*(o1)-60/100*(delta)$) -- ($(o1)+(o2)+(o3)-3/2*(delta)$);
  \draw ($(o2)+(o3)$) -- ($(o1)+(o2)+(o3)-3/2*(delta)$);
  \draw ($(o2)+(o3)$) -- ($(o1)+(o2)+(o3)-1/2*(delta)$);
  \draw ($(o2)+(o3)$) -- ($(o1)+(o2)+(o3)+1/2*(delta)$);
  \draw ($(o2)+(o3)$) -- ($(o1)+(o2)+(o3)+3/2*(delta)$);
  \node[circle,inner sep=0.9*1.8pt,fill=zero] at (0,0,0){};
  \node[triangle,inner sep=0.9*1.6pt,fill=one] at (o1){};
  \node[square,inner sep=0.9*1.8pt,fill=two] at (o2){};
  \node[pentagon,inner sep=0.9*2.0pt,fill=three] at (o3){};
  \node[square,inner sep=0.9*1.8pt,fill=two] at ($2*(o1)$){};
  \node[pentagon,inner sep=0.9*2.0pt,fill=three] at ($77/100*(o1)+(o2)$){};
  \node[circle,inner sep=0.9*1.8pt,fill=zero] at ($(o1)+(o3)$){};
  \node[circle,inner sep=0.9*1.8pt,fill=zero] at ($2*(o2)+35/1000*(o1)-60/100*(delta)$){};
  \node[circle,inner sep=0.9*1.8pt,fill=zero] at ($2*(o2)+1/2*(delta)$){};
  \node[triangle,inner sep=0.9*1.6pt,fill=one] at ($(o2)+(o3)$){};
  \node[square,inner sep=0.9*1.8pt,fill=two] at ($(o1)+(o2)+(o3)-3/2*(delta)$){};
  \node[square,inner sep=0.9*1.8pt,fill=two] at ($(o1)+(o2)+(o3)-1/2*(delta)$){};
  \node[square,inner sep=0.9*1.8pt,fill=two] at ($(o1)+(o2)+(o3)+1/2*(delta)$){};
  \node[square,inner sep=0.9*1.8pt,fill=two] at ($(o1)+(o2)+(o3)+3/2*(delta)$){};
\end{tikzpicture}

%% file: figures/koorn-D-4-4.tex
\begin{tikzpicture}[anchorbase]
\def\a{1} \def\b{4} \def\scale{1}
\draw[cyan!30,very thin] ({-\b*\scale},{-\b*\scale}) grid ({\b*\scale},{\b*\scale});
\draw[->] ({-\b*\scale},0) -- ({\b*\scale},0);
\draw[->] (0,{-\b*\scale}) -- (0,{\b*\scale});
\draw (0,0) circle ({\b*\scale});
\draw[line width=2pt,red] plot[samples=100,domain=0:\a*360,smooth,variable=\t] ({\scale*((\b-\a)*cos(\t)+\a*cos((\b-\a)*\t/\a)},{\scale*((\b-\a)*sin(\t)-\a*sin((\b-\a)*\t/\a)});
\draw[blue,fill=blue] (-1.082392200292392,0.0) circle (4*0.200000000000000pt);
\draw[blue,fill=blue] (0.0,0.0) circle (4*1.20000000000000pt);
\draw[blue,fill=blue] (0.0,-2.613125929752753) circle (4*0.200000000000000pt);
\draw[blue,fill=blue] (0.0,-1.0823922002923938) circle (4*0.200000000000000pt);
\draw[blue,fill=blue] (2.6131259297527523,0.0) circle (4*0.200000000000000pt);
\draw[blue,fill=blue] (0.0,2.613125929752753) circle (4*0.200000000000000pt);
\draw[blue,fill=blue] (1.082392200292392,0.0) circle (4*0.200000000000000pt);
\draw[blue,fill=blue] (0.0,1.0823922002923938) circle (4*0.200000000000000pt);
\draw[blue,fill=blue] (-2.6131259297527523,0.0) circle (4*0.200000000000000pt);
\end{tikzpicture}

%% file: figures/graph-2Ac-4.tex
\begin{tikzpicture}[xscale=1.5,yscale=1.5,anchorbase]
  \coordinate (1) at (0,0);
  \coordinate (2) at (0,1);
  \coordinate (3) at (0,2);
  \coordinate (4) at (0,3);
  \coordinate (5) at (0,4);
  \coordinate (6) at (1,1);
  \coordinate (7) at (1,2);
  \coordinate (8) at (1,3);
  \coordinate (9) at (2,2);
  \coordinate (10) at (0.4,0.2);
  \coordinate (11) at (0.4,1.2);
  \coordinate (12) at (0.4,2.2);
  \coordinate (13) at (0.4,3.2);
  \coordinate (14) at (0.4,4.2);
  \coordinate (15) at (1.4,1.2);
  \coordinate (16) at (1.4,2.2);
  \coordinate (17) at (1.4,3.2);
  \coordinate (18) at (2.4,2.2);
  \draw (1) -- (2) -- (3) -- (4) -- (5);
  \draw (6) -- (7) -- (8);
  \draw (10) -- (11) -- (12) -- (13) -- (14);
  \draw (15) -- (16) -- (17);
  \draw (1) -- (11) -- (3) -- (13) -- (5);
  \draw (6) -- (16) -- (8);
  \draw (10) -- (2) -- (12) -- (4) -- (14);
  \draw (15) -- (7) -- (17);
  \draw (2) -- (6);
  \draw (11) -- (15);
  \draw (3) -- (7) -- (9);
  \draw (12) -- (16) -- (18);
  \draw (4) -- (8);
  \draw (13) -- (17);
  \draw (2) -- (15);
  \draw (11) -- (6);
  \draw (3) -- (16) -- (9);
  \draw (12) -- (7) -- (18);
  \draw (4) -- (17);
  \draw (13) -- (8);
  \draw (1) -- (6) -- (9) -- (8) -- (5) -- (14) -- (17) -- (18) -- (15) -- (10) -- (1);
  \draw (6) -- (3) -- (8) -- (17) -- (12) -- (15) -- (6);
  \draw (2) -- (7) -- (4);
  \draw (11) -- (16) -- (13);
  \draw[lightneon, thick] (2) -- (11);
  \draw[lightneon, thick] (3) -- (12);
  \draw[lightneon, thick] (4) -- (13);
  \draw[lightneon, thick] (7) -- (16);
  \node[circle,inner sep=1.1*1.8pt,fill=zero] at (1){};
  \node[circle,inner sep=1.1*1.8pt,fill=zero] at (3){};
  \node[circle,inner sep=1.1*1.8pt,fill=zero] at (5){};
  \node[circle,inner sep=1.1*1.8pt,fill=zero] at (15){};
  \node[circle,inner sep=1.1*1.8pt,fill=zero] at (17){};
  \node[circle,inner sep=1.1*1.8pt,fill=zero] at (9){};
  \node[triangle,inner sep=1.1*1.6pt,fill=one] at (2){};
  \node[triangle,inner sep=1.1*1.6pt,fill=one] at (4){};
  \node[triangle,inner sep=1.1*1.6pt,fill=one] at (16){};
  \node[square,inner sep=1.1*1.8pt,fill=two] at (10){};
  \node[square,inner sep=1.1*1.8pt,fill=two] at (12){};
  \node[square,inner sep=1.1*1.8pt,fill=two] at (14){};
  \node[square,inner sep=1.1*1.8pt,fill=two] at (6){};
  \node[square,inner sep=1.1*1.8pt,fill=two] at (8){};
  \node[square,inner sep=1.1*1.8pt,fill=two] at (18){};
  \node[pentagon,inner sep=1.1*2.0pt,fill=three] at (7){};
  \node[pentagon,inner sep=1.1*2.0pt,fill=three] at (11){};
  \node[pentagon,inner sep=1.1*2.0pt,fill=three] at (13){};
\end{tikzpicture}

%% file: figures/koorn-2Ac-4.tex
\begin{tikzpicture}[anchorbase]
\def\a{1} \def\b{4} \def\scale{1}
\draw[cyan!30,very thin] ({-\b*\scale},{-\b*\scale}) grid ({\b*\scale},{\b*\scale});
\draw[->] ({-\b*\scale},0) -- ({\b*\scale},0);
\draw[->] (0,{-\b*\scale}) -- (0,{\b*\scale});
\draw (0,0) circle ({\b*\scale});
\draw[line width=2pt,red] plot[samples=100,domain=0:\a*360,smooth,variable=\t] ({\scale*((\b-\a)*cos(\t)+\a*cos((\b-\a)*\t/\a)},{\scale*((\b-\a)*sin(\t)-\a*sin((\b-\a)*\t/\a)});
\draw[blue,fill=blue] (2.613125929752753,0.0) circle (4*0.200000000000000pt);
\draw[blue,fill=blue] (-1.082392200292394,0.0) circle (4*0.200000000000000pt);
\draw[blue,fill=blue] (-1.4142135623730951,0.0) circle (4*0.200000000000000pt);
\draw[blue,fill=blue] (-0.0,2.613125929752753) circle (4*0.200000000000000pt);
\draw[blue,fill=blue] (0.0,1.4142135623730951) circle (4*0.200000000000000pt);
\draw[blue,fill=blue] (0.0,1.082392200292394) circle (4*0.200000000000000pt);
\draw[blue,fill=blue] (0.0,0.0) circle (4*1.20000000000000pt);
\draw[blue,fill=blue] (1.4142135623730951,0.0) circle (4*0.200000000000000pt);
\draw[blue,fill=blue] (-0.0,-1.082392200292394) circle (4*0.200000000000000pt);
\draw[blue,fill=blue] (1.082392200292394,-0.0) circle (4*0.200000000000000pt);
\draw[blue,fill=blue] (0.0,-2.613125929752753) circle (4*0.200000000000000pt);
\draw[blue,fill=blue] (0.0,-1.4142135623730951) circle (4*0.200000000000000pt);
\draw[blue,fill=blue] (-2.613125929752753,-0.0) circle (4*0.200000000000000pt);
\end{tikzpicture}

%% file: figures/graph-2Ac-2-4.tex
\begin{tikzpicture}[xscale=2,yscale=2,anchorbase]
  \coordinate (1) at (0,0);
  \coordinate (2) at (0,1);
  \coordinate (3) at (0.1,2);
  \coordinate (3') at (-0.3,1.9);
  \coordinate (4) at (1,1);
  \coordinate (5) at (1.1,2);
  \coordinate (5') at (0.7,1.9);
  \coordinate (6) at (2.1,2);
  \coordinate (6') at (1.7,1.9);
  \coordinate (7) at (0.4,0.2);
  \coordinate (8) at (0.4,1.2);
  \coordinate (9) at (0.6,2.3);
  \coordinate (9') at (0.2,2.2);
  \coordinate (10) at (1.4,1.2);
  \coordinate (11) at (1.6,2.3);
  \coordinate (11') at (1.2,2.2);
  \coordinate (12) at (2.6,2.3);
  \coordinate (12') at (2.2,2.2);
  \draw (1) -- (2) -- (7) -- (8) -- (1);
  \draw (2) -- (3) -- (8) -- (9) -- (2);
  \draw (2) -- (3) -- (8) -- (9) -- (2);
  \draw (2) -- (3') -- (8) -- (9') -- (2);
  \draw (4) -- (5) -- (10) -- (11) -- (4);
  \draw (4) -- (5') -- (10) -- (11') -- (4);
  \draw (2) -- (4) -- (8) -- (10) -- (2);
  \draw (3') -- (5') -- (3) -- (5) -- (9') -- (11') -- (9) -- (11);
  \draw (11)  to [bend right=35] (3');
  \draw (6') -- (5') -- (6) -- (5) -- (12') -- (11') -- (12) -- (11) -- (6');
  \draw (1) -- (4) -- (10) -- (7) -- (1);
  \draw (6') -- (4) -- (12');
  \draw (6) -- (10) -- (12);
  \draw (3') -- (4) -- (9');
  \draw (3) -- (10) -- (9);
  \draw (5') -- (2) -- (11');
  \draw (5) -- (8) -- (11);
  \draw (3') -- (3) -- (9') -- (9) -- (3');
  \draw (5') -- (5) -- (11') -- (11) -- (5');
  \draw[lightneon, thick] (2) -- (8);
  \node[circle,inner sep=0.9*1.8pt,fill=zero] at (1){};
  \node[circle,inner sep=0.9*1.8pt,fill=zero] at (3'){};
  \node[circle,inner sep=0.9*1.8pt,fill=zero] at (6'){};
  \node[circle,inner sep=0.9*1.8pt,fill=zero] at (9'){};
  \node[circle,inner sep=0.9*1.8pt,fill=zero] at (10){};
  \node[circle,inner sep=0.9*1.8pt,fill=zero] at (12'){};  
  \node[triangle,inner sep=0.9*1.6pt,fill=one] at (2){};
  \node[triangle,inner sep=0.9*1.6pt,fill=one] at (5){};
  \node[triangle,inner sep=0.9*1.6pt,fill=one] at (11){};
  \node[square,inner sep=0.9*1.8pt,fill=two] at (3){};
  \node[square,inner sep=0.9*1.8pt,fill=two] at (4){};
  \node[square,inner sep=0.9*1.8pt,fill=two] at (6){};
  \node[square,inner sep=0.9*1.8pt,fill=two] at (7){};
  \node[square,inner sep=0.9*1.8pt,fill=two] at (9){};
  \node[square,inner sep=0.9*1.8pt,fill=two] at (12){};
  \node[pentagon,inner sep=0.9*2.0pt,fill=three] at (5'){};
  \node[pentagon,inner sep=0.9*2.0pt,fill=three] at (8){};
  \node[pentagon,inner sep=0.9*2.0pt,fill=three] at (11'){};
\end{tikzpicture}

%% file: figures/koorn-2Ac-2-4.tex
\begin{tikzpicture}[anchorbase]
\def\a{1} \def\b{4} \def\scale{1}
\draw[cyan!30,very thin] ({-\b*\scale},{-\b*\scale}) grid ({\b*\scale},{\b*\scale});
\draw[->] ({-\b*\scale},0) -- ({\b*\scale},0);
\draw[->] (0,{-\b*\scale}) -- (0,{\b*\scale});
\draw (0,0) circle ({\b*\scale});
\draw[line width=2pt,red] plot[samples=100,domain=0:\a*360,smooth,variable=\t] ({\scale*((\b-\a)*cos(\t)+\a*cos((\b-\a)*\t/\a)},{\scale*((\b-\a)*sin(\t)-\a*sin((\b-\a)*\t/\a)});
\draw[blue,fill=blue] (-2.613125929752753,-0.0) circle (4*0.200000000000000pt);
\draw[blue,fill=blue] (1.082392200292394,-0.0) circle (4*0.200000000000000pt);
\draw[blue,fill=blue] (-1.0,1.0) circle (4*0.200000000000000pt);
\draw[blue,fill=blue] (-0.0,2.613125929752753) circle (4*0.200000000000000pt);
\draw[blue,fill=blue] (0.0,1.082392200292394) circle (4*0.200000000000000pt);
\draw[blue,fill=blue] (-1.0,-1.0) circle (4*0.200000000000000pt);
\draw[blue,fill=blue] (0.0,0.0) circle (4*1.20000000000000pt);
\draw[blue,fill=blue] (-0.0,-1.082392200292394) circle (4*0.200000000000000pt);
\draw[blue,fill=blue] (1.0,-1.0) circle (4*0.200000000000000pt);
\draw[blue,fill=blue] (-1.082392200292394,0.0) circle (4*0.200000000000000pt);
\draw[blue,fill=blue] (0.0,-2.613125929752753) circle (4*0.200000000000000pt);
\draw[blue,fill=blue] (1.0,1.0) circle (4*0.200000000000000pt);
\draw[blue,fill=blue] (2.613125929752753,0.0) circle (4*0.200000000000000pt);
\end{tikzpicture}

%% file: figures/graph-E-4.tex
\begin{tikzpicture}[anchorbase]
  \coordinate (1) at (4,1);
  \coordinate (2) at (4,-1);
  \coordinate (3) at (3,1);
  \coordinate (4) at (3,-1);
  \coordinate (5) at (2,1);
  \coordinate (6) at (2.3,0);
  \coordinate (7) at (2,-1);
  \coordinate (8) at (1.5,1.4);
  \coordinate (9) at (1.5,-1.4);
  \coordinate (10) at (.5,1.4);
  \coordinate (11) at (.5,-1.4);
  \coordinate (12) at (0,0);
  \draw (1) -- (3);
  \draw (1) -- (4);
  \draw (1) -- (6);
  \draw (2) -- (3);
  \draw (2) -- (4);
  \draw (2) -- (6);
  \draw[lightneon, thick] (3) -- (4);
  \draw (3) -- (5);
  \draw[lightneon, thick] (3) -- (6);
  \draw (3) -- (7);
  \draw (3) -- (8);
  \draw (3) -- (10);
  \draw (4) -- (5);
  \draw[lightneon, thick] (4) -- (6);
  \draw (4) -- (7);
  \draw (4) -- (9);
  \draw (4) -- (11);
  \draw[lightneon, thick] (5) -- (6);
  \draw (5) -- (8);
  \draw (5) -- (9);
  \draw (5) -- (10);
  \draw (5) -- (11);
  \draw (5) -- (12);
  \draw[lightneon, thick] (6) -- (7);
  \draw (6) -- (8);
  \draw (6) -- (9);
  \draw (6) -- (10);
  \draw (6) -- (11);
  \draw (7) -- (8);
  \draw (7) -- (9);
  \draw (7) -- (10);
  \draw (7) -- (11);
  \draw (7) -- (12);
  \draw (8) -- (9);
  \draw (8) -- (11);
  \draw (8) -- (12);
  \draw (9) -- (10);
  \draw (9) -- (12);
  \draw (10) -- (11);
  \draw (10) -- (12);
  \draw (11) -- (12); 
  \node[circle,inner sep=0.9*1.8pt,fill=zero] at (1){};
  \node[circle,inner sep=0.9*1.8pt,fill=zero] at (2){};
  \node[triangle,inner sep=0.9*1.6pt,fill=one] at (3){};
  \node[pentagon,inner sep=0.9*2.0pt,fill=three] at (4){};
  \node[circle,inner sep=0.9*1.8pt,fill=zero] at (5){};
  \node[square,inner sep=0.9*1.8pt,fill=two] at (6){};
  \node[circle,inner sep=0.9*1.8pt,fill=zero] at (7){};
  \node[pentagon,inner sep=0.9*2.0pt,fill=three] at (8){};
  \node[triangle,inner sep=0.9*1.6pt,fill=one] at (9){};
  \node[pentagon,inner sep=0.9*2.0pt,fill=three] at (10){};
  \node[triangle,inner sep=0.9*1.6pt,fill=one] at (11){};
  \node[square,inner sep=0.9*1.8pt,fill=two] at (12){};
\end{tikzpicture}

%% file: figures/koorn-E-4.tex
\begin{tikzpicture}[anchorbase]
\def\a{1} \def\b{4} \def\scale{1}
\draw[cyan!30,very thin] ({-\b*\scale},{-\b*\scale}) grid ({\b*\scale},{\b*\scale});
\draw[->] ({-\b*\scale},0) -- ({\b*\scale},0);
\draw[->] (0,{-\b*\scale}) -- (0,{\b*\scale});
\draw (0,0) circle ({\b*\scale});
\draw[line width=2pt,red] plot[samples=100,domain=0:\a*360,smooth,variable=\t] ({\scale*((\b-\a)*cos(\t)+\a*cos((\b-\a)*\t/\a)},{\scale*((\b-\a)*sin(\t)-\a*sin((\b-\a)*\t/\a)});
\draw[blue,fill=blue] (2.613125929752753,0.0) circle (4*0.200000000000000pt);
\draw[blue,fill=blue] (0.0,2.613125929752753) circle (4*0.200000000000000pt);
\draw[blue,fill=blue] (0.0,1.082392200292394) circle (4*0.200000000000000pt);
\draw[blue,fill=blue] (0.0,0.0) circle (4*0.800000000000000pt);
\draw[blue,fill=blue] (-2.613125929752753,-0.0) circle (4*0.200000000000000pt);
\draw[blue,fill=blue] (1.082392200292394,-0.0) circle (4*0.200000000000000pt);
\draw[blue,fill=blue] (-1.082392200292394,0.0) circle (4*0.200000000000000pt);
\draw[blue,fill=blue] (-0.0,-1.082392200292394) circle (4*0.200000000000000pt);
\draw[blue,fill=blue] (0.0,-2.613125929752753) circle (4*0.200000000000000pt);
\end{tikzpicture}